\documentclass{article}

\usepackage{arxiv}

\usepackage[utf8]{inputenc} 
\usepackage[T1]{fontenc}    
\usepackage{hyperref}       
\usepackage{url}            
\usepackage{booktabs}       
\usepackage{amsfonts}       
\usepackage{nicefrac}       
\usepackage{microtype}      
\usepackage{lipsum}
\usepackage{graphicx}
\graphicspath{ {./images/} }

\usepackage{bbm}
\usepackage{enumitem}

\usepackage{crop}
\usepackage{caption}
\usepackage{amsmath}
\usepackage{mathrsfs}
\usepackage{array}
\usepackage{color}
\usepackage{xcolor}
\usepackage{amssymb}
\usepackage{flushend}
\usepackage{stfloats}
\usepackage[figuresright]{rotating}
\usepackage{chngpage}
\usepackage{totcount}
\usepackage{fix-cm}
\usepackage{amsthm}
\usepackage{dirtytalk}

\title{Continuous Time Locally Stationary\\ Wavelet Processes}

\newcommand{\toptitlebar}{
  \hrule height 2pt
  \vskip 0.25in
  \vskip -\parskip%
}
\newcommand{\bottomtitlebar}{
  \vskip 0.29in
  \vskip -\parskip
  \hrule height 2pt
  \vskip 0.09in%
}

\newcommand\numberthis{\addtocounter{equation}{1}\tag{\theequation}}

\allowdisplaybreaks

\newtheorem{theorem}{Theorem}
\newtheorem{proposition}{Proposition}
\newtheorem{remark}{Remark}
\theoremstyle{definition}
\newtheorem{definition}{Definition}
\newtheorem{example}{Example}

\author{
 Henry Antonio Palasciano \\
  Department of Mathematics\\
 Imperial College London\\
  SW7 2AZ,  London, UK \\
  \texttt{henry.palasciano17@imperial.ac.uk}
   \And
    Marina I. Knight \\
    Department of Mathematics \\
    University of York \\
    YO10 5DD, York, UK \\
    \texttt{marina.knight@york.ac.uk}\\
    \And
        Guy P. Nason \\
      Department of Mathematics\\
 Imperial College London\\
  SW7 2AZ,  London, UK \\
  \texttt{g.nason@imperial.ac.uk}\\
}

\begin{document}

\maketitle

\begin{abstract}
This article introduces the class of continuous time locally stationary wavelet processes. Continuous time models enable us to properly provide scale-based time series models for irregularly-spaced observations for the first time, while also permitting a spectral representation of the process over a continuous range of scales. We derive results for both the theoretical setting, where we assume access to the entire process sample path, and a more practical one, which develops methods for estimating the quantities of interest from sampled time series. The latter estimates are  accurately computable in reasonable time by solving the relevant linear integral equation using the iterative soft-thresholding algorithm due to Daubechies, Defrise and De~Mol. Appropriate smoothing techniques are also developed and applied in this new setting. Comparisons to previous methods are conducted on the heart rate time series of a sleeping infant. Additionally, we exemplify our new methods by computing spectral and autocovariance estimates on irregularly-spaced heart rate data obtained from a recent sleep-state study.
\end{abstract}

\keywords{
Continuous time models; Irregular spacing; Local autocovariance;  Non-stationary time series; Wavelets; Wavelet spectrum.
}

\section{Introduction}
It is often customary to assume that a time series is second-order stationary,
since plenty of theoretically sound and well-tested methods are readily available for its analysis, see, for example, \cite{brd} or \cite{wald}.
It is well known that a stationary stochastic process $X(t)$, $t\in\mathbb R$, has a spectral representation of the form
\begin{equation}\label{eq:spec_rep}
    X(t) = \int_{-\infty}^\infty \exp(i\omega t) dZ(\omega),
\end{equation}
where $Z(\omega)$ is a square-integrable orthogonal increment process with $\mathbb E \left[|dZ(\omega)|^2\right] = dS^{(I)}(\omega)$ and $S^{(I)}(\omega)$ is the integrated spectrum of $X(t)$. The autocovariance function of $X(t)$, denoted $c(\tau)$,  has a similar representation in the frequency domain,
\begin{equation} \label{eq:cov_spec_fourier}
    c(\tau) = \int_{-\infty}^\infty \exp(i\omega \tau) dS^{(I)}(\omega),
\end{equation}
see, for example, \cite{priestley}. However, in practice, many time series are unlikely to be stationary, such as is often the case in finance for instance \cite{finance, appl_fin}. To extend these quantities to the non-stationary case, \cite{priestley} suggested introducing a time-varying amplitude function $A(t,\omega)$ into \eqref{eq:spec_rep} giving 
\begin{equation} \label{eq:spectr}
    X(t) = \int_{-\infty}^\infty A(t,\omega)\exp(i\omega t) dZ(\omega).
\end{equation}
Models such as \eqref{eq:spectr} form the basis for
locally stationary time series,
later developed by \cite{dahlhaus1997} in discrete time.

\cite{LSWP}
introduced the locally stationary wavelet processes,
which replaced the set of functions $\{\exp(i\omega t)\}$ by a set of discrete non-decimated wavelets. As with locally stationary Fourier processes \cite{dahlhaus1997}, the local stationarity
is enforced via constraints on the amplitude functions.
From locally stationary wavelet processes, one can then obtain a time-scale decomposition of the spectral density function, known as the evolutionary wavelet spectrum, and a wavelet representation of the local autocovariance similar to \eqref{eq:cov_spec_fourier}.

However, discrete time locally stationary wavelet processes are constructed under the assumption that the time series data are obtained at a constant sampling rate, which is not always the case as, in practice, observations can either be missing or be irregularly-spaced. Attempts to provide a suitable alternative for missing data were made by~\cite{missing2} via the use of the non-decimated lifting transform from~\cite{missing1}. Unfortunately: 1.\ the~\cite{missing2} estimates only
resemble the true spectrum defined in~\cite{LSWP} in that they can only estimate an underlying spectrum at the
$n$ time points in the irregular sample and nowhere in-between;
2.\ the lifting methods of \cite{missing2} bestow no complete analytical formula connection, nor asymptotic
theory, between the
expected periodogram and the original spectrum, at each time, such as can be found in equation~(20) in Proposition~4 of~\cite{LSWP}
or Theorem~\ref{thm:expectation2} below;
3.\ the scales in \cite{missing2} are artificial and depend on the inter-point spacing
of the samples and, hence, do not have a one-to-one correspondence to the actual wavelet scales underlying 
both the discrete and continuous
locally stationary wavelet processes
in \cite{LSWP} and below, respectively. Ideally, spectral estimates should not depend on the
sampling scheme.
4.\
\cite{missing2} state that their methods are \say{highly computationally intensive} in that  the periodogram evaluated at each time point requires computation of a full lifting
(second generation wavelet) transform for multiple, $P$, runs corresponding to different
lifting trajectories. For good estimation, typically $P$ is large
and was set to $P=1000$ and $P=5000$ in~\cite{missing2}.
Hence, periodogram estimation in~\cite{missing2} requires $n P$ lifting wavelet transforms compared to our method which uses one, or certainly at least $O(1)$. \cite{missing2} do
not specify spectral computational effort,
but they
remark that it is considerably more
computationally intensive than their periodogram
computation.

Instead, we propose the use of genuine continuous time locally stationary wavelet processes, which are natural counterparts to the discrete time versions from~\cite{LSWP}. To our knowledge, our new methods are the only ones that can accurately estimate wavelet spectral properties of an irregularly-sampled process within a reasonable computational time frame. Our method estimates the spectrum itself, not a simulacrum as in~\cite{missing2},
and is among the few guaranteeing the non-negativity of the spectrum (the other we know of was for regularly spaced data
in~\cite{fryzlewicz_phd} Chapter~4).

Aside from the practical advantages that a continuous time model provides, our results are also interesting from a more theoretical perspective. In fact, many discrete time models and methods also possess a continuous time counterpart, see, for example,~\cite{brock2, brock1, ARMA, edo, lorenzo}. Recently,~\cite{bennet} adapted locally stationary Fourier processes from the discrete~\cite{dahlhaus1997, dahlhaus2012} to the continuous time setting. Since, so far, the theory of locally stationary wavelet processes exists in discrete time only, our objective is to introduce  parallel
continuous time models, study their theoretical properties and develop
a practical means for their estimation and analysis.
In doing so, we obtain a representation of the evolutionary wavelet spectrum across a continuous range of scales, providing a complete decomposition of the underlying frequencies of the process.

\section{Background and Notation} \label{sec:background}
We briefly review the relevant background material for our continuous time locally stationary wavelet processes. Supplementary Material Section~S1 contains a more detailed discussion.

Discrete time locally stationary wavelet processes~\cite{LSWP, nasonbook} are time series models that have found applications in various domains, including finance~\cite{appl_fin}, economics~\cite{appl_econ}, aliasing detection~\cite{aliastests, aliastests2D}, biology~\cite{wav_bio}, gravitational-wave detection~\cite{gravity} and energy~\cite{appl_energy1, appl_energy2}, to name but a few.
They are defined as a sequence of doubly indexed stochastic processes having the representation
\begin{equation*}
        X_{t,T} = \sum_{j=1}^\infty\sum_{k=-\infty}^\infty w_{j,k;T} \psi_{j,k-t}\xi_{j,k}
\end{equation*}
in the mean square sense, where $\{w_{j,k;T}\}_{j\in\mathbb N, k\in\mathbb Z}$ is a set of amplitudes, $\{\psi_{j,k-t}\}_{j\in\mathbb N, k\in\mathbb Z}$ are discrete wavelets and $\{\xi_{j,k}\}_{j\in\mathbb N, k\in\mathbb Z}$ is a collection of uncorrelated random variables with zero mean and unit variance. Conditions are imposed on the amplitude functions to prevent them from varying too quickly, thus preserving the local stationarity of the process. More importantly, one assumes that, as $T\rightarrow\infty$,  $w_{j,k;T}$ converges to $W_j(k/T)$, a rescaled-time Lipschitz continuous amplitude function defined on
$z = k/T \in (0,1)$. The idea of rescaled-time was first introduced by~\cite{dahlhaus1997} and is frequently used when working with locally stationary processes. The corresponding evolutionary wavelet spectrum, autocorrelation wavelet and local autocovariance of the locally stationary wavelet process, for $j\in\mathbb N$ and $\tau\in\mathbb Z$, are
\begin{equation*}
        S_j(z) = |W_j(z)|^2,\quad \Psi_j(\tau) =  \sum_k \psi_{j,k}\psi_{j,k-\tau}, \quad\text{and}\quad c(z,\tau) = \sum_{j=1}^\infty S_j(z)\Psi_j(\tau),
    \end{equation*}
respectively. To estimate the evolutionary wavelet spectrum,~\cite{LSWP} first estimate the raw wavelet periodogram using the non-decimated wavelet transform and then apply a correction with the inner product operator $A_{j,\ell} = \sum_{\tau}\Psi_j(\tau)\Psi_\ell(\tau)$, $j,\ell\in\mathbb{N}$. See~\cite{LSWP} for further details. Section~\ref{sec:lswtheory} defines the continuous time counterparts of all these quantities and Section~\ref{sec:practical} develops methods to estimate these in practice.

L\'{e}vy processes and bases play a crucial role in the definition of locally stationary wavelet processes in continuous time. The following overview draws heavily from \cite{walsh, spectrallevy, applebaum, almut}. We take $(\Omega, \mathcal{F}, \mathbb P)$ to be the underlying probability space, $\{S,\mathcal{B}(S)\}$ to be a Borel space and denote the bounded Borel sets of $\mathcal{B}(S)$ by $\mathcal{B}_b(S)$. A real valued L\'{e}vy basis $L$ on $S$ is a collection $\{L(A): A\in \mathcal{B}_b(S)\}$ of random variables such that the law of $L(A)$ is infinitely divisible for all $A\in \mathcal{B}_b(S)$ and, if $A$ and $B$ are disjoint subsets in $\mathcal{B}_b(S)$, $L(A)$ and $L(B)$ are independent and $L(A \cup B) = L(A) + L(B)$ almost surely. We restrict our attention to homogeneous L\'{e}vy bases as defined by~\cite{almut}, since this is synonymous with stationarity. We say that a homogeneous L\'{e}vy basis $G$ is Gaussian if $G(A)\sim\mathcal N \{ \mu \text{Leb}(A), \sigma^2\text{Leb}(A) \}$ for all $A\in\mathcal{B}(\mathbb R^n)$ and some $\mu\in\mathbb R$ and $\sigma\in\mathbb R^+$, where Leb denotes the Lebesgue measure. Homogeneous Gaussian L\'{e}vy bases are equivalent to the concept of white noise defined by~\cite{walsh}. Furthermore, if $L$ is a homogeneous L\'{e}vy basis on $\{\mathbb R, \mathcal{B}(\mathbb R)\}$, then the process $(L_t)_{t\geq 0}$ with $L_t = L([0,t])$ is a L\'{e}vy process. Conversely, given a L\'{e}vy process $(L_t)_{t\geq 0}$, we can obtain a L\'{e}vy basis $L$ by setting $L((s,t])=L_t-L_s$.

The iterative soft-thresholding algorithm of~\cite{daub} is a powerful regularization method which we find useful in Section~\ref{sec:practical} for estimating the continuous time evolutionary wavelet spectrum from an estimate of the raw wavelet periodogram. This algorithm was developed to obtain regularized solutions to linear inverse problems $Th=f$, where $f = g+\varepsilon$ is some noisy observation of the true signal $g$. In our setting, $T$ is a linear integral operator, $f,g$ and $h$ are some functions in a Hilbert space $\mathcal{H}$ and $\varepsilon$ is the noise. In a purely theoretical setting, Mercer's theorem~\cite{mercer} or methods from \cite{kress} can be used to obtain a solution. Instead, \cite{daub} suggest minimizing the functional \begin{equation}\label{eq:daub}
    \Phi_{\mu, p}(h) = ||Th-f||^2 + \sum_{\gamma\in\Gamma}\mu_\gamma |\langle h, e_\gamma\rangle|^p,
\end{equation}
for $1\leq p\leq 2$, where $\{e_\gamma\}_{\gamma\in\Gamma}$ is an orthogonal basis of $\mathcal{H}$ and  $\mu=(\mu_\gamma)_{\gamma\in\Gamma}$ is a vector of weights, which act as regularization parameters. In our setting, we are interested in penalizing positive and negative coefficients asymmetrically for the case $p=1$ and so replace each $\mu_\gamma\in\mu$ with the positive weights $(\mu_\gamma^+, \mu_\gamma^-)$. The solution is then given by the iterative soft-thresholding algorithm
\begin{equation} \label{eq:iter_scheme}
    h^{n} = \mathcal{S}_\mu\left\{ h^{n-1} + T^\dagger\left(f-Th^{n-1}\right)\right\},
\end{equation}
where
\begin{equation} \label{eq:shrink}
    \mathcal{S}_{\mu} (h) = \sum_{\gamma \in\Gamma} \mathcal{S}_{\mu_\gamma^+,\mu_\gamma^-}\left(\langle h, e_\gamma\rangle \right) e_\gamma,
    \quad
    \mathcal{S}_{\mu_\gamma^+, \mu_\gamma^-}(x)=\left\{\begin{array}{lll}
    x-\mu_\gamma^+/2, & & \text{if } x\geq\mu_\gamma^+/2, \\
    0, & & \text{if } -\mu^-_\gamma/2<x<\mu^+_\gamma/2,\\
    x+\mu^-_\gamma/2, & & \text{if } x\leq-\mu^-_\gamma/2
    \end{array}\right.
\end{equation}
and $T^\dag$ denotes the adjoint of $T$. As shown in~\cite{daub_inv}, this iterative scheme converges to the unique minimizer of~\eqref{eq:daub} for arbitrarily chosen $h^0\in \mathcal{H}$. 
More details can be found in \cite{daub_inv, daub_inv3} or~\cite{inv_book}.

\section{Locally Stationary Wavelet Processes in Continuous Time: Theory}
\label{sec:lswtheory}

\subsection{Process Definition}

Our continuous time locally stationary wavelet processes rely on  continuous wavelets $\psi \in L^2(\mathbb R)$ for which we assume that 
\begin{equation} \label{eq:wav_asump}
    \int_{-\infty}^\infty \psi(t) dt=0, \quad ||\psi||^2 = \int_{-\infty}^\infty |\psi(t)|^2 dt=1.
\end{equation}
From $\psi\in L^2(\mathbb R)$ an entire family of wavelets $\{\psi_{u,v}(t), u\in\mathbb R^+, v\in\mathbb R\}$ can be generated via
\begin{equation} \label{eq:wav_fam}
    \psi_{u,v}(t) = u^{-1/2}\psi\left\{u^{-1}(t-v) \right\} = \psi(u, t-v)
\end{equation}
preserving the properties in~\eqref{eq:wav_asump}, for all $u\in\mathbb R^+$ and $ v\in\mathbb R$. In~\eqref{eq:wav_fam} we have introduced the notation $\psi(\,\cdot\,,\,\cdot\,)$, making the dependence on scale more explicit, with $\psi(1, t)=\psi(t)$. For more information on the continuous wavelet transform, see \cite{daub}.

In continuous time, we assume knowledge of the entire path rather than a finite set of observations. Consequently, it becomes
irrelevant to define a sequence of stochastic processes corresponding to increasing amounts of data becoming available.
Hence, continuous-time locally stationary wavelet process are defined as follows.

\begin{definition} \label{def:clswp}
A continuous time locally stationary wavelet process, $\{X(t)\}_{t\in\mathbb R}$ is a stochastic process defined in the mean square sense as 
\begin{equation}\label{eq:lsw}
X(t) = \int_0^\infty\int_{-\infty}^\infty W(u, v)\psi(u, v-t)L(du, dv),
\end{equation}
where $L$ is a homogeneous square-integrable L\'{e}vy basis with $E\{L(A)\}=0$ and $E\{L(A)^2\}=\text{Leb}(A)$ for all $A\in \mathcal B_b(\mathbb R^+\times\mathbb R)$,  $\{\psi_{u,v}(x), u\in\mathbb R^+, v\in\mathbb R\}$
is a continuous wavelet family and $W\in L^2(\mathbb R ^+, \mathbb R)$ is an amplitude function.
Further, the following conditions on $W$ hold: (i) for each fixed scale $u\in\mathbb R^+$, $W(u,\,\cdot\,)$ is a Lipschitz continuous function with Lipschitz constant $K(u)$. The Lipschitz constant function $K(\,\cdot\,)$ is bounded and satisfies 
\begin{equation}\label{eq:lipschitz}
\int_0^\infty uK(u) du <\infty;
\end{equation}
(ii) for each fixed $v\in \mathbb R$, $W(\,\cdot\,,v)$ satisfies
\begin{equation*}
    \int_{0}^\infty |W(u,v)|^2 du<\infty.
\end{equation*}
\end{definition}

Locally stationary time series are stochastic processes whose statistical properties change gradually over time. The use of wavelets enables an implicit local representation of the process, while the smoothness assumptions imposed on $W$ ensure that locally the process behaves in an approximately stationary fashion. The condition on the Lipschitz constants in  \eqref{eq:lipschitz} guarantees that the rate at which the statistical properties are allowed to vary decreases as the scales widen. This ensures that the process behaves in a realistic fashion and is required for estimation.

For representation \eqref{eq:lsw} to converge in mean square, the amplitude functions cannot be allowed to grow arbitrarily with $u$. Hence, assumption (ii) ensures that $|W(u,v)|^2 \rightarrow 0$ as $u\rightarrow\infty$.

The requirement that the L\'{e}vy basis is homogeneous is essential, as this implies that the basis is stationary, and vice-versa, see Supplementary Material Section~S1 or~\cite{almut}. This ensures that the non-stationarity is controlled entirely by the amplitude functions $W$. Furthermore, since $L$ is by definition independent when acting on disjoint sets, any correlation structure in the process is captured entirely by the wavelets themselves.

\begin{remark}
    The amplitude function $W(u, v)$ is defined for all $v\in\mathbb R$ as opposed to its discrete counterpart,
    which is only defined on $(0,1)$.
    Later, in Section~\ref{sec:practical},
    when adapting this definition for practical applications, we will again return to the concept of rescaled-time.
\end{remark}

\begin{example}
A Haar moving average process of order $\alpha$, or Haar \textsc{ma}($\alpha$), is a continuous time locally stationary wavelet process
$X(t)$ with amplitude function $W(u,v) = 1$ for $u = \alpha \in\mathbb R^+$ and zero otherwise. The L\'{e}vy basis is $L(du,dv) = \delta_\alpha(du)G(dv)$, where $\delta_\alpha$ is the Dirac measure centred at $\alpha$ and $G$ is a homogeneous Gaussian L\'{e}vy basis on $\mathcal B_b(\mathbb R)$. The underlying wavelet family $\{\psi_{u,v} \}$ is derived from the Haar wavelet, see Supplementary Material S2, through $\psi(x) = \psi_H(-x)$, with the reflection about the y-axis chosen to maintain a backward-looking process. The Haar \textsc{ma}($\alpha$) can be written as
\begin{equation*}
    X(t) = \int_0^\infty\int_{-\infty}^{\infty} \psi_H(u,v-t)\delta_\alpha(du)G(dv) = \alpha^{-1/2}\left\{\left(B_t - B_{t-\alpha/2} \right) - \left(B_{t-\alpha/2} - B_{t-\alpha}\right)\right\},
\end{equation*}
where $B_t$ is a standard Brownian motion.
Essentially, Haar moving average processes capture the difference between consecutive Brownian increments, each of width $\alpha/2$. See~\cite{LSWP} for the discrete-time analogue. There, for example the Haar \textsc{ma}(1) and \textsc{ma}(2) processes are
\begin{equation*}
    X_t^1 = 2^{-1/2} (\varepsilon_t - \varepsilon_{t-1}), \quad X_t^2 = 2^{-1} (\varepsilon_t + \varepsilon_{t-1} - \varepsilon_{t-2} - \varepsilon_{t-3}),
\end{equation*}
respectively, where $\{\varepsilon_t\}_{t\in\mathbb Z}$ is a sequence of independent and identically distributed random variables with zero mean and unit variance.
\end{example}

\subsection{Continuous Time Evolutionary Wavelet Spectrum}

The evolutionary wavelet spectrum
gives a location-scale decomposition of the spectral power of a continuous time
locally stationary wavelet process, and is defined as follows.
\begin{definition}
    The continuous time evolutionary wavelet spectrum
    $S$
    of a locally stationary wavelet process $X(t)$, with respect to the wavelet family $\psi$, is defined by 
    \begin{equation*}
        S(u, v) = |W(u,v)|^2,
    \end{equation*}
    for all scales $u\in\mathbb R^+$ and
    times $v\in\mathbb R$. 
\end{definition}
Assumption (ii) in Definition~\ref{def:clswp} implies that, for each fixed $v\in\mathbb R$,
\begin{equation*}
        \int_0^\infty S(u,v) du<\infty.
\end{equation*}
\begin{remark}
More generally and technically, we can define 
\begin{equation*}
    S^{(I)}(du, dv) = E \left\{|W(u,v)L(du,dv)|^2\right\},
\end{equation*}
where $S^{(I)}$ is the continuous time integrated evolutionary wavelet spectrum. We can write this as 
\begin{equation*}
    S^{(I)}(du, dv) = S(u,v)  du dv,
\end{equation*}
since $L$ is a homogeneous L\'{e}vy basis. We can then define 
\begin{equation}\label{eq:power_spec}
    S^{(I)}_v(u) = \int_0^u S(x, v)  dx = \frac{\partial S^{(I)}}{\partial v} (u, v),
\end{equation}
representing the contribution to the total power of the process in a small neighbourhood of $v$ for scales less than or equal to $u$, mirroring the evolutionary power spectrum of \cite{priestley}. 
\end{remark}

\setcounter{example}{0}
\begin{example}
    \textit{continued.} The evolutionary wavelet spectrum of a Haar \textsc{ma}($\alpha$) process is
\begin{equation*}
    S(u,v) = \left\{ \begin{array}{ccc}
     1,    &&  \text{for } u=\alpha,\\
     0,    &&  \text{for } u \neq \alpha.
    \end{array}
    \right.
\end{equation*}
\end{example}

\subsection{Continuous Time Local Autocovariance}
\cite{LSWP} show that the autocovariance of a locally stationary wavelet process has a wavelet representation, known as the local autocovariance. We investigate continuous time versions and begin with continuous autocorrelation wavelets,
see~\cite{autocorr_wavelets}.
\begin{definition}
    For scale $u\in\mathbb R^+$ and lag $\tau\in\mathbb R$, the continuous autocorrelation wavelets are
    \begin{equation*}
        \Psi(u, \tau) = \int_{-\infty}^{\infty}\psi(u,v)\psi(u,v-\tau)dv.
    \end{equation*}
\end{definition}

\begin{definition} 
    For time $t \in\mathbb R$
    and lag $\tau\in\mathbb R$, the continuous time local autocovariance is
    \begin{equation}\label{eq:autocov_form}
        c(t,\tau) = \int_0^\infty S(u,t)\Psi(u,\tau)du,
    \end{equation}
    where $S(u,t)$ and $\Psi(u,\tau)$, $u\in\mathbb R^+$, are the
    continuous time evolutionary wavelet spectrum and autocorrelation wavelets, respectively.
\end{definition}

The usual autocovariance function of a
continuous time locally stationary wavelet processes is
\begin{equation}
    \label{eq:regularacf}
    c_X(t,\tau) = \text{cov}\{X(t),X(t+\tau)\} = E\{X(t)X(t+\tau)\},
\end{equation}
since $X(t)$ has zero mean. This can also be written in terms of the evolutionary wavelet spectrum and the underlying wavelets of the process as
\begin{equation}
    \label{eq:autocov_standard}
        c_X(t, \tau) = \int_0^\infty\int_{-\infty}^\infty S(u,t+v)\psi(u,v)\psi(u,v-\tau) dv  du,
\end{equation}
see the proof of Proposition~\ref{prop:autocov} in Supplementary Material Section~S6. Comparing~\eqref{eq:autocov_form} and~\eqref{eq:autocov_standard}, we see that, at each time $t$, the local autocovariance effectively replaces $S(u, t+v)$ with the constant $S(u,t)$, as a function of $v$, over a finite interval centred around $t$ of width matching that of the wavelet. Given the slow evolution of the statistical properties of the process, $S(u, t+v)$ is similar to $S(u, t)$ in this region and, if the process is stationary, the two expressions are equivalent. Therefore the local autocovariance assumes that the spectrum is locally constant, or equivalently that the process is stationary at a local level, a concept which underpins the framework of locally stationary processes. This is similar to the discrete time local autocovariance defined in~\cite{LSWP}. The next proposition establishes a bound on the maximum difference between expressions~\eqref{eq:autocov_form} and~\eqref{eq:autocov_standard}.

\begin{proposition}\label{prop:autocov}
Let $\{X(t)\}_{t\in\mathbb R}$ be a continuous time locally stationary wavelet process, $c_X(t,\tau)$  be its autocovariance function and $c(t,\tau)$ its local autocovariance. Then,
\begin{equation*}
    |c_X(t, \tau) - c(t,\tau)| \leq \gamma \int_0^\infty u K(u)  du,
\end{equation*}
for all scales $u\in\mathbb R^+$ and times $v\in\mathbb R$, where $\gamma$ is a constant which depends on the decay structure of the wavelet $\psi$ and $K(\,\cdot\,)$ is the Lipschitz constant function from Definition~\ref{def:clswp}.
\end{proposition}

In essence, the more stationary the process $X(t)$ is, the closer $c_X(t, \tau)$ and $c(t, \tau)$ will be, since $ \int_0^\infty u K(u)  du$ can be interpreted as a measure of the non-stationarity of the process. For a more local bound, if the spectrum $S$ admits a Taylor series representation, we can write
\begin{equation*}
    S(u,t+v) =  S(u, t) + v \frac{\partial S}{\partial v}(u, t) + o(v^2),
\end{equation*}
for $v$ close to $t$. Using this expansion we can write 
\begin{equation*}
    |c_X(t,\tau) - c(t,\tau)| \leq \gamma\int_0^\infty u \left|\frac{\partial S}{\partial v}(u, t) \right| du + M,
\end{equation*}
where $M \geq 0$ is a constant which bounds the discrepancy due to the higher order terms. From the above expression it becomes clearer that the difference between the $c_X(t, \tau)$ and $c(t, \tau)$ depends on the local structure of the process.  A comparison of the local and standard autocovariance functions on a simple example can be found in Supplementary Material Section~S3.

\begin{definition}
    The continuous time local autocorrelation is defined as
    \begin{equation*}
        \rho(t, \tau) = \sigma^{-2}(t)c(t,\tau),
    \end{equation*}
    for all $t, \,\tau\in\mathbb R$, where $\sigma^2(t) = c(t, 0) = \int_0^\infty S(u,t) du$ is the local variance of the process $X(t)$ and $c(t,\tau)$ is the local autocovariance defined above. 
\end{definition}

\setcounter{example}{0}
\begin{example} \textit{continued.}
    For the Haar \textsc{ma}($\alpha$) processes, we have 
\begin{equation}\label{eq:haar_ac}
    c(t, \tau) = \Psi_H(\alpha,\tau) = c_X(t,\tau),
\end{equation}
where $\Psi_H$ is the Haar autocorrelation wavelet, see Supplementary Material Section~S2.
Since Haar \textsc{ma} processes are stationary, the local and standard autocovariance functions are equivalent.
Further,
$\rho(t, \tau) =  c(t, \tau)$, since $\Psi_H(\alpha,0) = 1$, as shown in Proposition~\ref{prop:properties1} below.
\end{example} 

\begin{remark}
    Strictly speaking, since the Dirac measure is not continuous with respect to the Lebesgue measure,~\eqref{eq:haar_ac} is calculated using $S^{(I)}_t$ from~\eqref{eq:power_spec} as follows
    \begin{equation*}
    c(t, \tau) = \int_0^\infty \Psi(u,\tau)\,S^{(I)}_t(du) = \int_0^\infty \Psi(u,\tau)\delta_\alpha(du) = \Psi(\alpha,\tau).
    \end{equation*}
\end{remark}

To conclude this section we list some  useful properties of autocorrelation wavelets.
\begin{proposition} \label{prop:properties1}
    For all $u\in\mathbb R^+$, with $\Psi(u,\tau)$ denoting the autocorrelation wavelets and $\Psi(\tau) = \Psi(1,\tau)$, we have:
        (i) $\Psi(u,\tau)$ is symmetric as a function of $\tau$, $\Psi(u,\tau) = \Psi(u,-\tau)$;
        (ii) $\Psi(u,0)=1$;
        (iii) $\Psi(u, \tau) = \Psi(\tau/u)$;
        (iv) $\hat\Psi(u,\omega) = u |\hat\psi(u\omega)|^2$, where $\hat\Psi(u,\omega)$ is the Fourier transform of $\Psi(u, t)$;
        (v) $\int_{-\infty}^\infty\Psi(u,\tau) d\tau=0$.
\end{proposition}

\subsection{Relationship between the Wavelet Periodogram and Spectrum}

We now define the raw wavelet periodogram
and its relationship to the evolutionary wavelet spectrum. First,
we need to define the wavelet coefficients.

\begin{definition}
    The wavelet coefficients of the process $X(t)$ in terms of the wavelet basis $\psi$ are
    \begin{equation}\label{eq:coeffs}
        d(u,v) = \int_{-\infty}^\infty X(t) \psi(u,t-v) dt,
    \end{equation}
    for all $u\in\mathbb R^+$ and $v\in\mathbb R$.
\end{definition}

\begin{definition}\label{def:rwp_ct}
    The raw wavelet periodogram of the process $X(t)$ is
    \begin{equation} \label{eq:rwp}
        I(u,v) = |d(u,v)|^2,
    \end{equation}
    for all scales $u\in\mathbb R^+$ and times $v\in\mathbb R$, where $d(u,v)$ are the wavelet coefficients. We define the notation $\beta(u, v) = E \{ I(u, v) \}$, which turns out to be useful.
\end{definition}

Before going further, we introduce the continuous counterpart of the inner product matrix defined in~\cite{LSWP}, the inner product kernel. 

\begin{definition}
    For scales $u,x\in\mathbb R^+$, the inner product kernel  is
    \begin{equation*}
        A(u,x) = \int_{-\infty}^\infty \Psi(u,\tau)\Psi(x,\tau)d\tau.
    \end{equation*}
\end{definition}

\begin{proposition}\label{prop:properties2}
    The inner product kernel $A$ is continuous, symmetric and non-negative definite. Furthermore, for all $b >0$ and $x,y \in\mathbb R^+$, $A(b x, b y) = b A(x,y)$.
\end{proposition}

The following theorem establishes a relationship between the raw wavelet periodogram and the evolutionary wavelet spectrum, paralleling and extending that of~\cite{LSWP}.

\begin{theorem} \label{thm:expectation}
    Let $\beta(u, v)$ be the expectation of the raw wavelet periodogram of the locally stationary wavelet process $X(t)$, $S(u,v)$ the corresponding evolutionary wavelet spectrum and $A(u,x)$ the inner product kernel. Then
    \begin{equation}\label{eq:spec_raw}
        \left|\beta(u, v) - \int_0^\infty A(u,x)S\left(x, v\right)  dx \right| \leq \gamma \int_0^\infty x K(x)  dx,
    \end{equation}
    for all scales $u\in\mathbb R^+$ and times $v\in\mathbb R$, where $\gamma$ is a constant which depends on the decay structure of the wavelet $\psi$ and $K(\,\cdot\,)$ is the Lipschitz constant function from Definition~\ref{def:clswp}.
\end{theorem}

As for Proposition \ref{prop:autocov}, if $X(t)$ is stationary, 
then the two terms on the left hand side 
of~\eqref{eq:spec_raw} are equivalent. The expectation of the raw wavelet periodogram can be written as
\begin{equation*} 
\beta (u, v)= \int_0^\infty\int_{-\infty}^\infty S(x, y+v)\left\{\int_{-\infty}^\infty\psi(x,y-t)\psi(u,t) dt\right\}^2 dy dx,
\end{equation*}
see the proof of Theorem~\ref{thm:expectation} in Supplementary Material Section~S6, containing terms resembling the 
cross-scale autocorrelation wavelets from the discrete setting~\cite{cross-scale,replicated2}.
Hence $\int_0^\infty S\left(x, v\right)A(u,x)  dx$ can be viewed as approximating the expectation of the raw wavelet periodogram under the assumption that $S(u,t+v) = S(u,t)$ at a local level, in the same way that the local autocovariance approximates the autocovariance. Because a similar discussion follows Proposition~\ref{prop:autocov}, we omit a more detailed analysis here for brevity.

\subsection{Solving for the Spectrum and Related Quantities}

Theorem~\ref{thm:expectation} motivates the following definition:
%
%
\begin{definition}  \label{def:operatorT} The inner product operator $T_A:L^2(\mathbb R^+)\mapsto L^2(\mathbb R^+)$  is a  linear integral operator defined by
\begin{equation}\label{eq:operatorT}
    (T_A f)(x) = \int_0^\infty A(x,y) f(y) dy, 
\end{equation}
for $f\in L^2(\mathbb R^+)$ and where $A(x,y)$  is the inner product kernel.
\end{definition}

Theorem \ref{thm:expectation} and Definition \ref{def:operatorT} suggest that one could approximate
$\beta(u, v)$ by
\begin{equation} \label{eq: approx_spec}
    \beta(u, v) \approx \int_0^\infty A(u,x)S\left(x, v\right)  dx = \{ T_A S(\,\cdot\,,\, v)\}(u).
\end{equation}
Therefore, given $\beta(u, v)$ we could invert the above linear integral equation to obtain the evolutionary wavelet spectrum as a solution $\tilde S$. We now show how this could be done using the properties of the inner product kernel and Mercer's theorem~\cite{mercer}:
this is a key difference with the earlier
discrete time version in~\cite{LSWP}. 

Since $T_A$ has a continuous symmetric kernel, it is compact, bounded and self-adjoint. Further:
%
%
\begin{proposition} \label{prop:nontrivial}
    Let $T_A$ be the linear integral operator from Definition \ref{def:operatorT}. Then $T_A$ has a trivial null space, i.e. $N(T_A) = \{f\in L^2(\mathbb R^+): T_A f = 0\}=\{0\}$.
\end{proposition}

Since the inner product kernel is continuous, symmetric and non-negative definite we can apply Mercer's theorem~\cite{mercer}, see Supplementary Material Section~S1, to obtain a decomposition of $A$
of the form
\begin{equation*}
    A(x,y) = \sum_{n=1}^\infty \lambda_n \varphi_n(x)\varphi_n(y),
\end{equation*}
where $\lambda_n$ and $\varphi_n$ are the eigenvalues and eigenfunctions of $T_A$.  Therefore
\begin{equation}\label{eq:op}
    (T_A f)(x) = \sum_{n=1}^\infty \lambda_n \langle\varphi_n,f\rangle\varphi_n(x),\quad  \langle f, g\rangle = \int_0^\infty f(x)g(x) dx,
\end{equation}
for $f, g\in L^2\left(\mathbb R^+\right)$.
Since $T_A$ is self-adjoint, this allows us to obtain the solution
\begin{equation} \label{eq:approx_sol}
        \tilde S(u,t) = \sum_{n=1}^\infty \lambda_n^{-1} \langle \beta(\cdot, t), \varphi_n(\,\cdot\,)\rangle \varphi_n(u),
\end{equation}
to the linear integral equation~\eqref{eq: approx_spec},
for all $u\in\mathbb R^+$ and $t\in \mathbb R$, which is unique, since $N(T_A)=\{0\}$.

The evolutionary wavelet spectral solution~\eqref{eq:approx_sol}
can be used to obtain an estimate of the local autocovariance using~\eqref{eq:autocov_form}. Furthermore, the ability to invert $T_A$ also permits
us to invert~\eqref{eq:autocov_form}
to obtain a representation of the spectrum in terms of the local autocovariance as follows:

\begin{theorem}\label{thm:cov_inv}
Let $c(t,\tau)$ be the local autocovariance of a continuous time locally stationary wavelet process. Then the evolutionary wavelet  spectrum can be written as
    \begin{equation*}
        S(u, t) = \int_{-\infty}^\infty B(u, \tau) c(t,\tau) d\tau,\quad
        B(u, \tau) = \sum_{n=1}^\infty \lambda_n^{-1}\langle \Psi(\,\cdot\,,\tau), \varphi_n(\,\cdot\,)\rangle\varphi_n(u),
    \end{equation*}
    where $\Psi$ is an autocovariance wavelet and $\lambda_n$ and $\varphi_n$ are the eigenvalues and eigenvectors of the inner product operator $T_A$ from Definition \ref{def:operatorT}. 
\end{theorem}

\begin{example}\label{ex:white_noise}
    The expectation of the Haar raw wavelet periodogram of a Gaussian white noise process $G((s,t])\sim\mathcal{N}\{0, (t-s)\sigma^2\}$ is $ \beta(u,v) = \sigma^2$,
    for all $u\in\mathbb{R}^+$ and $v\in\mathbb{R}$.
    From 
    \begin{equation*}
        \sigma^2 = \int_0^\infty A_H(u,x) S(x, v)dx,
    \end{equation*}
    where $A_H$ is the Haar inner product kernel, we obtain
    \begin{equation*}
        S(u,v) = \frac{\sigma^2}{\log(2)u^2},
    \end{equation*}
    for all $u\in\mathbb{R}^+$ and $v\in\mathbb{R}$. Further details can be found in Supplementary Material Section S6.
\end{example}

So far, we have
ignored the fact that generally one does not work directly with $\beta(u, v)$, but rather with a realization, or in some cases a set of samples, of the locally stationary wavelet process $X(t)$, which can then be used to obtain an estimate of $\beta(u, v)$. This, of course, introduces noise into the linear integral equation. To address this problem, regularization methods such as the spectral cut-off method discussed in Supplementary Material Section~S1
could be employed.

\section{Implementation of Spectral and Autocovariance Estimators}
\label{sec:practical}

\subsection{Sampled Process Definition}
In practice, one typically has access to a discrete and finite realization of the process, rather than the entire sample path. Therefore, estimation of spectral and autocovariance quantities requires adaptations of the previous section's theoretical framework. To begin, we provide a new definition of locally stationary wavelet processes in continuous time, reintroducing a sequence of stochastic processes and the concept of rescaled-time.

\begin{definition} \label{def:clswp2}
Let $n(T)$ be a non-decreasing integer valued function of $T\in\mathbb{R}^+$, such that $n(T)\rightarrow\infty$ as $T\rightarrow\infty$, and $\mathcal{T} = \{t_0, t_1,\ldots,t_{n(T)-1}\}$ be a set of times, such that $0=t_0<t_1<\cdots<t_{n(T)-1}=T$.
A sequence of stochastic processes $\{X_T(t)\}_{t\in\mathcal{T}}$ indexed by $T$ is said to be a sampled continuous time locally stationary wavelet process if it admits the representation
\begin{equation}\label{eq:lsw2}
X_T(t) = \int_0^{J_\psi(T)}\int_{-\infty}^\infty w_T(u, v)\psi(u, v-t)L(du, dv),
\end{equation}
in the mean square sense, where $L$ is a square-integrable homogeneous L\'{e}vy basis with $E\{L(A)\}=0$ and $E\{L(A)^2\}=\text{Leb}(A)$ for all $A\in \mathcal B_b(\mathbb R^+\times\mathbb R)$,  $\{\psi_{u,v}(x), u\in\mathbb R^+, v\in\mathbb R\}$
is a continuous wavelet family, $w_T\in L^2(\mathbb R ^+, \mathbb R)$ is a sequence of location-scale amplitude functions and $J_\psi(T)$ is an increasing function which depends on the wavelet $\psi$. Further, there exists a function $W:\mathbb{R}^+\times \mathbb(0,1)\mapsto \mathbb R$ such that: (i) for each fixed scale $u\in\mathbb R^+$, $W(u,\,\cdot\,)$ is a Lipschitz continuous function with Lipschitz constant $K(u)$. The Lipschitz constant function $K$ is bounded and satisfies
\begin{equation*}
\int_0^\infty uK(u) du <\infty;
\end{equation*}
(ii) for each fixed $z\in \mathbb (0,1)$, $W(\,\cdot\,,z)$ satisfies
\begin{equation*}
\int_{0}^\infty |W(u,z)|^2 du<\infty;
\end{equation*}
(iii) for each fixed scale $u\in\mathbb R^+$ there exists a constant $C(u)$ such that for each $T\in\mathbb R^+$
\begin{equation*}
    \sup_{v\in\mathcal{T}}\left|w_T(u,v) - W\left(u,\frac{v}{T}\right)\right| \leq\frac{C(u)}{T}, \quad
    \int_0^\infty C(u) du < \infty.
\end{equation*}
\end{definition}

We have now reintroduced the concept of rescaled-time, $z=v/T\in(0,1)$, as is standard in the literature on locally stationary processes, see  \cite{dahlhaus1997} for example. In essence, this means that as $T\rightarrow \infty$ we are collecting increasing amounts of data and so are learning more and more about the local structure of the process. Assumption (c) formalizes this notion, ensuring that $w_T(u,v)$ converges to its rescaled-time limit $W(u,v/T)$, while also becoming smoother in the process. In our setting, the growing number of observations is $n(T)$, whereas $T$ now denotes the total length of the data. Thus, we still rescale the process by $T$, with the increased flexibility of being able to handle irregularly spaced observations.

The function $J_\psi(T)$ determines the maximum scale that contributes to the process. As more data becomes available, the model is able to include increasingly coarser scales and so $J_\psi(T)\rightarrow \infty$ as $T\rightarrow\infty$. We make the dependence of $J_\psi(T)$ on the structure of the underlying wavelet explicit by including $\psi$ in the subscript.

Continuous wavelets can be defined for any scale in $\mathbb R^+$ and can be shifted to any location in $\mathbb R$. This makes the continuous wavelet transform particularly well suited for handling missing observations or irregularly spaced data.

\begin{definition}
    The evolutionary wavelet spectrum of the sampled continuous time locally stationary wavelet process, $X_T(t)$, with respect to $\psi$ is 
    \begin{equation*}
        S(u,z) = |W(u,z)|^2,
    \end{equation*}
    for all $u\in\mathbb R^+$ and $z = v/T \in(0,1)$.
\end{definition}

Assumptions (ii) and (iii) of Definition~\ref{def:clswp2}, imply that, for each fixed $z\in(0,1)$,
\begin{equation*}
    \int_0^\infty S(u,z)du<\infty, \quad  S(u,z) = \lim_{T\rightarrow\infty}|w_T(u,v)|^2,
\end{equation*}
for all $u\in\mathbb R^+$ and $v\in\mathbb R$.

The local autocovariance and autocorrelation remain as before, except that we replace
$t$ with its rescaled-time version $z$, as was done for the evolutionary wavelet spectrum above. To emphasize the dependence on $T$, we now denote the autocovariance by $c_{X_T}(t,\tau)$.

Proposition~\ref{prop:autocov} highlighted the discrepancy between the `regular' autocovariance from~\eqref{eq:regularacf} and the `wavelet' autocovariance
from~\eqref{eq:autocov_form} for the theoretical process
with equality in the stationary case. In the sampled case, the autocovariances become close asymptotically.

\begin{proposition}\label{prop:autocov2} Let $\{X_T(t)\}_{t\in\mathcal{T}}$ be a sampled continuous time locally stationary wavelet process, $c_{X_T}(t,\tau)$ its autocovariance function, $c(z,\tau)$ the local autocovariance and $z=t/T$. Then 
    \begin{equation*}
        \left|c_{X_T}(t,\tau) - c(z,\tau)\right| = O\left(T^{-1}\right)
    \end{equation*}
    In other words, the autocovariance converges to the local autocovariance as $T\rightarrow\infty$.
\end{proposition}

\subsection{Estimation Theory}\label{sec:estimation}

The wavelet coefficients are defined as in~\eqref{eq:coeffs}, while the raw wavelet periodogram is now defined in rescaled-time $z=v/T$ as 
\begin{equation*}
    I\left(u,z\right) = I \left(u,\frac{v}{T}\right) = |d(u,v)|^2 = |d(u,zT)|^2.
\end{equation*}
As in Definition~\ref{def:rwp_ct}, we again use the notation $\beta(u,z) = E\left[I(u,z)\right]$. This leads to:

\begin{theorem}\label{thm:expectation2}
Let $I(u,z)$ be the raw wavelet periodogram of the sampled continuous time locally stationary wavelet process $X_T(t)$, $S(u,z)$ its corresponding evolutionary wavelet spectrum and $A(u,x)$ the inner product kernel. Then
\begin{equation}
    \label{eq:asympIS}
    \left|\beta(u, z) - \int_0^{J_\psi(T)} A(u,x)S(x,z)  dx\right|=O\left(T^{-1}\right)
\end{equation}
and
\begin{equation}
    \label{eq:asympVar}
    \left|\textup{var}\left\{I(u,z)\right\} - 2\left\{\int_0^{J_\psi(T)} A(u,x)S(x,z)  dx\right\}^2\right|=O\left(u T^{-1}\right),
\end{equation}
for all $u\in(0, J_\psi(T)]$ and $z\in(0,1)$.
\end{theorem}

Equation~\eqref{eq:asympIS} suggests 
estimating spectrum, $S$, by solving the linear integral equation
\begin{equation}\label{eq:invert}
    \beta(u, z) =\int_0^{J_\psi(T)} A(u,x) S(x,z)  dx = \{T_A S(\,\cdot\,,z)\}(u),
\end{equation}   
using estimates $\hat\beta$ of $\beta$, on which we elaborate in Supplementary Material Section~S4.2, where
\begin{equation} \label{eq:operatorT2}
    (T_A f)(x) = \int_0^{J_\psi(T)} A(x,y) f(y)  dy.
\end{equation}
The same properties hold for this linear integral operator as for its theoretical counterpart defined in~\eqref{eq:operatorT}. The use of an estimate of $\beta$ naturally introduces noise into~\eqref{eq:invert}. Hence, to compute an estimate $\hat S$ of the spectrum $S$, we use the iterative soft-thresholding algorithm developed by~\cite{daub_inv} and described in Section~\ref{sec:background}. This involves reformulating~\eqref{eq:invert} as a minimization problem with the aim of minimizing the discrepancy
\begin{equation}\label{eq:discrep}
    \Delta S(\,\cdot\,, z) = \left|\left|T_A S(\,\cdot\,,z) - \hat{\beta} (\,\cdot\,,z)\right|\right|^2_2 +\left|\left| \mu(\,\cdot\,) S(\,\cdot\,,z)\right|\right|_1,
\end{equation}
for each $z\in(0,1)$, where $||\,\cdot\,||_1$ and $||\,\cdot\,||_2$ denote the $L^1(\mathbb R^+)$ and $L^2(\mathbb R^+)$ norms respectively and $\mu > 0$ is the regularization function. By definition the evolutionary wavelet spectrum is a non-negative quantity and so we use asymmetric soft thresholding. Hence, in this setting, we reformulate the iterative scheme as
\begin{equation} \label{eq:spec_iter}
    \hat S^n(\,\cdot\,,z) = \mathcal{S}_{\mu}\left[\hat S^{n-1}(\,\cdot\,,z) + T_A\{\hat\beta(\,\cdot\,,z) - T_A \hat S^{n-1}(\,\cdot\,,z)\}\right], \quad n = 1,\ldots, N,
\end{equation}
where $N$ is the number of iterations and
\begin{equation*} \label{eq:spec_shrink}
    (\mathcal{S}_{\mu} f)(u) = \left\{\begin{array}{lll}
    f(u) - \mu(u)/2,  &&  \text{if }  f(u) \geq \mu(u)/2,\\
    0,    &&  \text{if } f(u) < \mu(u)/2,
    \end{array}\right.
\end{equation*}
for $u\in\mathbb{R}^+$, which, as shown in Section 2 of~\cite{daub_inv}, converges strongly to the unique minimizer of~\eqref{eq:discrep}.  
Unfortunately, an exact closed-form equation for estimating the evolutionary wavelet spectrum from the raw wavelet periodogram does not exist, which is a key difference from the discrete time setting of~\cite{LSWP} and other later estimators in the literature.
More details on the convergence of this algorithm in practice can be found in Supplementary Material Section~S4.3.

\begin{remark}
We have implicitly chosen to penalize the function values directly when imposing the constraint in~\eqref{eq:discrep} and, subsequently, when defining the $\mathcal{S}_{\mu}$ operator. This makes guaranteeing the non-negativity of the spectrum straightforward and is equivalent to selecting the standard basis in~\eqref{eq:shrink}. The regularization in~\eqref{eq:discrep} takes the form of a function $\mu$, allowing us to penalize scales differently if necessary. We have replaced $T^\dag_A$ with $T_A$ in~\eqref{eq:spec_iter}, as $T_A$ is by definition self-adjoint.
\end{remark}

\subsection{Smoothing the Discretized Raw Wavelet Periodogram Across Time} \label{sec:smoothing}

As for the discrete setting \cite{LSWP}, the raw wavelet periodogram is not a consistent estimator of $\beta$. A key advantage of the continuous wavelet transform is that we can compute the wavelet coefficients $d(u,v)$ on a regularly spaced grid, regardless of the time series original structure. Hence, we assume that we have computed the continuous wavelet transform at the evenly spaced discrete locations $v_j$, $j=1,\ldots,M_v = 2^\eta$ for some $\eta\in\mathbb R^+$, such that $M_v\rightarrow\infty$ as $T\rightarrow \infty$, and the corresponding estimate of the raw wavelet periodogram $\beta(u, z_j)$ for $z_j = v_j /T$. Full details can be found in Supplementary Material Section~S4.2. Standard wavelet smoothing techniques using the discrete wavelet transform~\cite{gao, DJ94, NvS97, neumann95, vSS96, nasonbook} can then be applied to each level of $\hat\beta$ across location in a straightforward manner, as in \cite{LSWP}.

\begin{remark}
    In practice, one discretizes the scales as well, however since we smooth each level of $\beta$ independently, we omit this here. Full details can be found in Supplementary Material~S4.2. 
\end{remark}

To apply wavelet shrinkage to the discretized periodogram, we introduce a separate orthonormal discrete wavelet basis of $L^2(0, 1)$. We refer to these as wavelet shrinkage basis and coefficients to distinguish them from the continuous wavelets defined in~\eqref{eq:wav_fam} and used throughout the article. 
\begin{definition}
    Let $\{\tilde\phi_{\ell_0, m}, \tilde\psi_{\ell, m}\}$ be an orthonormal discrete wavelet basis of $L^2(0,1)$ and define the mother and father wavelet coefficients of the discretized raw wavelet periodogram as
\begin{equation}\label{eq:sec_wav1}
    \tilde d^u_{\ell, m} = \frac{1}{M_v}\sum_{j=1}^{M_v} \hat\beta(u,z_j) \tilde\psi_{\ell, m}(z_j), \quad \ell = \ell_0,\ldots, \eta, \;\; m=0,\ldots, 2^\ell - 1,
\end{equation}
and
\begin{equation}\label{eq:sec_wav2}
    \tilde c^u_{\ell_0, m} = \frac{1}{M_v}\sum_{j=1}^{M_v} \hat\beta(u,z_j) \tilde\phi_{\ell_0, m}(z_j), \quad m=0,\ldots, 2^{\ell_0} - 1,
\end{equation}
respectively, where $\tilde\psi_{\ell,m}(z) = 2^{\ell/2}\tilde\psi(2^\ell z -m)$, $\tilde{\psi}$ is a wavelet function satisfying~\eqref{eq:wav_asump}, $\ell = \ell_0,\ldots, \eta$, $m=0,\ldots, 2^\ell - 1$ and $\ell_0$ is the coarsest scale included in the scheme. The superscript $u$ is used to indicate that each level $u$ is smoothed separately as a function of $z$.
\end{definition}

Equations~\eqref{eq:sec_wav1} and~\eqref{eq:sec_wav2} parallel equation (39) of~\cite{LSWP1997}. Following~\cite{nasonbook}, we set $\ell_0=3$, since the finer scales are the ones most susceptible to noise. We can then apply a soft threshold, which we define below, to the wavelet coefficients $\tilde d_{\ell, m}^{u}$ to obtain estimates $\tilde d_{\ell, m}^{u,S}$. The resulting non-linear estimate of the raw wavelet periodogram is then 
\begin{equation*}
    \tilde{\beta}(u, z_j) = \sum_{m=0}^{2^{\ell_0} - 1} \tilde{c}^u_{\ell_0, m} \tilde{\phi}_{\ell_0, m}(z_j) + \sum_{\ell=\ell_0}^{\eta}\sum_{m=0}^{2^\ell - 1} \tilde{d}^{u,S}_{\ell, m} \tilde{\psi}_{\ell, m}(z_j).
\end{equation*}
The following theorem parallels that of~\cite{LSWP}.

\begin{theorem} \label{thm:smoothing}
For locally stationary wavelet processes with homogeneous Gaussian Lévy bases, the wavelet coefficients $\tilde{d}^u_{\ell, m}$ are asymptotically normally distributed with
\begin{equation*}
    \left|E\left(\tilde d^u_{\ell, m}\right) - \int_0^1\int_0^\infty A(u,x) S(x,z) \, dx \, \tilde\psi_{\ell,m}(z) \,dz\right| = O\left(2^{\ell/2} T^{-1}\right)
\end{equation*}
and 
\begin{equation} \label{eq:smoothing_var}
    \left|\textup{var}\left(\tilde d^u_{\ell, m}\right) - 2 T^{-1}\int_0^1\left\{\int_0^\infty A(u,x) S(x,z) dx\right\}^2 \, \tilde\psi^2_{\ell, m}(z) dz\right| = O\left(2^{\ell} u T^{-2}\right).
\end{equation}
\end{theorem}

Hence, following~\cite{gao, vSS96, LSWP}, we obtain the following result.
\begin{theorem}\label{thm:smoothing_thresh}
    Define the threshold 
    \begin{equation}\label{eq:threshold_def}
        \lambda(\ell, m, u, T) = \tilde\sigma^u_{\ell, m}\log\left(T\right),\quad \tilde\sigma^u_{\ell, m} = \left\{\textup{var}\left(\tilde{d}^u_{\ell, m}\right)\right\}^{1/2}.
\end{equation} 
Then, under the assumptions of Theorem~\ref{thm:smoothing}, for each $u\in\mathbb R^+$,
\begin{equation}\label{eq:smoothing_rate}
    E\left\{\int_0^1 \left|\tilde\beta(u,z) - \beta(u,z)\right|^2 \right\} = O\left\{T^{-\frac{2}{3}}\log^2(T)\right\}
\end{equation}
\end{theorem}

The above relies on the assumption that the underlying Lévy basis is Gaussian. In other settings, one could apply the techniques described in~\cite{NvS97}. The methods relying on the Haar-Fisz transform, as developed by~\cite{FryzlewiczA, FryzlewiczB} for the discrete setting, could be extended to work in the continuous case as well.

In general, the smoothing techniques described here are applied to the logarithm of the raw wavelet periodogram instead, which tends to stabilize the variance of the process,
see~\cite{LSWP} or Remark~4.11 of~\cite{vSS96}. This mitigates the varying spectral variance and also pulls the distribution closer to normality. Therefore, for the remainder of this article, we apply any smoothing to the log raw wavelet periodogram. 

    
To estimate the noise level $\tilde\sigma^u_{\ell, m}$ in the threshold~\eqref{eq:threshold_def},~\cite{DJ94} suggest using the median absolute deviation of the wavelet shrinkage coefficients at the finest scale, divided by 0.6745, since, the finest scale coefficients are generally just noise. The scaling term of 0.6745 is
used to ensure the median absolute deviation is a consistent estimator in the normal case,
but often used for other symmetric distributions.
Alternatively,~\cite{nasonbook} propose computing the median absolute deviation on the wavelet coefficients ranging from scales 3 to the finest scale. Both methods estimate a threshold that does not depend on 
$\ell$ and $m$, but does depend on 
$u$, as the smoothing is performed separately at each scale.

\subsection{Estimating the Spectrum from the Smoothed Raw Wavelet Periodogram}
The smoothed raw wavelet periodogram $\tilde\beta$ is then used to compute an estimate of the evolutionary wavelet spectrum $\hat S$ using~\cite{daub_inv}'s iterative soft-thresholding algorithm. The initial value in~\eqref{eq:spec_iter}, is set to $S^0(u, z) = \hat\beta(u,z)$, although this can be chosen arbitrarily.
As shown in Section~4 of \cite{daub_inv}, one should select $\mu=\mu(\varepsilon)$ such that
\begin{equation*}
        \lim_{\varepsilon\rightarrow 0}\mu_\varepsilon = 0\quad\text{and}\quad \lim_{\varepsilon\rightarrow 0}\varepsilon^2 / \mu_\varepsilon = 0,
\end{equation*}
to ensure convergence to the true solution as the noise level $\varepsilon$ approaches zero. Therefore, we set $\mu_\varepsilon = \varepsilon$. In practice, we estimate this quantity using the same method as for noise estimation in wavelet shrinkage, as described in~\cite{DJ94}, enabling us to penalize each scale individually. Moreover, the iterative scheme in~\eqref{eq:spec_iter} operates on each location separately, adding greater flexibility by allowing $\mu$ to vary by location.

\subsection{Computational Complexity}
Our method uses the continuous wavelet transform, which, in a naive implementation
expends $O(n^2)$ computational effort per scale, where $n$ is the number of data points.
However, several methods have been developed that can achieve $O(n)$ instead~\cite{fast_cwt}. Smoothing the raw wavelet periodogram adds another $O(n)$ operations per scale. The iterative scheme~\eqref{eq:spec_iter}, requires $O(m^2)$ operations per iteration and location, where $m$ is number of scales. This can be reduced to $O(m)$ if the inner product operator $T_A$ is sparse~\cite{daub_inv}, which is often the case. On a standard laptop it takes approximately  thirty seconds to
simulate $1000$ realizations, compute the raw wavelet periodograms and average the 
results, for example the realizations in Section~\ref{sec:haar} next. It takes
a further five to ten seconds to smooth those estimates and compute the evolutionary
wavelet spectrum estimates.

The earlier qualitative and approximate methods of~\cite{missing2}
require the application of
$P$ lifting wavelet transforms, where $P$ corresponds to the number of lifting algorithm
trajectories, with each lifting (wavelet) transform being $O(n)$
computational effort. Hence, the total computational effort for the methods
of \cite{missing1} is $O(nP)$, where $P$ is typically in the order of thousands.
The spectral estimation methods of \cite{missing2} build on \cite{missing1}
resulting in a method which
is considerably more computationally intensive, but that work does not establish the precise
order and the code can take many hours to run.

\subsection{Non-Stationary Haar Moving Average Example} \label{sec:haar}

 We now apply our methods to realizations from a non-stationary
continuous time Haar moving average process.
The central panel of Figure~\ref{fig:haar_ma} shows
the true underlying evolutionary wavelet spectrum. See Supplementary Material Section~S5
for a complete functional specification 
 and considerably more detail.
We also display
$n(T) = 1500$ samples drawn from a single realization
and the true local autocorrelation, which, in this case, is equal to the local autocovariance.
\begin{figure}
    \centering
    \includegraphics[width=0.95\linewidth]{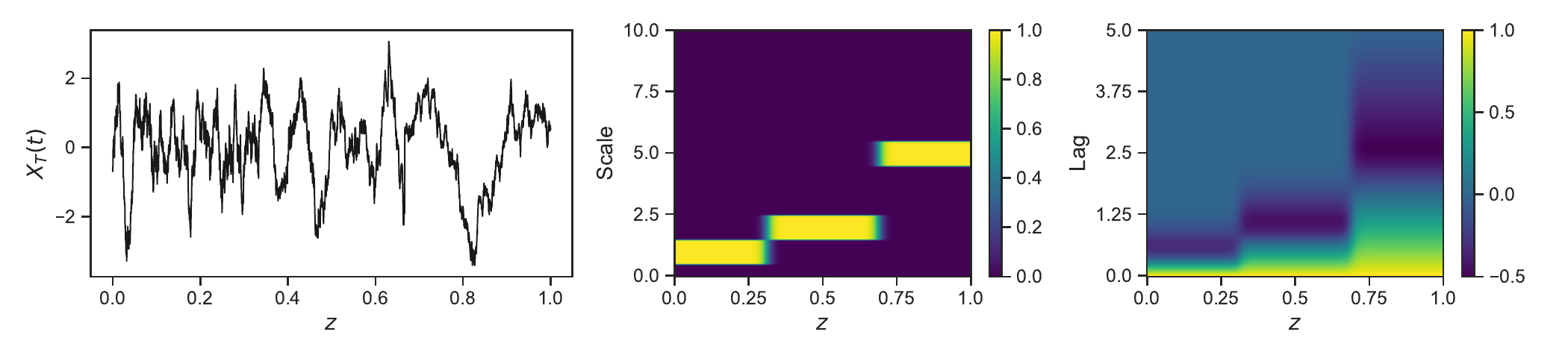}
    \caption{ Left: $n(T) = 1500$ equally-spaced samples from a single realization of continuous time
    non-stationary Haar \textsc{ma} process.
    Centre: the underlying evolutionary wavelet spectrum.
    Right: the corresponding local autocorrelation.}
    \label{fig:haar_ma}
\end{figure}

\begin{figure}
    \centering
    \includegraphics[width=1\linewidth]{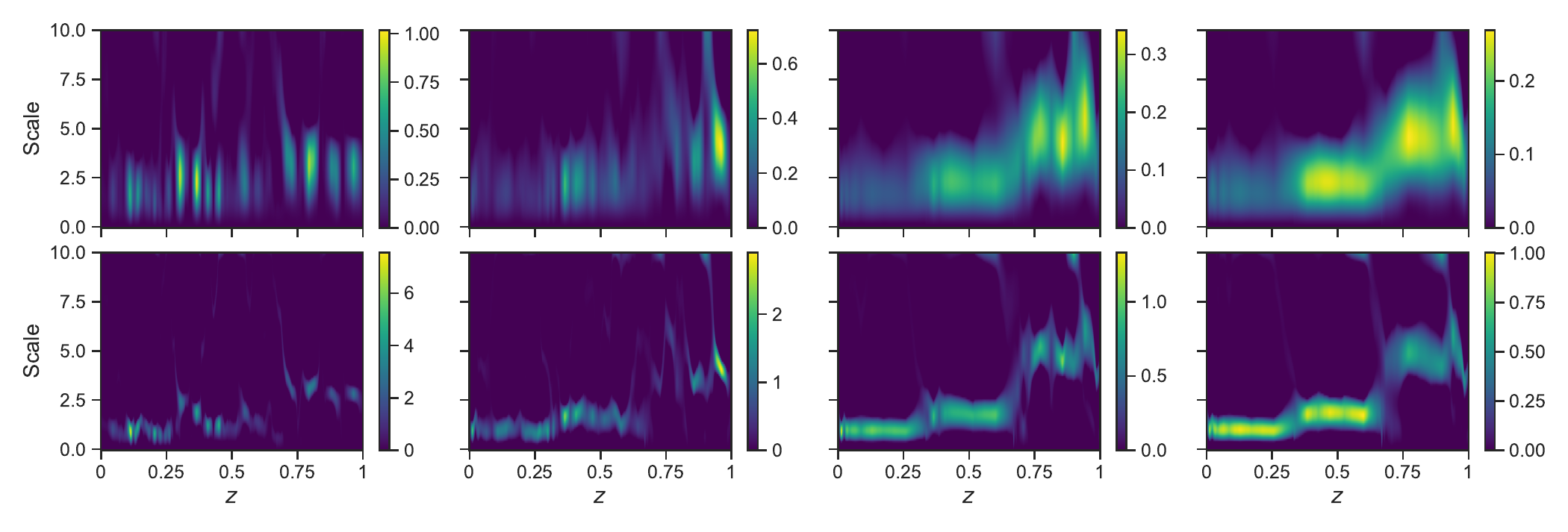}
    \caption{Haar moving average evolutionary wavelet spectral estimates using~\cite{daub_inv}'s iterative soft-thresholding algorithm. Top row: $N=100$ iterations. Bottom row: $N=10000$ iterations. From left to right: estimates computed using $R=1, 10, 100$ and $1000$ realizations.}
    \label{fig:spec_results}
\end{figure}

The estimates of the evolutionary wavelet spectrum from a single realization of $X_T(t)$ are shown on the far left of Figure~\ref{fig:spec_results}. The top row corresponds to running the
iterative soft-thresholding algorithm of~\cite{daub_inv} for $N=100$ iterations and the bottom row for $N=10000$.

Figure~\ref{fig:spec_results} also shows the results of drawing $R$ realizations each
with $n(T) = 1500$ sample, estimating a raw wavelet periodogram for each,
averaging them and then using \cite{daub_inv}'s iterative soft-thresholding algorithm to obtain a spectral estimate. Section S4.2 of the Supplementary Material defines this
mathematically in equation (S8) and discusses the practical importance of this
averaging when replicated samples are available. The second, third and fourth columns of Figure~\ref{fig:spec_results}
show estimates averaged over $R= 10$, $100$ and $1000$ realizations.

\begin{figure}
    \centering
    \includegraphics[width=1\linewidth]{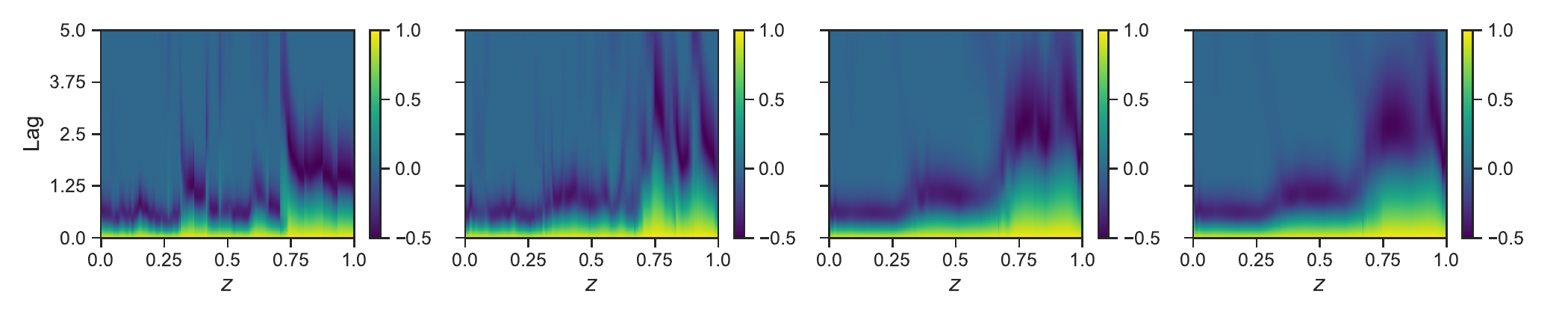}
    \caption{Haar moving average local autocorrelation 
    estimates. From left to right: estimates computed using $R=1$, $10$, $100$ and $1000$ realizations
    and $N=10000$.}
    \label{fig:estim_cov}
\end{figure}

Each estimate $\hat\beta$ is smoothed as described in Section~\ref{sec:smoothing}
using~\cite{daub} wavelets with three vanishing moments. As in ~\cite{nasonbook}, the threshold in~\eqref{eq:threshold_def} is estimated using and, applied to, the coefficients ranging from level three to the finest scale. The local autocorrelations obtained according to equation (S9) from
Supplementary Material S4.2 are displayed
in Figure~\ref{fig:estim_cov}.

As expected, results improve
as we increase both the number of samples used to estimate the raw wavelet periodogram and the
number of iterations in the~\cite{daub_inv} iterative soft-thresholding algorithm. We observe that finer scale information requires a larger number of iterations to become as visible in the estimate as information at larger scales. Spectral information also tends to deviate more from its true value at larger scales, aligning with the results in Theorem~\ref{thm:expectation2}. Even when only a single realization is employed, it is still possible to distinguish between different frequency regions. However, all estimates display a higher degree of variability in such cases, and so it may be beneficial to further investigate more appropriate smoothing techniques, an avenue left for future research.

\section{Real World Applications}
\label{sec:application}

\subsection{Infant Electrocardiogram}
\begin{figure}
    \centering
    \includegraphics[width=0.85\linewidth]{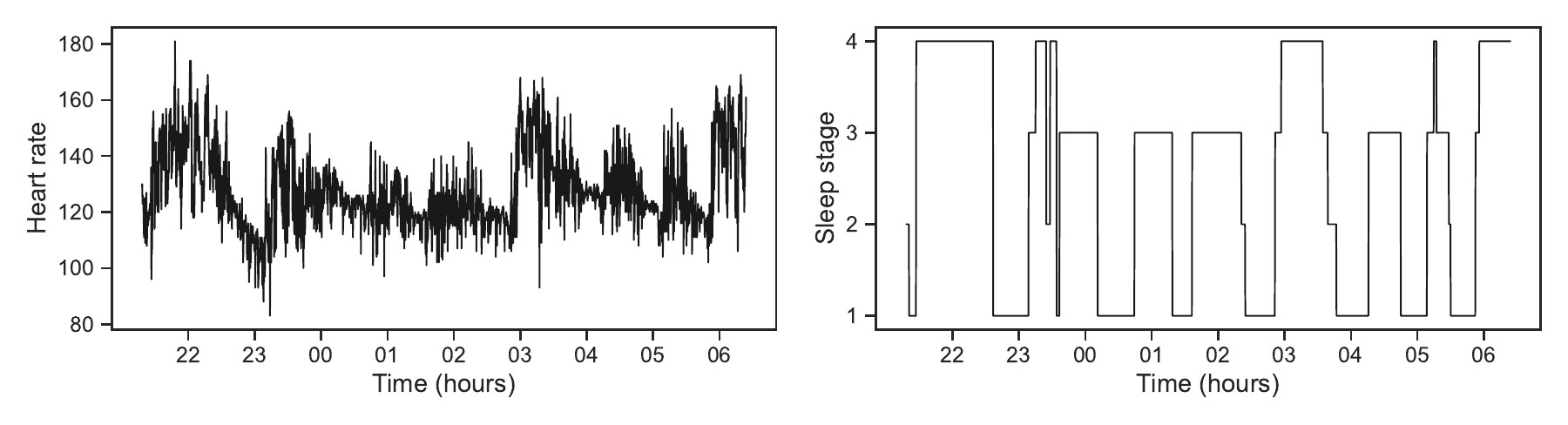}
    \caption{Left: electrocardiogram recording of a sleeping infant. Right: corresponding sleep stages.}
    \label{fig:babyECG}
\end{figure}

We apply our methods to the electrocardiogram recording of an infant, consisting of $2048$ regularly-spaced observations sampled at 1/16 Hz, used in the original discrete time work by~\cite{LSWP} and displayed in Figure~\ref{fig:babyECG}.
The {\tt BabyECG} time series used to construct this plot is freely available from
the {\tt wavethresh} \cite{wavethresh} package on CRAN.
To compare our approach to that of~\cite{missing2}, we replace $25\%$ of the observations with missing values.

\begin{figure}
    \centering
    \includegraphics[width=0.9\linewidth]{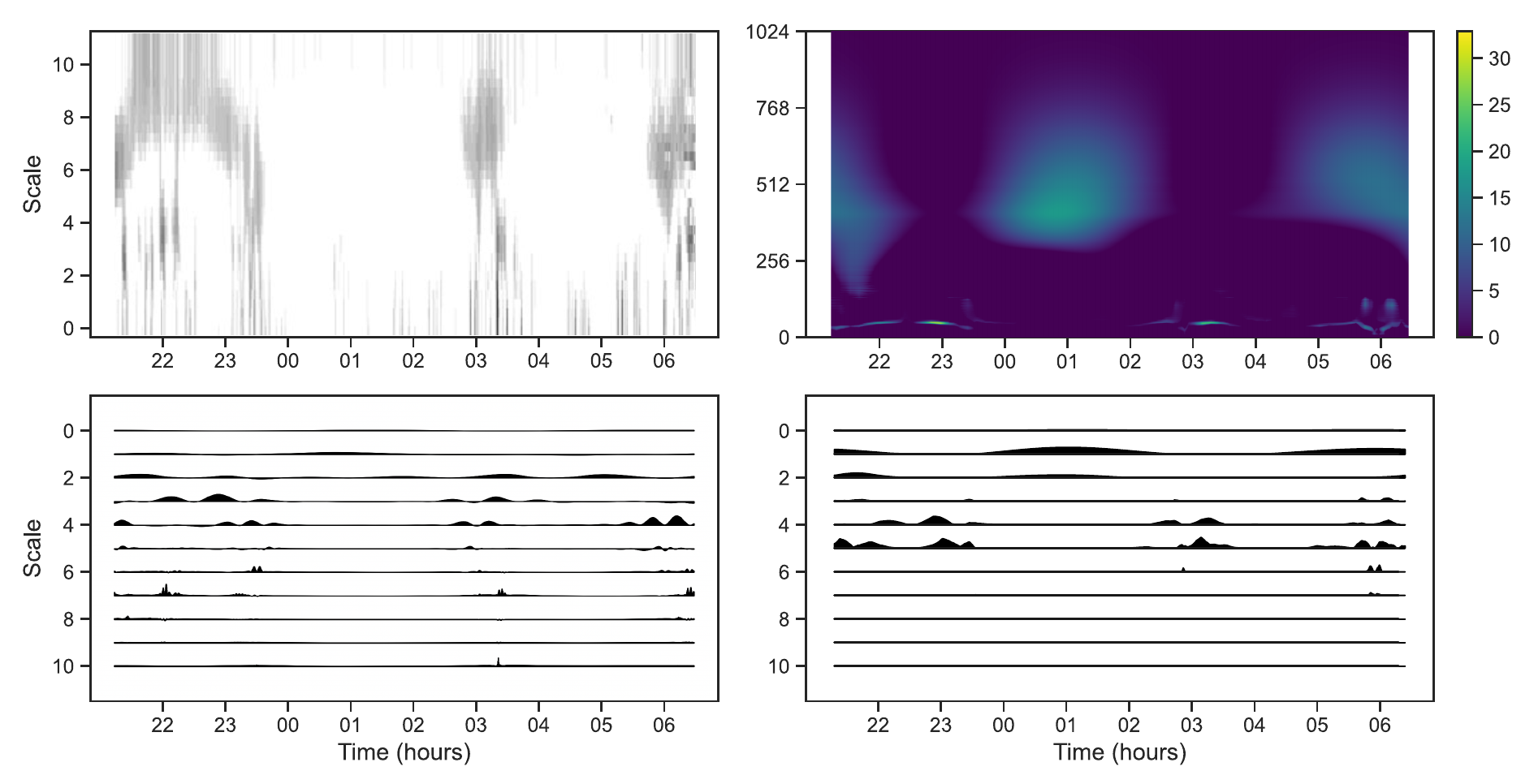}
    \caption{Sleeping infant heart rate spectral estimates. Top left: wavelet periodogram, reproduced with permission from~\cite{missing2}. Top right: continuous wavelet spectral estimate.
    Bottom left: wavelet spectrum estimated with the methods of~\cite{LSWP}. Bottom right: our discretized continuous wavelet spectrum.}
    \label{fig:baby_ecg_specs}  
\end{figure}

\cite{LSWP} and~\cite{missing2}, use the~\cite{daub} least asymmetric wavelets of order $10$ for estimation and we use the similar Ricker wavelets, see Supplementary Material S2.
We estimate our spectrum for scales ranging from $1$ to $1024 = 2^{10}$, smoothing the log periodogram using~\cite{daub} wavelets with four vanishing moments. We also discretize our spectral estimate into dyadic bins, but only to compare it with that of~\cite{LSWP} estimated on the complete time series. Both continuous and discretized estimates are displayed in the right column of Figure~\ref{fig:baby_ecg_specs}. Since different methods display scale numbering differently, all spectra in Figure~\ref{fig:baby_ecg_specs} are displayed so that the scales become coarser when moving up the y-axis.

Another appeal of our approach over the lifting-based method of~\cite{missing2} lies in the natural scale representation afforded by our construction. In contrast the scales associated to the lifting periodogram are discrete and unevenly distributed, as dictated by the time-sampling regime that governs the scale construction. Therefore,
\cite{missing2} must carry out regression over scale to produce a continuous-like lifting spectral representation. We contrast the lifting periodogram in Figure~\ref{fig:baby_ecg_specs} (top left) to its continuous spectral counterpart (top right).

Reassuringly, all estimators capture the activity burst time-localization. 
Our new spectrum captures relevant frequency behaviour in similar regions to~\cite{LSWP}, which benefits from
complete data, whereas ours does not, whilst also providing a  continuous scale representation. However, these features are less clearly elicite by the lifting scale periodogram mapping,
possibly because of scale-smoothing artefacts.
For example, the top left plot of Figure~\ref{fig:baby_ecg_specs} does not seem to have
any noticeable power at any scale around time 0100, whereas the
estimate from \cite{LSWP} and our methods (right column) do so, particularly at mid-
to coarser scales.

\subsection{Application to Heart Rate Series}
The data consists of
$n(T) = 4007$ irregularly spaced observations spanning $21569$ seconds
(approximately six hours). Mostly,
the intervals between consecutive observations range from $1$ to $10$ seconds, though a few points are separated by longer gaps, up to 130 seconds.
The discrete time methods of~\cite{LSWP} simply cannot be used, while existing methods \cite{missing2}, besides being computationally demanding, are unable to maintain a theoretical frame of reference for truly irregular data, as opposed to sampled with missingness.
Our methods are the only ones to our knowledge that can produce a genuine estimate
in reasonable time.
 
\begin{figure}
    \centering
    \includegraphics[width=0.85\linewidth]{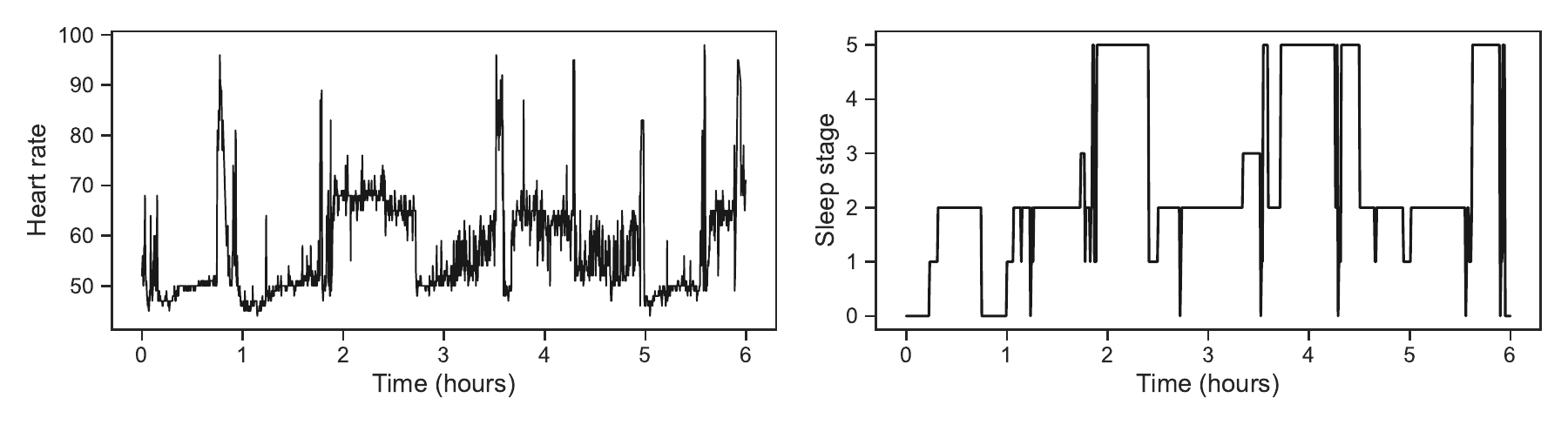}
    \caption{Left: heart rate of a sleeping individual. Right: corresponding sleep stage. Data from~\cite{hr_sl_data}.}
    \label{fig:hr_sl}
\end{figure}

\cite{sleep_paper} used wearable devices to collect heart rate data from several individuals for a study aimed at predicting sleep stages. The
data are available at~\cite{hr_sl_data}. The sleep stages are labeled~0 to~5, representing progressively deeper stages of sleep: stage~0 corresponds to wakefulness, stages~1 to~3 correspond to
non-rapid eye movement sleep, and stage~5 corresponds to rapid eye movement sleep. Figure~\ref{fig:hr_sl} displays the recordings of the chosen individual.

Due to the notable non-stationarity exhibited by the heart rate time series, Haar wavelets are used in our analysis. Whilst our approach directly models the series using the locally stationary wavelet framework, an alternative approach might involve separating a local mean component from the series first, for example by extending~\cite{vsm00, vogt, trendLSW1} to continuous time. The raw wavelet periodogram is estimated at evenly spaced locations every five seconds, across scales ranging from $10$ to $10000$ with symmetric boundary conditions.
The log periodogram is smoothed using~\cite{daub} least asymmetric wavelets with two vanishing moments. Finally, the iterative soft-thresholding algorithm of~\cite{daub_inv}
computes a spectral estimate, from which we obtain the local autocorrelation.

\begin{figure}
    \centering
    \includegraphics[width=0.95\linewidth]{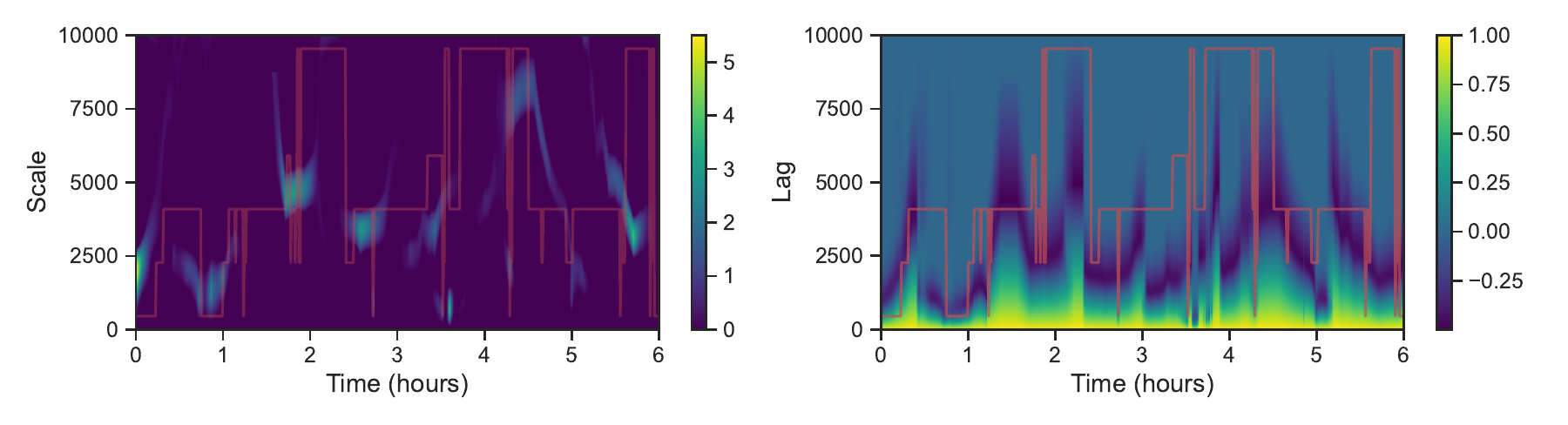}
    \caption{Estimates for the heart rate data, with the sleep stage series overlaid on top in red. Left: evolutionary wavelet spectrum. Right: local autocorrelation.}
    \label{fig:haar_results}
\end{figure}

Figure~\ref{fig:haar_results}  shows that the spectrum is concentrated primarily at coarser scales in regions where the heart rate corresponds to deeper sleep stages, such as the non-rapid eye movement stage. We can see three distinct peaks, which approximately align to
the sleep stage. A less clear
relationship might be seen in the local autocorrelation estimate. This suggests that spectral information at wider scales can be associated to deeper sleep stages.

\section{Discussion}
\label{sec:conc}

For future work:
promising directions could be improved smoothing techniques for the raw wavelet periodogram,
or alternative approaches  to evolutionary wavelet spectrum estimation.
This could involve extending the framework in Section~S4.3 of the Supplementary Material or adapting other algorithms, such as~\cite{daub_inv2, daub_inv3, fista} or~\cite{other_algos}.
One  could investigate the impact of discrete approximations on continuous time models, particularly in relation to aliasing. Aliasing occurs when high-frequency components contaminate lower frequencies due to a sampling rate which is too low. Locally stationary wavelet processes have already shown promise in detecting its absence~\cite{aliastests, aliastests2D}. Continuous time models could yield further insights into these behaviours, especially since the implementation of such methods inherently requires some form of down-sampling.

\section*{Acknowledgement}
We would like to thank the Associate Editors, reviewers and Professor Almut Veraart of Imperial College London for helpful discussions and feedback on an earlier version of this work. 

\section*{Supplementary Material}
\label{SM}
The Supplementary Material provides  background material in~S1,
autocorrelation wavelet and inner product kernels in~S2, comparison of  local and standard autocovariances on a simple example in~S3, practical considerations in~S4, an extension of the Haar moving average example in~S5, and proofs in~S6.
The code can be found at: \url{https://github.com/henrypalasciano/CLSWP}

\bibliographystyle{unsrt}
\bibliography{article}

\newpage 

\setcounter{section}{0}
\setcounter{figure}{0}
\setcounter{table}{0}
\setcounter{equation}{0}

\renewcommand{\thesection}{S\arabic{section}}
\renewcommand{\theequation}{S\arabic{equation}}
\renewcommand{\thedefinition}{S\arabic{definition}}
\renewcommand{\thetheorem}{S\arabic{theorem}}
\renewcommand{\theproposition}{S\arabic{proposition}}
\renewcommand{\theexample}{S\arabic{example}}

\vbox{%
    \hsize\textwidth
    \linewidth\hsize
    \vskip 0.1in
    \toptitlebar
    \centering
    {\LARGE\sc Supplementary Material\par}
    \bottomtitlebar
    \vskip 0.1in
}

\section{Background Material}

\subsection{Review of Locally Stationary Wavelet Processes in Discrete Time}
Locally stationary wavelet processes~\cite{LSWP, nasonbook} are time series models that have found applications in various domains, including finance~\cite{appl_fin}, economics~\cite{appl_econ}, aliasing detection~\cite{aliastests, aliastests2D}, biology and neurosciences~\cite{wav_bio,pcalsw}, energy~\cite{appl_energy1, appl_energy2}, and gravitational wave detection~\cite{gravity}. They are constructed from discrete wavelets which we now define.

\begin{definition}
    Let $\{g_k\}_{k\in\mathbb Z}$ and $\{h_k\}_{k\in\mathbb Z}$ be the high and low-pass quadrature mirror filters used in the construction of the~\cite{daub} compactly supported orthogonal continuous wavelets. Let $N_h$ be the number of non-zero elements of $\{h_k\}_{k\in\mathbb Z}$ and $N_j = (2^j - 1)(N_h-1) + 1$ for scales $j\in\mathbb N$. Then the discrete wavelets~\cite{LSWP} $\psi_j = (\psi_{j,0},\psi_{j,1},\ldots,\psi_{j, N_j-1})$ are vectors of length $N_j$ obtained using the formulae
    \begin{align*}
        \psi_{1,n} &= \sum_k g_{n-2k}\delta_{0,k} = g_n,\quad n = 0,\ldots, N_1-1,\\
        \psi_{j+1,n} &= \sum_k h_{n-2k}\psi_{j,k}, \quad n=0,\ldots,N_{j+1} - 1,
    \end{align*}
    where $\delta_{0,k}$ denotes the Kronecker delta.
\end{definition}

Locally stationary wavelet processes can then be defined as follows.

\begin{definition}
    A locally stationary wavelet process~\cite{LSWP} is a sequence of  stochastic processes $\{X_{t,T}\}_{t=0,\ldots,T-1}$, $T=2^J$ for some $j\in\mathbb N$, having the following representation in the mean square sense:
    \begin{equation*}
        X_{t,T} = \sum_{j=1}^\infty\sum_{k=-\infty}^\infty w_{j,k;T} \psi_{j,k-t}\xi_{j,k},
    \end{equation*}
    where $\{w_{j,k;T}\}_{j\in\mathbb N, k\in\mathbb Z}$ is a set of amplitudes, $\{\psi_{j,k-t}\}_{j\in\mathbb N, k\in\mathbb Z}$ are discrete wavelets and $\{\xi_{j,k}\}_{j\in\mathbb N, k\in\mathbb Z}$ is a collection of uncorrelated random variables with zero mean and unit variance. Furthermore, for each $j\in\mathbb N$ there exists a Lipschitz continuous function $W_j:(0,1)\mapsto \mathbb R$ such that: (i) $\sum_{j=1}^\infty|W_j(z)|^2 <\infty$ uniformly in $z\in (0,1)$;
    (ii) the Lipschitz constants, $L_j$, are uniformly bounded in $j$ and $\sum_{j=1}^\infty 2^jL_j<\infty$;
    (iii) there exist $\{C_j\}_{j\in\mathbb N}$ such that, for each $T$,
    \begin{equation*}
        \sup_k |w_{j,k;T} - W_j(k/T)|\leq C_j/T \quad\text{and}\quad \sum_{j=1}^\infty C_j < \infty,
    \end{equation*}
    where, for each $j$, the supremum is taken over $k=0,1,\ldots T-1$.
\end{definition}

The amplitude functions $W_j$ are defined in rescaled-time, meaning that the time index $t\in\mathbb N$ is rescaled to $z = t/T\in(0,1)$. This concept, introduced by~\cite{dahlhaus1997}, is commonly employed when working with locally stationary processes. The limit $T\rightarrow\infty$ is interpreted as increasing amounts of information becoming available about the local properties of the process, allowing for meaningful asymptotic analysis and statistical estimation of model parameters.

We now present spectral quantities:

\begin{definition}
    The evolutionary wavelet spectrum~\cite{LSWP}  of the locally stationary wavelet process $\{X_{t,T}\}_{t=0,\ldots,T-1}$ is 
    \begin{equation*}
        S_j(z) = |W_j(z)|^2
    \end{equation*}
    with respect to the wavelet family $\{\psi_{j,k}\}$, for $j\in\mathbb N$ and $z\in(0,1)$.
\end{definition}

The evolutionary wavelet spectrum provides a time-scale decomposition of the contributions to the variance of the process $X_{t,T}$. In fact, $S_j(z)$ is approximately equivalent to the variance of the process within the frequency band $[2^{-j}\pi, 2^{1-j}\pi]$, although the exact band depends on the choice of wavelet~\cite{aliastests}.

One of the fundamental concepts in the study of time series is the autocovariance function of a process. Fortunately, the autocovariance of a locally stationary wavelet process can be expressed in terms of the underlying wavelets and the evolutionary wavelet spectrum, mirroring the representation
\begin{equation*}
    c(\tau) = \int_{-\infty}^\infty \exp(i\omega \tau) dS(\omega)
\end{equation*}
of stationary time series in the Fourier domain.
This representation relies on autocorrelation wavelets, defined next.

\begin{definition}
    For $j,\ell\in\mathbb N$ and $\tau\in\mathbb Z$, the autocorrelation wavelets~\cite{LSWP} are 
    \begin{equation*}
        \Psi_j(\tau) =  \sum_k \psi_{j,k}\psi_{j,k-\tau},
    \end{equation*}
    and the corresponding inner product operator~\cite{LSWP} is
    \begin{equation}\label{eq:ipo}
        A_{j,\ell} = \sum_\tau \Psi_j(\tau)\Psi_\ell(\tau).
    \end{equation}
\end{definition}

\begin{definition}
    The local autocovariance~\cite{LSWP} of the locally stationary wavelet process $\{X_{t,T}\}_{t=0,\ldots,T-1}$ is 
    \begin{equation*}
        c(z,\tau) = \sum_{j=1}^\infty S_j(z)\Psi_j(\tau),
    \end{equation*}
    for $\tau\in\mathbb Z$ and $z\in(0,1)$.
\end{definition}

\cite{LSWP} show that the local autocovariance $c(z, \tau)$ is the limit, as $T\rightarrow\infty$, of the \say{standard} autocovariance 
\begin{equation*}
    c_T(z, \tau) = \text{cov}(X_{t,T}, X_{t+\tau,T}) = E(X_{t,T}, X_{t+\tau,T})
\end{equation*}
of the process $X_{t,T}$.

We conclude this section with the definition of the raw wavelet periodogram, which is a biased estimator of the evolutionary wavelet spectrum. To obtain an unbiased estimator, it is necessary to correct the raw periodogram using the inner product operator defined in equation~\eqref{eq:ipo}, see~\cite{LSWP} for further information.

\begin{definition}
    For $\ell\in\mathbb N$ and $m\in\mathbb Z$, the raw wavelet periodogram~\cite{LSWP} of the locally stationary wavelet process $\{X_{t,T}\}_{t=0,\ldots,T-1}$ is 
    \begin{equation*}
        I_{\ell,m} = d_{\ell,m}^2,
    \end{equation*}
    where $\{d_{\ell,m}\}$ are the non-decimated wavelet coefficients of $X_{t,T}$ given by
    \begin{equation*}
        d_{\ell,m} = \sum_{t=0}^{T-1} X_{t,T} \psi_{\ell, m-t}.
    \end{equation*}
\end{definition}

\subsection{Review of Lévy Processes}

Lévy processes and bases play a crucial role in the definition of locally stationary wavelet processes in continuous time. The contents of this section draw heavily from~\cite{walsh, spectrallevy, applebaum, almut}. 
We take $(\Omega, \mathcal{F}, \mathbb P)$ to be the underlying probability space, $\{S,\mathcal{B}(S)\}$ to be a Borel space and denote the bounded Borel sets of $\mathcal{B}(S)$ by $\mathcal{B}_b(S)$. 

\begin{definition} \label{def:le1}
    A real-valued random variable $X$ is said to be infinitely divisible~\cite{applebaum, almut} if, for all $n\in\mathbb N$, there exist independent and identically distributed random variables $Y_1^{(n)},\ldots, Y_n^{(n)}$ such that $X \stackrel{d}{=} Y_1^{(n)}+\cdots+Y_n^{(n)}$, where $\stackrel{d}{=}$ denotes equivalence in distribution.
\end{definition}

\begin{definition} A real-valued stochastic process $L=\{L(t)\}_{t>0}$ on $(\Omega, \mathcal{F}, \mathbb{P})$ is called a Lévy process~\cite{almut} if it exhibits the following properties: (i) $L(0)=0$ almost surely; 
    (ii) $L$ has independent increments: for any $n\in\mathbb N$ and any sequence of times $t_0 < t_1 < \cdots < t_n$, the random variables $L(t_0)$, $L(t_1) - L(t_0), \ldots, L(t_n) - L(t_{n-1})$ are independent;
    (iii) $L$ has stationary increments: the law of $L(t+s)-L(t)$ does not depend on $t$;
    (iv) $L$ is continuous in probability: $\forall\varepsilon >0, \displaystyle\lim_{h \rightarrow 0} \mathbb P(|L(t+h) - L(t)| \geq \varepsilon)=0 $;
    (v) $L$ has càdlàg paths.
\end{definition}

\begin{theorem} \label{thm:LK}
    Lévy-Khintchine~\cite{applebaum, almut}. Let $X$ be an $\mathbb R^n$-valued random variable. Then $X$ is infinitely divisible if there exists a vector $b\in\mathbb R^n$, a positive definite symmetric $n\times n$ matrix $A$ and a Lévy measure $\nu$ on $\mathbb R^n\backslash\{0\}$, which is to say a measure satisfying $\int_{\mathbb R^{n}\backslash\{0\}}\max(|y|^2, 1)\nu(dy)<\infty$, such that, for all $u\in R^n$,
    \begin{align*}
        \phi_X(u) &= E\{\exp(i\langle u,X\rangle)\}\\
        &=\exp{i\langle b,u \rangle - \frac{1}{2}\langle u, Au\rangle + \int_{\mathbb{R}^n\backslash\{0\}} \left\{e^{i\langle u, y\rangle} - 1 - i\langle u,y\rangle \mathbbm{1}_{B_1(0)}(y)\right\}\nu(dy)}.
    \end{align*}
    Here, $B_1(0)$ is the open ball of radius 1 centred at 0 and $\langle\,\cdot\,,\,\cdot\,\rangle$ denotes the inner product in $\mathbb R^d$. Conversely, any mapping of this form is the characteristic function of an infinitely divisible random variable on $\mathbb R^n$.
\end{theorem}

Since every Lévy process is infinitely divisible, its characteristic function also takes the above form. Therefore, a Lévy process can be viewed as a Brownian motion with drift $b$ interspersed with jumps of arbitrary size~\cite{applebaum}.

\begin{definition} A random measure~\cite{applebaum, almut} $M$ is a collection of real-valued random variables $\{M(A):A\in\mathcal{B}_b(S)\}$ such that if $A_1, A_2, \ldots$ is a sequence of disjoint elements of $\mathcal{B}_b(S)$ satisfying $\bigcup_{j=1}^\infty A_j\in \mathcal{B}_b(S)$, then 
\begin{equation*}
M\Bigg(\;\bigcup_{j=1}^\infty A_j\Bigg) = \sum_{j=1}^\infty M(A_j),
\end{equation*}
where the right hand side converges almost surely.
\end{definition}

\begin{definition}
    A random measure $M$ is said to be stationary~\cite{almut} if for all $\boldsymbol s\in S$ and any finite collection $A_1, A_2, \ldots, A_n\in\mathcal{B}_b(S)$, the random vector $\{M(A_1 + \boldsymbol s), M(A_2 + \boldsymbol s),\dots, M(A_n + \boldsymbol s)\}$ has the same law as $\{M(A_1), M(A_2),\dots, M(A_n)\}$.
\end{definition}

\begin{definition}\label{def:levybasis}
A real valued Lévy basis~\cite{almut}  on $S$ is a collection $\{L(A):    A\in \mathcal{B}_b(S)\}$ of random variables satisfying the following properties: (i) if $A$ and $B$ are disjoint subsets in $\mathcal{B}_b(S)$, then
$L(A \cup B) = L(A) + L(B)$
almost surely;
(ii) if $A$ and $B$ are disjoint subsets in $\mathcal{B}_b(S)$, then $L(A)$ and $L(B)$ are independent;
(iii) the law of $L(A)$ is infinitely divisible for all $A\in \mathcal{B}_b(S)$.
\end{definition}

Equivalently, we could define a Lévy basis $L$  to be an independently scattered, infinitely divisible random measure. Here, independently scattered means that, if $A_1,A_2,\ldots$ is a sequence of disjoint elements of $\mathcal{B}_b(S)$, then the random variables $L(A_1),L(A_2),\ldots$ are independent. See~\cite{almut} for a proof of the equivalence between the two definitions.

\begin{proposition}\label{prop:firstone} \cite{almut}.
    Let $L$ be a Lévy basis and $u\in\mathbb R$. Then for all $A\in \mathcal{B}_b(S)$, 
    \begin{align*}
        \phi_{L(A)}(u) &= E[\exp\{iuL(A)\}]\\
        & = \exp\left[iu\xi^*(A) - \frac{1}{2} u^2 a^*(A) + \int_{\mathbb{R}}\left\{e^{iux} - 1 - iux\mathbbm{1}_{[-1,1]}(x) \right\} n(dx, A)\right]
    \end{align*}
    where $\xi^*$ is a signed measure on $\mathcal{B}_b(S)$, $a^*$ is a measure on $\mathcal{B}_b(S)$ and $n(\,\cdot\,,\,\cdot\,)$ is the generalized Lévy measure.
\end{proposition}

Since the law of $L(A)$ is infinitely divisible for all $A\in \mathcal{B}_b(S)$, it is not surprising that it possesses a Lévy-Khintchine representation.

\begin{definition} 
    Let $L$ denote a Lévy basis on $\{\mathbb R^n, \mathcal{B}(\mathbb R^n)\}$. If $\xi^*$, $a^*$, $n(dx,\,\cdot\,)$, as defined in Proposition \ref{prop:firstone}, are all absolutely continuous with respect to the Lebesgue measure and their Radon-Nikodym derivatives do not depend on $\boldsymbol{z} \in S$, then $L$ is a homogeneous Lévy basis~\cite{almut}.
\end{definition}

\begin{proposition} \label{prop:le1} \cite{almut}.
    Let $L$ be a Lévy basis on $\{\mathbb R^n, \mathcal{B}(\mathbb R^n)\}$. Then $L$ is homogeneous if and only if it is stationary.
\end{proposition}

When defining continuous time locally stationary wavelet processes, we restrict our attention to homogeneous Lévy bases as the above proposition makes the concepts of homogeneity and stationarity synonymous. This is a desirable property, as we would like to control the non-stationary of the process entirely through the deterministic amplitude functions of the process.

\begin{example}
    Let $G$ be a Lévy basis. If $G(A)\sim\mathcal N\{\mu \text{Leb}(A), \sigma^2\text{Leb}(A)\}$ for all $A\in\mathcal{B}(\mathbb R^n)$, for some $\mu\in\mathbb R$ and $\sigma\in\mathbb R^+$, then $G$ is a homogeneous Gaussian Lévy basis~\cite{almut}. Gaussian Lévy bases are equivalent to the concept of white noise of~\cite{walsh}.
\end{example}

If $L$ is a homogeneous Lévy basis on $\{\mathbb R, \mathcal{B}(\mathbb R)\}$, then the process $(L_t)_{t\geq 0}$ with $L_t = L([0,t])$ is a Lévy process. Conversely, given a Lévy process $(L_t)_{t\geq 0}$, we can obtain a Lévy basis $L$ by setting $L((s,t])=L_t-L_s$. It is easy to see a similar connection between a homogeneous Gaussian Lévy basis and Brownian motion~\cite{almut}.

\subsection{Review of Linear Integral Equations and Regularization Methods}\label{subsec:review3}

Lastly, we review some relevant concepts related to linear integral equations and regularization methods. Although not strictly necessary for defining locally stationary wavelet process in continuous time, these tools are useful for obtaining an estimate of the evolutionary wavelet spectrum.
The material presented here draws from~\cite{mercer, kress, daub_inv, zemyan}. In what follows, we denote the adjoint of a linear operator $T$ by $T^\dagger$, its eigenvalues and eigenfunctions by $\lambda_n$ and $\varphi_n$ respectively, and the null space of $T$ by $N(T)$. 

\begin{definition}
    Let $S\subset \mathbb R^n$, $K:S\times S\mapsto \mathbb R$ be a continuous function and $h\in L^2(S)$. Then the linear operator $T:L^2(S)\mapsto L^2(S)$ defined by 
    \begin{equation*} 
        (T h)(x) = \int_S K(x,y) h(y)\, dy,\quad x\in S,
    \end{equation*}
    is called a linear integral operator with continuous kernel $K$~\cite{kress}. Integral operators are compact and bounded.
\end{definition}

We focus on linear integral equations of the first kind,
\begin{equation}\label{eq:int1}
Th = f,
\end{equation}
where the function $f\in L^2(S)$ and the integral operator $T$ are known, and the function $h\in L^2(S)$ is unknown. However, compact linear operators do not posses a bounded inverse, which can lead to problems when attempting to invert the above equation to recover $h$. 

When $T$ is a self-adjoint compact linear operator on a Hilbert space, $T$ has at most countably many different eigenvalues $\lambda_n$ accumulating at zero.  Denoting the orthogonal projections onto the eigenspaces by $P_n$ and assuming that the sequence of non-zero eigenvalues is ordered such that $|\lambda_1| \geq |\lambda_2| \geq\ldots$, we can write
    \begin{equation*}
        T = \sum_{n=1}^\infty \lambda_n P_n
    \end{equation*}
in the sense of norm convergence. Under certain conditions, Mercer's theorem~\cite{mercer}  gives us a decomposition of the kernel in a similar form.
\begin{theorem}\label{thm:mercer}
Mercer's Theorem~\cite{mercer,buescu}. Let $T$ be a linear integral operator defined by
\begin{equation*}
    (T h)(x) = \int_S K(x,y)  h(y)\,dy
\end{equation*}
where $K$ is a continuous symmetric non-negative definite kernel. Then the kernel $K$ admits an expansion in terms of the eigenvalues $\lambda_n$ and eigenfunctions $\varphi_n$ of $T$ as
\begin{equation*}
       K(x,y) = \sum_{n=1}^\infty  \lambda_n\varphi_n(x)\varphi_n(y),
\end{equation*}
where the convergence is absolute and uniform. The non-negative definite kernel implies non-negative eigenvalues, which decay at a rate of $O(n^{-1})$.
\end{theorem}

When $T$ is a self-adjoint linear integral operator, we can use Mercer's theorem~\cite{mercer} to obtain a solution to~\eqref{eq:int1} of the form
\begin{equation*}
        h = \sum_{n=1}^\infty \frac{1}{\lambda_n} \langle f, \varphi_n\rangle \varphi_n,
\end{equation*}
provided that $f\in N(T)^\perp$, where $\perp$ denotes the orthogonal complement, and 
    \begin{equation*}
        \sum_{n=1}^\infty \frac{1}{\lambda_n^2} |\langle f, \varphi_n\rangle |^2<\infty.
    \end{equation*}

Typically, one encounters noise contamination of the form $Th = f = g + \varepsilon$, where $g$ denotes the true function and $\varepsilon$ the noise. The effects of this noise on the solution can be amplified due to the unboundedness of $T^{-1}$, particularly for smaller eigenvalues. To address this issue, a simple regularization method is the spectral cut-off, which replaces the solution with
\begin{equation*}
        h_\mu = \sum_{\lambda_n^2>\mu}\frac{1}{\lambda_n}\langle f,\varphi_n\rangle\varphi_n,
    \end{equation*}
for some $\mu>0$, basically eliminating the noise amplification due to the smaller eigenvalues.

Alternatively,~\eqref{eq:int1} can be reformulated as a minimization problem, with the aim of minimizing
\begin{equation*}
\Delta h = ||Th - f||^2.
\end{equation*}
A pseudo-solution of the form $\tilde{h} = (T^\dagger T)^{-1} T^\dagger f$ can be obtained, assuming that $T$ has a trivial null space, $N(T) = \{0\}$. However, the unboundedness of $T^{-1}$ and the presence of noise can again lead to difficulties. To attempt to deal with this problem, Tikhonov regularization, also known as ridge regression~\cite{elements}, is commonly employed. This involves minimizing 
\begin{equation*}
\Delta h = ||Th - f||^2 + \mu ||h||^2,
\end{equation*}
for some $\mu>0$, in order to obtain a solution of the form
\begin{equation*}
    h_\mu = (T^\dagger T+\mu I)^{-1}T^\dagger f,
\end{equation*}
essentially penalizing solutions with large norms. When $T$ is self-adjoint, this can be written as
\begin{equation*}
    h_\mu = \sum_{n=1}^\infty \frac{\lambda_n}{\lambda_n^2 + \mu}\langle f, \varphi_n\rangle \varphi_n.
\end{equation*}

\cite{daub_inv} instead propose to minimize
\begin{equation*}
    \Phi_{\mu, p}(h) = ||Th-f||^2 + \sum_{\gamma\in\Gamma}\mu_\gamma |\langle h, e_\gamma\rangle|^p,
\end{equation*}
where $\{e_\gamma\}_{\gamma\in\Gamma}$ is an orthogonal basis of the Hilbert space $H$,  $\mu=\{\mu_\gamma\}_{\gamma\in\Gamma}$ is a vector of weights, which act as regularization parameters, and $1\leq p\leq 2$. The solution for the case $p=1$ can be obtained via an iterative scheme given by
\begin{equation} \label{eq:iter_scheme2}
    h^n = \mathcal{S}_{\mu}\left\{ h^{n-1} + T^\dagger(f-Th^{n-1})\right\},
\end{equation}
where 
\begin{equation*}
    \mathcal{S}_{\mu} (h) = \sum_\gamma \mathcal{S}_{\mu_\gamma}\left(\langle h, e_\gamma\rangle \right) e_\gamma
\end{equation*}
and 
\begin{equation} \label{eq:sym2}
    \mathcal{S}_{\mu_\gamma}(x)=\left\{\begin{array}{lll}
    x-\mu_\gamma/2, && \text{if } x\geq\mu_\gamma/2, \\
    0,  && \text{if } |x|<\mu_\gamma/2,\\
    x+\mu_\gamma/2, && \text{if } x\leq-\mu_\gamma/2.
    \end{array}\right.
\end{equation}
In the literature, this method is known as the iterative soft-thresholding algorithm. \cite{daub_inv} also discuss an extension of this algorithm for setting in which one may want to penalize the positive and negative coefficients differently. This is done by replacing each $\mu_\gamma\in\mu$ with the positive weights $(\mu_\gamma^+, \mu_\gamma^-)$ and~\eqref{eq:sym2} with
\begin{equation*} 
    \mathcal{S}_{\mu_\gamma^+,\mu_\gamma^-}(x)=\left\{\begin{array}{lll}
    x-\mu_\gamma^+/2, && \text{if } x\geq\mu_\gamma^+/2, \\
    0,  && \text{if } -\mu^-_\gamma/2<x<\mu^+_\gamma/2,\\
    x+\mu^-_\gamma/2, && \text{if } x\leq-\mu^-_\gamma/2.
    \end{array}\right.
\end{equation*}
We conclude this section by stating the two main results of~\cite{daub_inv}.
\begin{theorem}
    Consider the iterative scheme defined by~\eqref{eq:iter_scheme2}, for $n=1,2,\ldots$ and arbitrarily chosen $h^0\in H$. Then $h^n$ converges strongly to the unique minimizer of $\Phi_{\mu, p}$.
\end{theorem}

\begin{theorem}\label{thm:conv_cond}
    If $\mu = \mu_\varepsilon$, where $\varepsilon$ denotes the noise level, satisfies
    \begin{equation*}
        \lim_{\varepsilon\rightarrow 0}\mu_\varepsilon = 0\quad\text{and}\quad \lim_{\varepsilon\rightarrow 0}\varepsilon^2 / \mu_\varepsilon = 0,
\end{equation*}
    then
    \begin{equation*}
        \lim_{\varepsilon\rightarrow 0} \left[\sup_{||f - Th||\leq\varepsilon} ||h^*_{\mu_\varepsilon} - h||\right] = 0,
    \end{equation*}
    for any $h\in H$, where $h^*_{\mu_\varepsilon}$ denotes the unique minimizer of $\Phi_{\mu, p}$.
\end{theorem}

\section{Some Autocorrelation Wavelets and Inner Product Kernels}
\subsection{Haar Wavelet} 
The Haar wavelet is one of the simplest and most commonly used wavelets. There are very well localized in time, however, at the cost of being discontinuous~\cite{wav_book}.

\begin{definition}
    The Haar wavelet~\cite{wav_book} at scale $u\in\mathbb R^+$ and location $t\in\mathbb R$ is 
    \begin{equation*}
      \psi_H(u, t) = u^{-1/2} \left\{ \mathbbm{1}_{[0,u/2]}(t) - \mathbbm{1}_{[u/2, u]}(t)\right\},
    \end{equation*}
    where $\mathbbm{1}$ is the usual indicator function.
\end{definition}

\begin{proposition} \label{prop:haar1}
    For scale $u\in\mathbb R^+$ and lag $\tau\in\mathbb R$, the Haar autocorrelation wavelet is
    \begin{equation*}
        \Psi_H(u,\tau) = \left\{\begin{array}{lll}
         1-3|\tau|/u,   &&  0<|\tau| \leq u/2,\\
         |\tau|/u - 1,  &&  u/2 < |\tau| \leq u,\\
         0, && \text{otherwise.}
        \end{array}\right.
    \end{equation*}
\end{proposition}

\begin{proposition} \label{prop:haar2}
    For scales $u,x\in\mathbb R^+$, the Haar inner product kernel is
    \begin{equation*}
        A_H(u,x) = \left\{\begin{array}{lll}
        \frac{x^2}{2u},  &&  0<x \leq u/2,\\[1mm]
        2x - u + \frac{u^2}{6x} - \frac{5x^2}{6u},    && u/2 < x \leq u,\\[1mm]
        2u - x + \frac{x^2}{6u} - \frac{5u^2}{6x},     && u < x \leq 2u, \\[1mm]
        \frac{u^2}{2x}, && x \geq 2u.
        \end{array}\right.
    \end{equation*}
\end{proposition}

\begin{figure}
    \centering
    \includegraphics[width=0.8\linewidth]{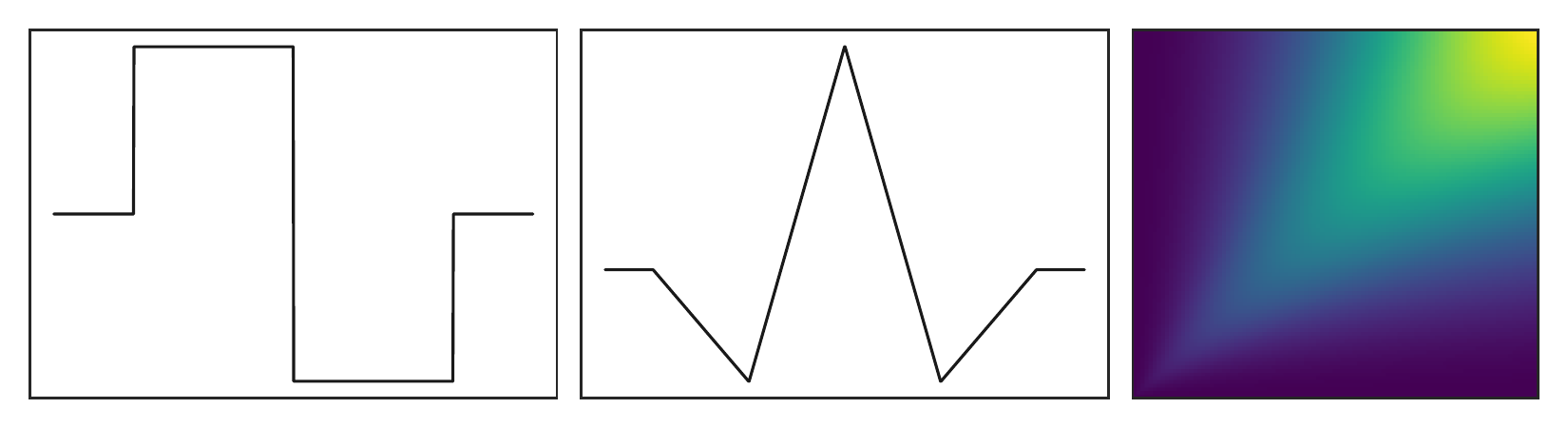}
    \caption{Left: Haar wavelet. Centre: Haar autocorrelation wavelet. Right: Haar inner product kernel.}
    \label{fig:haar}
\end{figure}
The Haar wavelet, autocorrelation wavelet and inner product kernel
are shown in Figure~\ref{fig:haar}.

\subsection{Ricker Wavelet}
The Ricker wavelet is defined as the negative normalized second derivative of a Gaussian function. Ricker wavelets do not have compact support, however, its tails decay rapidly to zero, similar to a Gaussian probability density function~\cite{ricker}.
\begin{definition}
    The Ricker wavelet~\cite{ricker} at scale $u\in \mathbb R^+$ and location $t\in\mathbb R$ is
    \begin{equation*}
        \psi_R(u,t) = \frac{2}{\pi^{1/4}(3 u)^{1/2}} \left(1-\frac{t^2}{u^2}\right)e^{-\frac{t^2}{2u^2}}.
    \end{equation*}
\end{definition}

\begin{proposition} \label{prop:ricker1}
    For scale $u\in\mathbb R^+$ and lag $\tau\in\mathbb R$, the Ricker autocorrelation wavelet is
    \begin{equation*}
        \Psi_R(u,\tau) = \left(1+\frac{\tau^4}{12 u^4} - \frac{\tau^2}{u^2}\right) e^{-\frac{\tau^2}{4u^2}}.
    \end{equation*}
\end{proposition}

\begin{proposition} \label{prop:ricker2}
    For scales $u,x\in\mathbb R^+$, the Ricker inner product kernel is
    \begin{equation*}
        A_R(u,x) = \frac{70\pi^{1/2}\, u^5 x^5}{3(u^2 + x^2)^{9/2}}.
    \end{equation*}
\end{proposition}

\begin{figure}
    \centering    \includegraphics[width=0.8\linewidth]{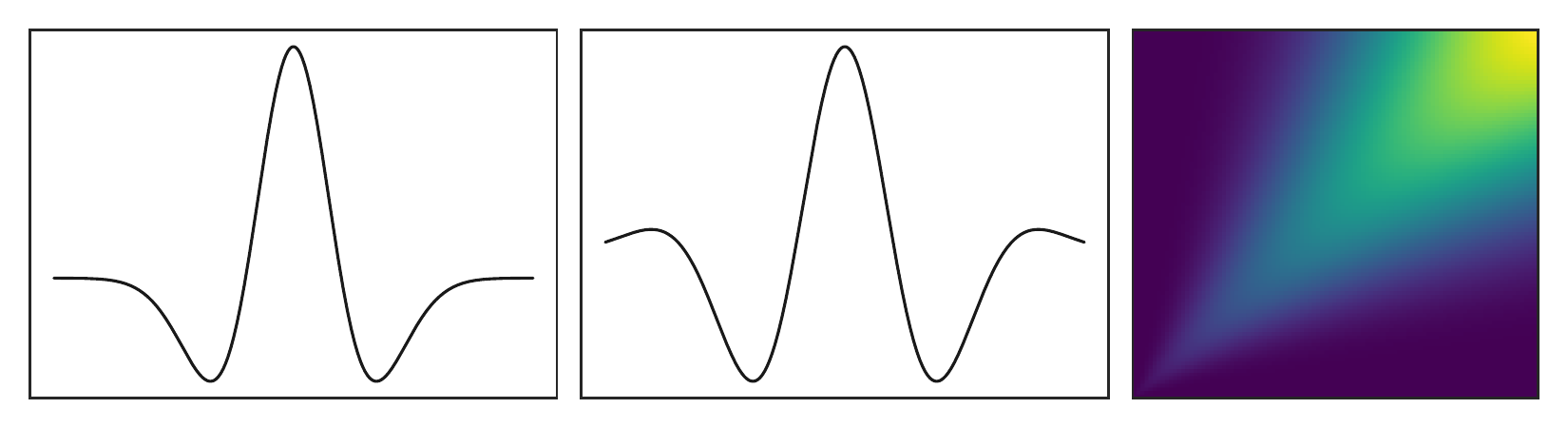}
    \caption{Left: Ricker wavelet. Centre: Ricker autocorrelation wavelet. Right: Ricker inner product kernel.}
    \label{fig:ricker}
\end{figure}

The Ricker wavelet, autocorrelation wavelet and inner product kernel
are shown in Figure~\ref{fig:ricker}.

\begin{remark}
    Several other wavelets~\cite{wavelets_chui} can be constructed from the derivatives of the Gaussian function $C\exp(-t^2)$, with the normalization constant $C$ chosen such that the wavelet function $\psi_G(t)$ satisfies~\eqref{eq:wav_asump}. For example, the wavelet obtained from the first derivative of the Gaussian is
    \begin{equation*}
        \psi_{G_1}(u, t) = -\frac{2^{5/4} t}{\pi^{1/4} u}e^{-\frac{t^2}{u^2}},
    \end{equation*}
    which can be interpreted as a continuous version of the Haar wavelet. Autocorrelation wavelets and inner product kernels can be derived in a similar fashion to those of the Ricker wavelet.
\end{remark}

\subsection{Morlet Wavelet}
The real valued Morlet wavelet is constructed by combining an oscillatory cosine term with a Gaussian envelope, which is localized in both time and frequency. There also exist complex Morlet wavelets, which combine both sine and cosine oscillatory functions, however, in this article, our focus is on real-valued wavelets only.

\begin{definition}
    The Morlet wavelet~\cite{morlet1, morlet2} at scale $u\in \mathbb R^+$ and location $t\in\mathbb R$ is
    \begin{equation*}
        \psi_M(u,t) = \frac{2^{1/2}}{\pi^{1/4}} \cos\left(\frac{\eta^{1/2}t}{u}\right)e^{-\frac{t^2}{2 u^2}}, 
    \end{equation*}
    where $\eta = 2\pi^2 / \log(2)$.
\end{definition}

\begin{proposition}\label{prop:morlet}
    For scale $u\in\mathbb R^+$ and lag $\tau\in\mathbb R$, the Morlet autocorrelation wavelet is
    \begin{equation*}
        \Psi_M(u, \tau) = \left\{e^{-\eta} + \cos\left(\frac{\eta^{1/2}\tau}{u}\right)\right\}e^{-\frac{\tau^2}{4u^2}}, 
    \end{equation*}
    where $\eta = 2\pi^2 / \log(2)$.
\end{proposition}

\begin{proposition} \label{prop:morlet2}
    For scales $u,x\in\mathbb R^+$, the Morlet inner product kernel is
    \begin{equation*}
        A_M(u,x) = 2e^{-2\eta}xu\left(\frac{\pi}{u^2+x^2}\right)^{1/2}\left[1 + 2e^{\eta/2}\cosh\left\{\frac{\eta(u^2-x^2)}{2(u^2 + x^2)}\right\} + e^{\eta}\cosh\left(\frac{2\eta ux}{u^2 + x^2}\right)\right],
    \end{equation*}
    where $\eta = 2\pi^2 / \log(2)$.
\end{proposition}

\begin{figure}
    \centering
    \includegraphics[width=0.8\linewidth]{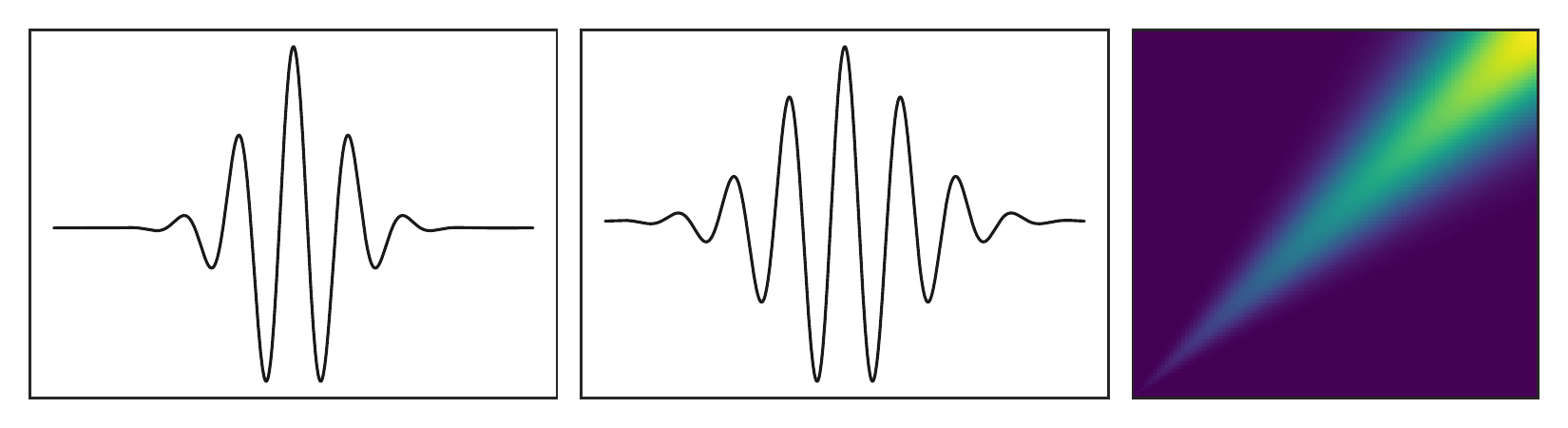}
    \caption{Left: Morlet wavelet. Centre: Morlet autocorrelation wavelet. Right: Morlet inner product kernel.}
    \label{fig:morlet}
\end{figure}
The Morlet wavelet, autocorrelation wavelet and inner product kernel
are shown in Figure~\ref{fig:morlet}.

\subsection{Shannon Wavelet}
The Shannon wavelet is the complete opposite of the Haar wavelet. It is not well localized in time and, in fact, does not have a compact support. However, it is band-limited in the frequency domain~\cite{wav_book}.
\begin{definition}
     The Shannon wavelet~\cite{wav_book} at scale $u\in \mathbb R^+$ and location $t\in\mathbb R$ is
    \begin{equation*}
        \psi_S(u, t) = \frac{u^{1/2}}{\pi t} \left\{\sin{\left(\frac{2\pi t}{u}\right)} - \sin{\left(\frac{\pi t}{u}\right)}\right\}.
    \end{equation*}
    Its Fourier transform~\cite{wav_book} is
    \begin{equation*}
        \hat\psi_S(u, \omega) = u^{1/2}\left\{\mathbbm{1}_{\left[-2\pi/u, -\pi/u\right]}(\omega) + \mathbbm{1}_{\left[\pi/u, 2\pi/u\right]}(\omega) \right\}.
    \end{equation*}
\end{definition}

\begin{proposition}\label{prop:shannon1}
    For scale $u\in\mathbb R^+$ and lag $\tau\in\mathbb R$, the Shannon autocorrelation wavelet is
    \begin{equation*}
        \Psi_S(u, \tau) = u^{1/2} \psi_S(u, \tau).
    \end{equation*}
\end{proposition}

\begin{proposition} \label{prop:shannon2}
    For scales $u,x\in\mathbb R^+$, the Shannon inner product kernel is
    \begin{equation*}
        A_S(u, x) = \left\{\begin{array}{lll}
        2x - u, && u/2 \leq x \leq u,\\
        2u-x,  &&  u < x \leq 2u,\\
        0, && \text{otherwise.}
        \end{array}\right.
    \end{equation*}
\end{proposition}

\begin{figure}
    \centering
    \includegraphics[width=0.8\linewidth]{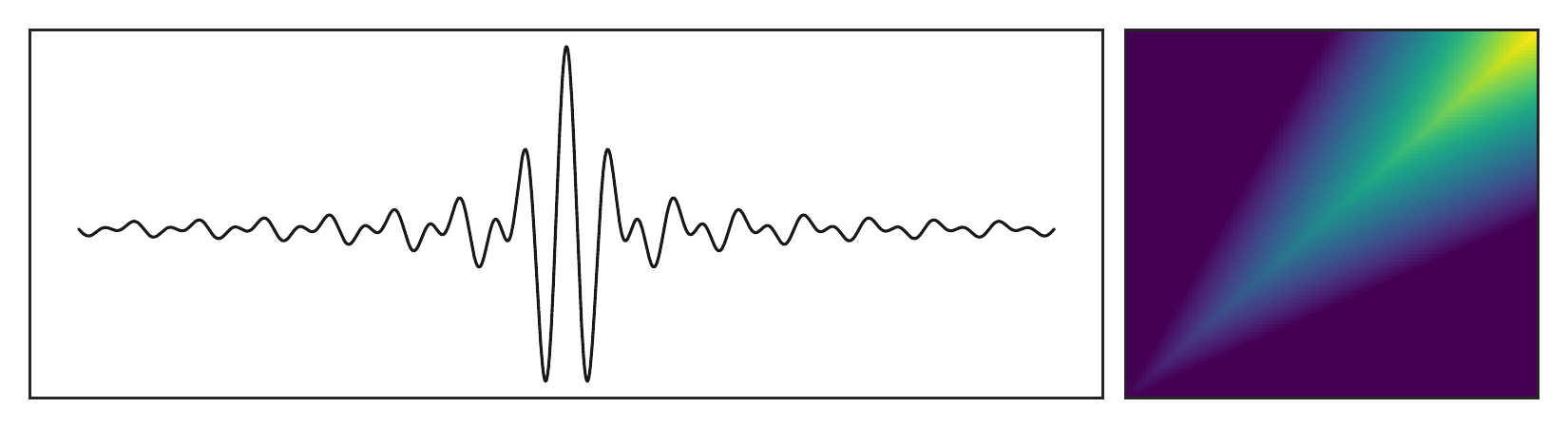}
    \caption{Left: Shannon wavelet and autocorrelation wavelet. Right: Shannon inner product kernel.}
    \label{fig:shannon}
\end{figure}

The Shannon wavelet, autocorrelation wavelet and inner product kernel
are shown in Figure~\ref{fig:shannon}.

\section{A Comparison of the Local and Standard Autocovariance Functions}
For $t,\tau\in\mathbb R$, the local and standard  autocovariance functions are 
\begin{equation}\label{eq:a1}
        c(t, \tau) = \int_0^\infty S(u,t)\Psi(u,\tau) du
\end{equation}
and
\begin{equation}\label{eq:a2}
        c_X(t, \tau) = \int_0^\infty\int_{-\infty}^\infty S(u,t+v)\psi(u,v)\psi(u,v-\tau) dv  du,
\end{equation}
respectively, as defined in Section~\ref{sec:lswtheory}. Proposition~\ref{prop:autocov} established a global bound on the discrepancy between \eqref{eq:a1} and \eqref{eq:a2}. One can also obtain a bound which depends on the lag $\tau$,
\begin{equation}\label{eq:tau_b}
    |c_X(t, \tau) - c(t,\tau)| \leq \gamma \left|\int_0^\infty u K(u) \Psi(u,\tau) du\right|,
\end{equation}
as derived in the proof of Proposition~\ref{prop:autocov}. From this bound it becomes clear that the area  over which the discrepancy has an effect decreases to zero as $\tau$ increases.

\begin{figure}
    \centering
    \includegraphics[width=0.9\linewidth]{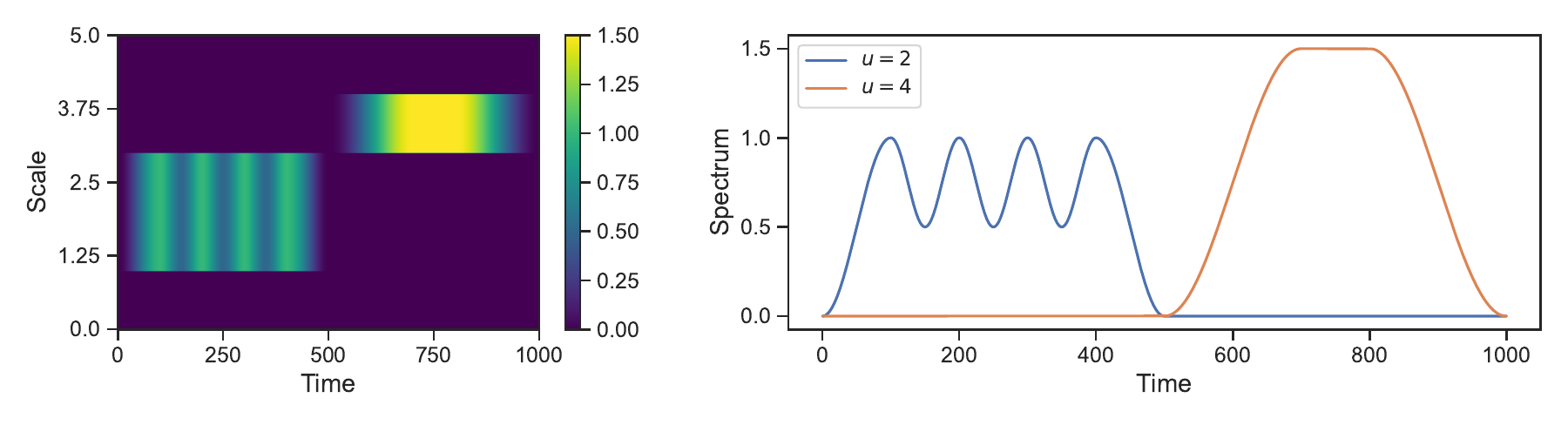}
    \caption{The full evolutionary wavelet spectrum and the evolutionary wavelet spectrum at scales $u=2$ and $u=4$.}
    \label{fig:spec_example}
\end{figure}

We now study the difference between the two quantities through the use of a simple example. The evolutionary wavelet spectrum we consider is
\begin{equation*}
    S(u,t) = \left\{\begin{array}{lll}
    \frac{1}{2}\left\{\cos\left(\frac{\pi t}{100} - 1\right) + 1\right\}, &&  \text{for }  0\leq t\leq 100,\\[1mm]
    \frac{1}{4}\left\{\cos\left(\frac{\pi t}{50}\right) + 3\right\}, &&  \text{for } 100< t\leq 400,\\[1mm]
    \frac{1}{2}\left\{\cos\left(\frac{\pi t}{100}\right) + 1 \right\}, &&  \text{for }   400 < t\leq 500,
\end{array}\right.
\end{equation*}
for $1\leq u \leq 3$, 
\begin{equation*}
S(u,t) = \left\{\begin{array}{lll}
    \frac{3}{4} \left\{\cos\left(\frac{\pi t}{200} + \frac{1}{2}\right) + 1\right\},  &&  \text{for } 500< t\leq 700,\\[1mm]
    \frac{3}{2}, && \text{for } 700< t\leq 800,\\[1mm]
    \frac{3}{4}\left\{\cos\left(\frac{\pi t}{200}\right) + 1\right\}, && \text{for } 800< t\leq 1000,
    \end{array}\right.
\end{equation*}
for $3 < u \leq 4$ and zero otherwise. This is displayed in Figure~\ref{fig:spec_example}, as well as a plot of the spectrum at scales $u=2$ and $u=4$. The evolutionary wavelet spectrum displays a wide variety of different behaviours. To compare \eqref{eq:a1} and \eqref{eq:a2} over a wider range of scales, we also redefine the same spectrum but with the scales replaced by factors of $4$ and $10$. We take the underlying wavelets to be Ricker wavelets. There is no particular reason for this choice and this comparison can easily be repeated for other wavelets as well. Although in theory $\psi_R(x) > 0$ for all $x\in\mathbb{R}$, the tails decay rapidly enough that we can assume $\psi_R(x)\approx 0$ for $x\notin (-4, 4)$.

\begin{figure}
    \centering
    \includegraphics[width=0.9\linewidth]{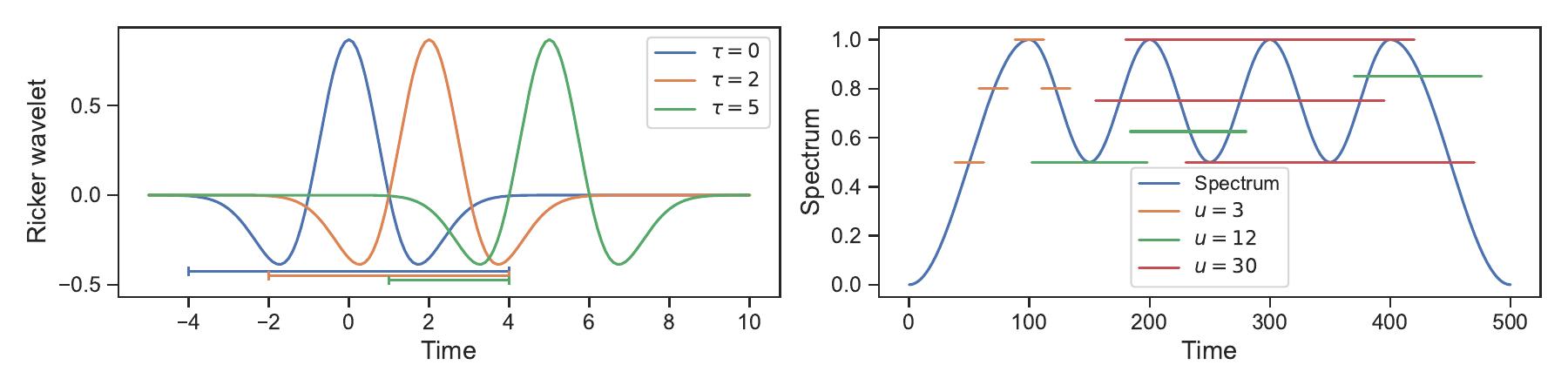}
    \caption{Left: the Ricker wavelet for scale $u = 1$ and different lags $\tau$, with horizontal lines denoting the region of overlap. Right: the local autocovariance approximates the spectrum with a constant in an interval of length at most equal to the width of the underlying wavelets. The horizontal lines display this for different scales $u$ and times $t$.}
    \label{fig:visual}
\end{figure}

As discussed in the main text, the local autocovariance at each time $t\in\mathbb R$ is computed under the assumption that the spectrum remains constant at $S(u,t)$ over the domain of the underlying wavelet. Figure \ref{fig:visual} provides further insight into this. On the left we depict the Ricker wavelet for different values of $\tau$, highlighting that only the region over which they overlap contributes to the autocovariance of the process, this being maximal at $\tau = 0$. This also offers a visual interpretation to the bound in \eqref{eq:tau_b}. On the right, we display a portion of the spectrum. Horizontal lines represent the value used in computing the local autocovariance for a range of scales $u$, with the width of this region increasing linearly with the scales.

\begin{figure}
    \centering
    \includegraphics[width=0.9\linewidth]{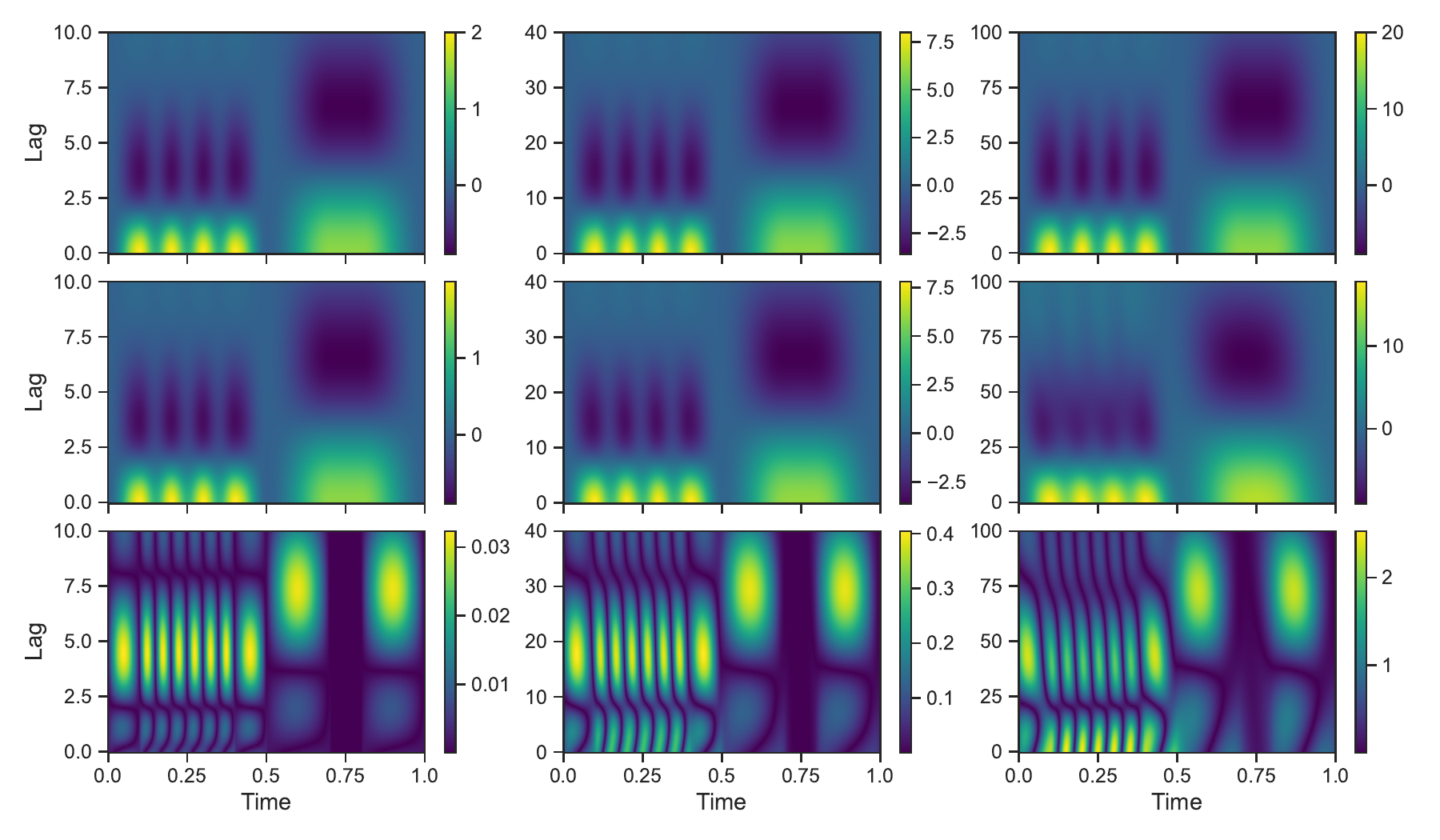}
    \caption{Top row: the local autocovariance. Central row: the standard autocovariance.  Bottom row: absolute difference between local and standard autocovariance. From left to right: computations carried out using the evolutionary wavelet spectrum defined for scales $0< u\leq 5$, $0< u\leq 20$ and $0< u\leq 100$.}
    \label{fig:spec_example_results}
\end{figure}

We compute the local autocovariance and the standard autocovariance for the evolutionary wavelet spectra defined above and present the results in Figure~\ref{fig:spec_example_results}. The two quantities are almost indistinguishable in all cases, with the two computed from the spectrum defined using the coarsest scales displaying the most significant differences, as anticipated. As illustrated in the right panel of Figure~\ref{fig:visual}, the wavelet width at these scale covers a significant portion of the spectrum. This scenario is unrealistic in practice, as we expect the spectrum to vary at a much slower rate at such wide scales. Despite this, the two quantities remain nearly identical, suggesting that in practice the local autocovariance and the standard autocovariance are generally indistinguishable, as also implied by the theoretical properties. The bottom row of Figure~\ref{fig:spec_example_results} displays the absolute difference between these two quantities. What can be clearly seen is that a greater proportion of the plot is contaminated as the scales become wider, from left to right. Proportionally to the original quantities however, the magnitudes of the errors remain quite low.

The comparison conducted in this brief section carries over quite naturally when examining the difference between the expectation of the raw wavelet periodogram and its approximation, for which a bound was established in Theorem~\ref{thm:expectation}. Since the spectrum is approximated in a similar manner to the local autocovariance, we expect the comparison to yield similar results.

\section{Practical Considerations}
\label{sec:practicalconsiderations}

\subsection{Regularization Methods}
\begin{figure}
    \centering
    \includegraphics[width=1\linewidth]{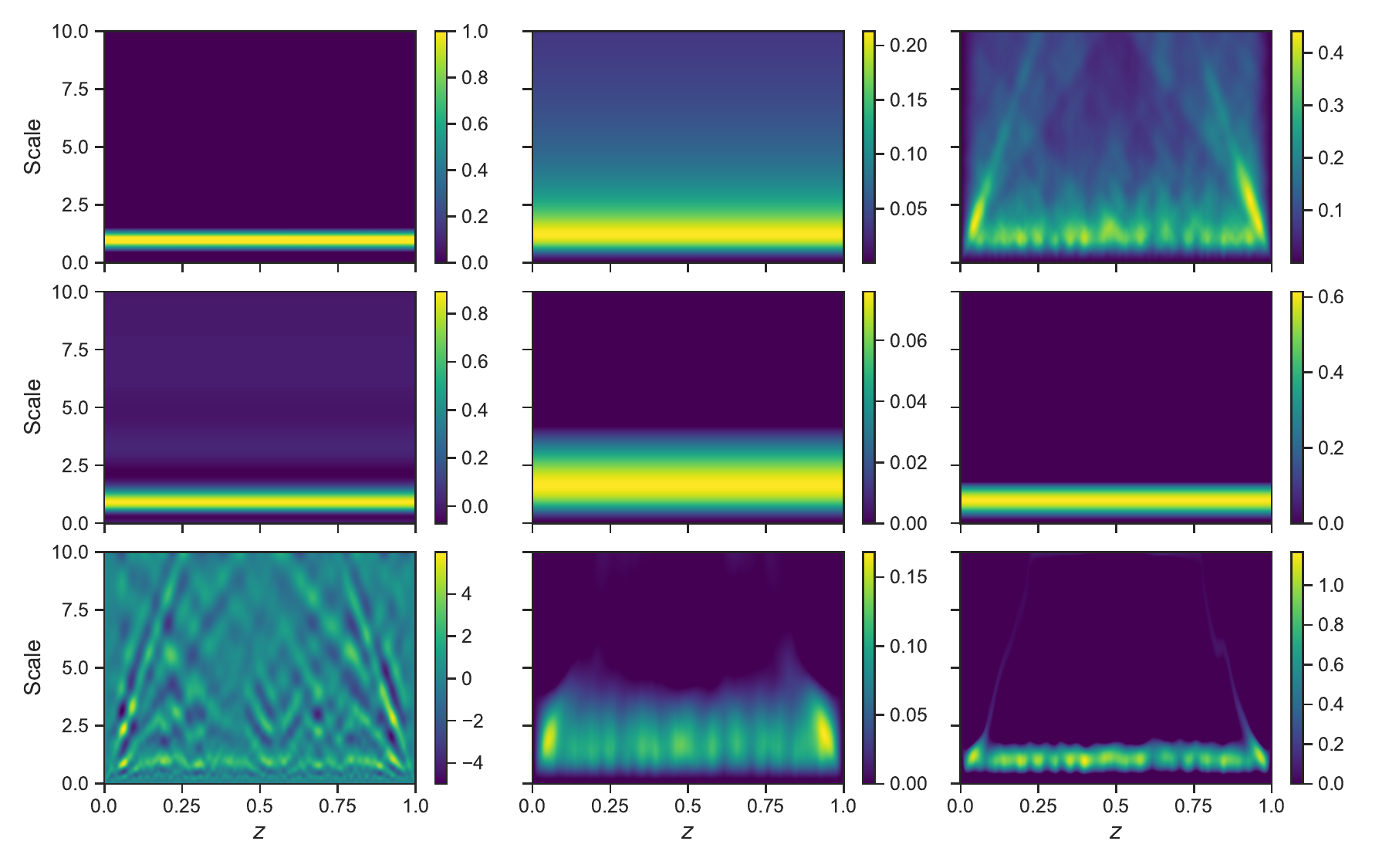}
    \caption{From left to right. Top row: true evolutionary wavelet spectrum, true raw wavelet periodogram $\beta$ and an estimate of the raw wavelet periodogram $\hat\beta$. Central row: evolutionary wavelet spectra computed from $\beta$ using Tikhonov regularization and \cite{daub_inv}'s iterative soft-thresholding algorithm for $N=100$ and $N=10000$. Bottom row: the same as the central row, but using $\hat\beta$ instead of $\beta$.}
    \label{fig:tikvsdaub}
\end{figure}

As discussed in Section~\ref{subsec:review3}, there are several regularization methods to choose from. We now provide some justification for our choice of the iterative soft-thresholding algorithm of~\cite{daub_inv} through a simple comparison with Tikhonov regularization. The spectrum used here is displayed in the top left panel of Figure~\ref{fig:tikvsdaub} with Haar wavelets as the underlying wavelet family. The top centre and right panels of Figure~\ref{fig:tikvsdaub} show the true expectation of the raw wavelet periodogram $\beta$ and a noisy estimate $\hat\beta$, respectively.

The central row of Figure~\ref{fig:tikvsdaub} displays the spectrum obtained using Tikhonov regularization on the left and using~\cite{daub_inv}'s iterative soft-thresholding algorithm for $N=100$ and $N=10000$ in the middle and on the right respectively, when we have access to the actual theoretical expectation $\beta$. Tikhonov regularization is the natural choice for inverting~\eqref{eq:invert}, both in terms of accuracy and computational efficiency, although it fails to preserve the non-negativity of the spectrum even in this setting.

However, access to the true $\beta$ is rarely available and as soon as we work with estimates, the performance of Tikhonov regularization rapidly deteriorates.  On the other hand,~\cite{daub_inv}'s iterative thresholding remains robust against noise contamination and preserves the non-negativity, as can be seen in the bottom row of Figure~\ref{fig:tikvsdaub}, which displays the spectra in the same order as the middle row but using $\hat\beta$ instead of $\beta$.

\subsection{Discretization Procedure}\label{sec:discretization}
When discretizing the various quantities of interest, we assume that we have a set of observations of $X_T(t)$ at times $0=t_0<t_1<\cdots<t_{n(T)-1}=T$, with $n(T)$ from Definition~\ref{def:clswp2}. The continuous wavelet transform can be computed for arbitrary scales $u$ and shifts $v$, without being restricted by the structure of the time series data. Hence, we opt to do so using a regularly spaced grid, consisting of $M_u\in\mathbb{N}$ scales and $M_v\in\mathbb{N}$ locations. We select scales $u_1 < u_2 < \cdots < u_{M_u} = J_\psi(T)$ and locations $0 = v_1 < v_2 < \cdots < v_{M_v} = T$, such that $u_i - u_{i-1} = \Delta u$ and $v_j - v_{j-1} = \Delta v$, for $i = 2,\ldots, M_u$, $j=2,\ldots, M_v$. In general, $M_u$ and $M_v$ grow linearly with $T$. The wavelet coefficients~\eqref{eq:coeffs} can then be approximated as
\begin{align*}
    d(u_i,v_j) \approx \frac{1}{2} \sum_{k=1}^{n(T)-1} \{ X_T(t_k)\psi(u_i, t_k-v_j) + X_T(t_{k-1})\psi(u_i,t_{k-1}-v_j)\}\Delta t_k,
\end{align*}
or, equivalently,
\begin{align*}
    d(u_i,v_j) \approx \frac{1}{2}\sum_{k=0}^{n(T)-1} X_T(t_k)\psi(u_i, t_k-v_j) \tilde\Delta t_k,
\end{align*}
where $\Delta t_k = t_k - t_{k-1}$ and $\tilde\Delta t_k = t_{k+1} - t_{k-1}$.

    When computing the continuous wavelet transform, or wavelet transforms in general, the question often arises as to which wavelet family to use. The bound in Theorem~\ref{thm:expectation} depends on the decay properties of the wavelet and the Lipschitz constant function of the process, which can be interpreted as a measure of non-stationarity. Hence, the choice of wavelet depends on the degree of non-stationarity exhibited by the process. For example, Shannon wavelets would be more appropriate for analyzing signals with limited non-stationarity, as they can leverage a significantly wider range of data to estimate the coefficients at each location. However, they may struggle with highly non-stationary data, for which the Haar wavelet is probably a better choice.

The corresponding discretized raw wavelet periodogram is then 
\begin{equation*}
    I(u_i, z_j) = \left|d(u_i,v_j)\right|^2,
\end{equation*}
where $z_j=v_j/T$. 

\begin{remark}
    It is not strictly necessary to restrict the discretization procedure to a regularly spaced grid. However, doing so for the locations $v$ proves particularly useful when applying smoothing techniques, as described in Section~\ref{sec:smoothing}.
\end{remark}

Throughout this work, we have mentioned that, in practice, we generally work with an estimate $\hat{\beta}$ of the expectation of the raw wavelet periodogram $\beta$. In most situations, we only have a single realization of the process, $X_T(t)$, and thus our best estimate would simply be our discretized raw wavelet periodogram $\hat{\beta}(u_i,z_j) = I(u_i,z_j)$. However, there are situations where we may have access to several realizations of the same process, known as replicated time series, see \cite{replicated} and \cite{replicated2} for example, and so can set 
\begin{equation} \label{eq:estimeq}
    \hat{\beta}(u_i,z_j) = R^{-1}\sum_{r=1}^R I_r(u_i,z_j),
\end{equation}
where $R$ is the number of replicates and $I_i(u,z)$ denotes the raw wavelet periodogram for the $i$-th realization. Here we again see the strength of the continuous wavelet transform, in that~\eqref{eq:estimeq} can be computed regardless of whether the specific time points of the individual time series match. The estimate $\beta$ can then be smoothed using the methods described in Section~\ref{sec:smoothing} of the main text.

Using the same grid, the inner product operator in~\eqref{eq:operatorT2} can be approximated as 
\begin{equation*}
    \{T_A S(\,\cdot\,,z_j)\}(u_i) \approx \Delta u \sum_{k=1}^{M_u} A(u_i, u_k) S(u_k, z_j)
\end{equation*}
and then used to obtain an estimate of the evolutionary wavelet spectrum $\hat S(u_i, z_j)$ from $\hat\beta(u_i, v_j)$ using \cite{daub_inv}'s iterative scheme in~\eqref{eq:spec_iter}. This can be written in matrix form as 
\begin{equation*}
    T_A (S_{z_j} ) = A S_{z_j},
\end{equation*}
where $S_{z_j} = \{S(u_1, z_j), \ldots, S(u_{M_u})\}$ is a vector of size $M_u$ and $A$ is a matrix with entries $A_{i,k} = \Delta u A(u_i, u_k)$, for $i,k = 1, \ldots, M_u$.
From $\hat S(u_i, z_j)$ we can also compute an estimate of the local autocovariance using
\begin{equation}
    \hat c(z_j,\tau) \approx \Delta u \sum_{i=1}^{M_u} \hat S(u_i, z_j) \Psi(u_i, \tau),
\end{equation}
for appropriately chosen lags $\tau$, which again can also be written in matrix form.

\subsection{Implementation of the Iterative Soft-Thresholding Algorithm}

In this section we investigate the convergence behaviour of \cite{daub_inv}'s iterative soft-thresholding algorithm and discuss its implementation in practice. Empirically, we observe that convergence of the evolutionary wavelet spectrum at finer scales typically requires more time. The reason for this is two-fold: 1.  wavelet coefficients at finer scales are naturally smaller than those at coarser scales, which means it takes longer for this information to appear in the spectrum; 2. in order to use \cite{daub_inv}'s iterative scheme, we need to normalize the inner product operator by its norm, replacing $T_A$ with $T_A / ||T_A||$ and $\hat\beta$ with $\hat\beta / ||T_A||$ in~\eqref{eq:spec_iter}. Formally, the iterative process becomes
\begin{equation}\label{eq:norm_iter_scheme}
    \hat S^n(\,\cdot\,,z) = \mathcal{S}_{\mu}\left[\hat S^{n-1}(\,\cdot\,,z) + \frac{T_A}{||T_A||}\left\{\frac{\hat\beta(\,\cdot\,,z)}{||T_A||} - \frac{T_A}{||T_A||} \hat S^{n-1}(\,\cdot\,,z)\right\}\right], \quad n = 1,\ldots, N,
\end{equation}
for $z\in(0,1)$, with $\mathcal{S}_{\mu}$ from equation~\eqref{eq:spec_shrink}. This is equivalent to solving the integral equation
\begin{equation*}
    \frac{\beta(u, v)}{||T_A||} = \left\{ \frac{T_A}{||T_A||} S(\,\cdot\,,\, v)\right\}(u),
\end{equation*}
which has the same solution as~\eqref{eq: approx_spec}, with the added benefit of ensuring that the iterative scheme remains stable and does not explode. However, since the norm $||T_A||$ increases  with the width of the largest scale, $u_{M_u} = J_\psi(T)$,  the size of the iterative updates at finer scales decreases, resulting in slower convergence at those scales.

\begin{figure}
    \centering
    \includegraphics[width=0.65\linewidth]{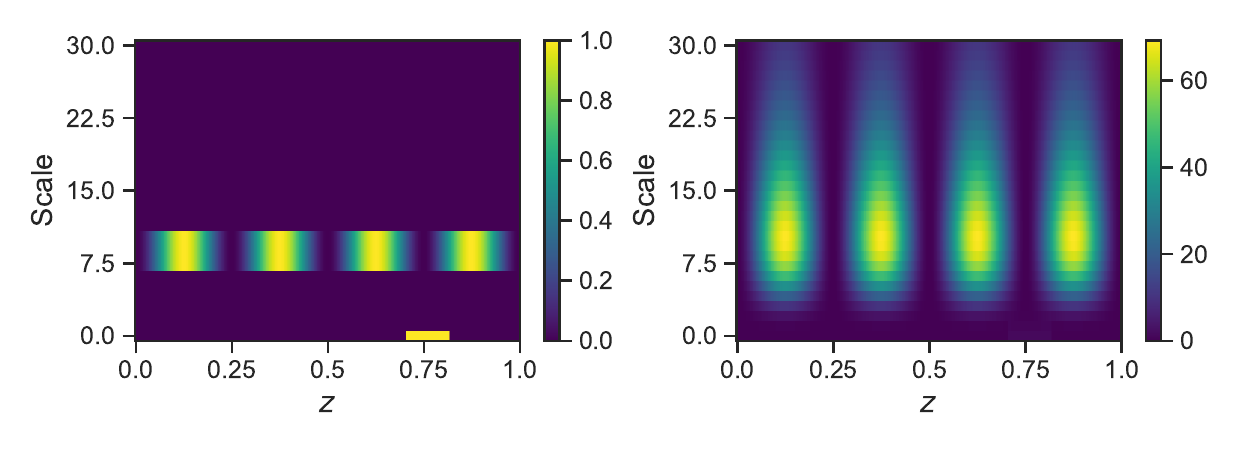}
    \caption{Left: square-sine and burst evolutionary wavelet spectrum. Right: Corresponding raw wavelet periodogram .}
    \label{fig:sine_burst}
\end{figure}

We demonstrate this behaviour on an example. The underlying spectrum is given by 
\begin{equation*}
    S(u, z) = \left\{\begin{array}{lll}
    \sin^2(4\pi z),     &&  u\in (8, 12), \;z\in(0,1),\\
    1,     && u\in(0.75, 1.25), \;z\in(0.7, 0.825), \\
    0,     && \text{otherwise},
    \end{array}\right.
\end{equation*}
and can be viewed as a continuous analogue to the square sine and burst spectrum that appears in Chapter~5 of~\cite{nasonbook}. The underlying wavelets used in this example are Ricker wavelets. The evolutionary wavelet spectrum and corresponding raw wavelet periodogram are displayed in Figure~\ref{fig:sine_burst}. The high frequency information is clearly visible in the evolutionary wavelet spectrum around scale 1, where it has the same spectral power as the square sine wave at scales 8 to 12. However, in the raw wavelet periodogram, the power at scale 1 is significantly weaker, making it invisible in the figure, despite its presence. This aligns with the fact that wavelet coefficients at finer scales are inherently smaller.

\begin{figure}
    \centering
    \includegraphics[width=1\linewidth]{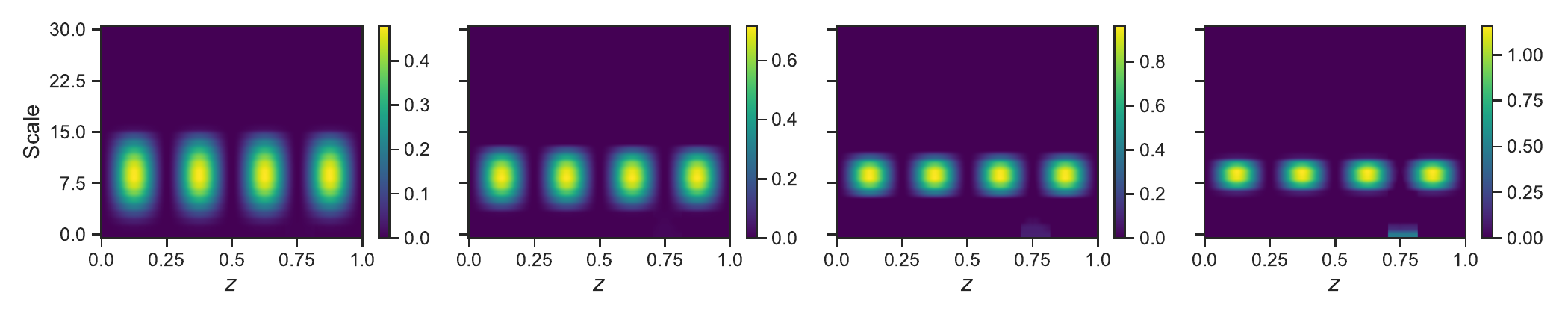}
    \caption{Reconstructions of the square-sine and burst spectrum from the true raw wavelet periodogram. From left to right: $N=100$, $1000$, $10000$ and $100000$.}
    \label{fig:naive_run}
\end{figure}

Figure~\ref{fig:naive_run} displays the results of running~\cite{daub_inv}'s iterative soft-thresholding algorithm for $N=100$, $1000$, $10000$ and $100000$ iterations, starting with the true raw wavelet periodogram $\beta$. Even for $N=100000$, the high frequency information remains barely visible in the reconstructed spectrum. This is primarily due to the significantly smaller step size at finer scales, due to the normalization of the inner product operator in the iterative scheme.

With only $N=100$ iterations the upper part of the spectrum has already converged to its true value of zero. Similarly, the square-sine region also seems to have practically converged after $N=10000$ iterations. These observations suggest an alternative approach: instead of updating the entire spectrum with each iteration, we could instead update the spectrum up to some scale $u$, renormalising the inner product operator accordingly. This would result in larger step sizes at finer scales and faster convergence.

\begin{figure}
    \centering
    \includegraphics[width=0.85\linewidth]{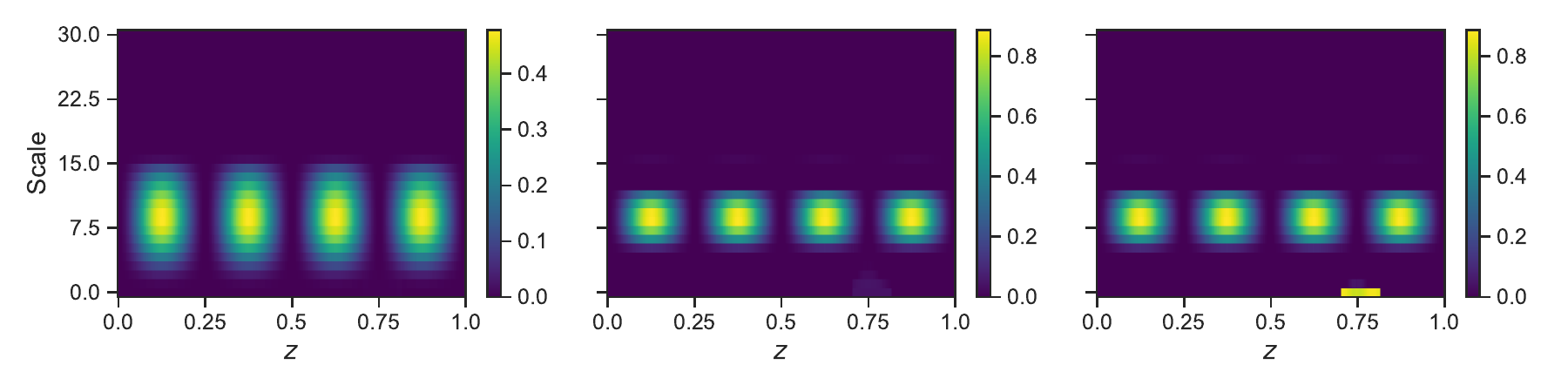}
    \caption{Reconstructions of the square-sine and burst spectrum from the true raw wavelet periodogram. Left: $N_1=100$ iterations of the thresholding algorithm applied to the full spectrum. Centre: a further $N_2=250$ iterations applied only up to scale $u=15$. Right: a final $N_3=100$ iterations applied only up to scale $u=4$.}
    \label{fig:smart_run}
\end{figure}

Figure~\ref{fig:smart_run} displays the results of initially running the iterative soft-thresholding algorithm for $N_1=100$ iterations, then applying this up to scale $u=15$ for another $N_2=250$ iterations, and then up to scale $u=4$ for a final $N_3=100$ iterations. This results in significantly faster convergence, requiring only $N=450$ in total compared to $N=100000$ when applying the scheme to the entire spectrum, while also achieving better results.

Using the notation of the previous section, the iterative scheme~\eqref{eq:norm_iter_scheme} in practice is 
\begin{align*}
    \hat S^n_{z_j} &= \mathcal{S}_{\mu}\left\{\hat S^{n-1}_{z_j} + \frac{A}{||A||}\left(\frac{\hat\beta_{z_j}}{||A||} - \frac{A}{||A||} \hat S^{n}_{z_j}\right)\right\}, \\
    & = \mathcal{S}_{\mu}\left(\hat S^{n-1}_{z_j} + \frac{A\hat\beta_{z_j}}{||A||^2} - \frac{A^2}{||A||^2} \hat S^{n}_{z_j}\right), \quad n = 0,\ldots, N,
\end{align*}
for $j=1,\ldots, M_v$, where $||A||$ denotes the norm of $A$, $\hat\beta_{z_j} = \{\hat\beta(u_1, z_j),\ldots,\hat\beta(u_{M_u}, z_j)$ is a vector of length $M_u$ and $\mathcal{S}_{\mu}$ is the usual shrink operator. The restriction of the iterative soft-thresholding algorithm to scales up to $u_i$, for some $1 < i <M_u$, is 
\begin{equation*}
    \left(\hat S^n_{z_j}\right)_{1:i}  = \mathcal{S}_{\mu}\left\{\left(\hat S^{n-1}_{z_j}\right)_{1:i} + \frac{A_{1:i}\hat\beta_{z_j}}{||A_{1:i,1:i}||^2} - \frac{\left(A^2\right)_{1:i}}{||A_{1:i,1:i}||^2} \hat S^{n-1}_{z_j}\right\},
\end{equation*}
where  $S_{1:i}$ corresponds to taking the first $i$ entries of the vector $S$, $A_{1:i}$ corresponds to taking the first $i$ rows of $A$ and $A_{1:i, 1:i}$ corresponds to taking the first $i$ rows and columns of $A$. Hence, $A_{1:i}$ is of dimension $i\times M_u$ instead of $M_u\times M_u$, allowing each step to make use of the entire periodogram estimate $\hat\beta_{z_j}$ and previous iterate $\hat S^{n-1}_{z_j}$, while updating the spectrum $\hat S^{n}_{z_j}$ only up to scale $u_i$ and therefore only applying the normalization up to this scale in the form $||A_{1:i,1:i}||$. 

The number of iterations $N$ depends on the specific context. Typically, low-frequency information becomes well localized with only 
$N=100$ to $1000$ iterations. In practice, a stopping criterion can be implemented to interrupt the iterative algorithm once the difference between consecutive iterates falls below a predefined threshold. Alternatively, the algorithm could use such a criterion to progressively reduce the number of scales it is applied to as the larger scales converge.

We now turn our attention to the problem of setting the regularisation parameter(s) $\mu$ in practice. Chapter~4 of~\cite{inv_book} provides some guidance on this topic, discussing methods such as the L-curve, which plots the norm of the regularized solution against the norm of its residuals for different values of $\mu$. However, such methods require running the regularization algorithm multiple times. Alternatively, $\mu$ can be chosen based on the noise level, as specified by the conditions of~Theorem~\ref{thm:conv_cond}. The noise can be estimated using the mean absolute deviation of the finest scale wavelet coefficients as proposed by~\cite{DJ94} for wavelet shrinkage or~\cite{mad}. This approach can also be adapted to be scale and location dependent, following the suggestions made in Section 4.4 of the main text.

\section{Non-Stationary Haar Moving Average Example}
\subsection{Evolutionary Wavelet Spectrum}
\begin{figure}
    \centering
    \includegraphics[width=0.9\linewidth]{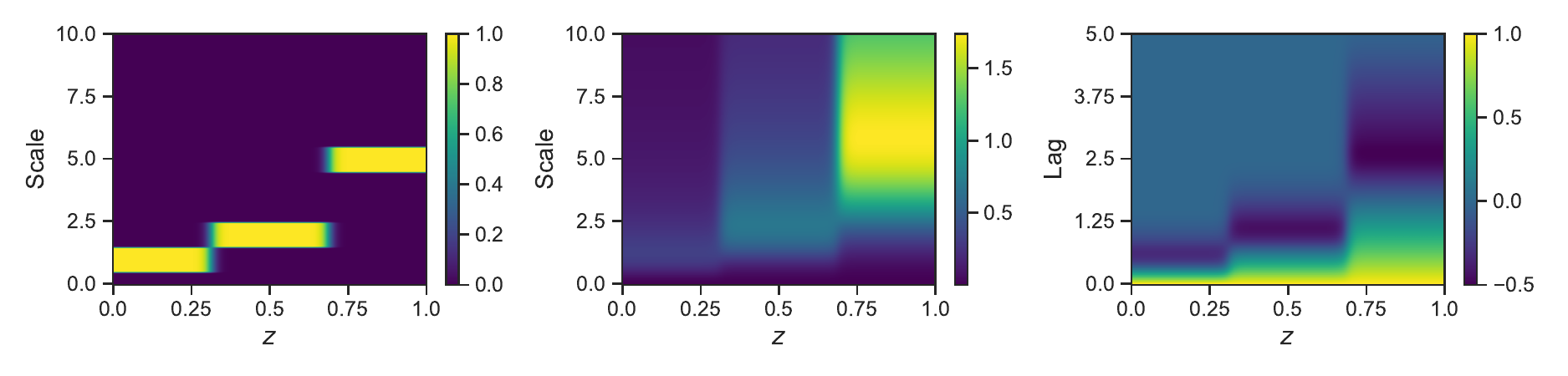}
   \caption{Haar moving average theoretical quantities. Left: evolutionary wavelet spectrum. Centre: raw wavelet periodogram. Right: local autocovariance.}
    \label{fig:haar2}
\end{figure}

The full functional form of the evolutionary wavelet spectrum  on the non-stationary Haar moving average process used in Section~\ref{sec:haar} is
\begin{equation*}
S(u, z) = \left\{\begin{array}{lll}
1,  && \text{for } z \in \left(0, 1/4\right] \text{ and }  u\in(1/2, 3/2), \\
f(15-48z),  && \text{for } z \in \left(1/4, 3/8\right) \text{ and } u\in(1/2, 3/2), \\
f(48z-15), && \text{for } z \in \left(1/4, 3/8\right) \text{ and }  u\in(3/2, 5/2), \\
1, && \text{for } z \in \left[3/8, 5/8\right] \text{ and } u\in(3/2, 5/2),\\
f(33-48z), \quad  && \text{for } z \in \left(5/8, 3/4\right) \text{ and } u\in(3/2, 5/2), \\
f(48z-33), \quad  && \text{for } z \in \left(5/8, 3/4\right) \text{ and } u\in(9/2, 11/2), \\
1, && \text{for } z \in \left[3/4, 1\right) \text{ and } u\in(9/2, 11/2),\\
0, && \text{otherwise},
\end{array}\right.
\end{equation*}
where $f(z) = \tanh(z)/2 + 1/2.$ This is displayed in Figure~\ref{fig:haar2}, alongside the expectation of the raw wavelet periodogram $\beta$ and the local autocovariance.

\subsection{Other Estimates}
\begin{figure}
    \centering
    \includegraphics[width=1\linewidth]{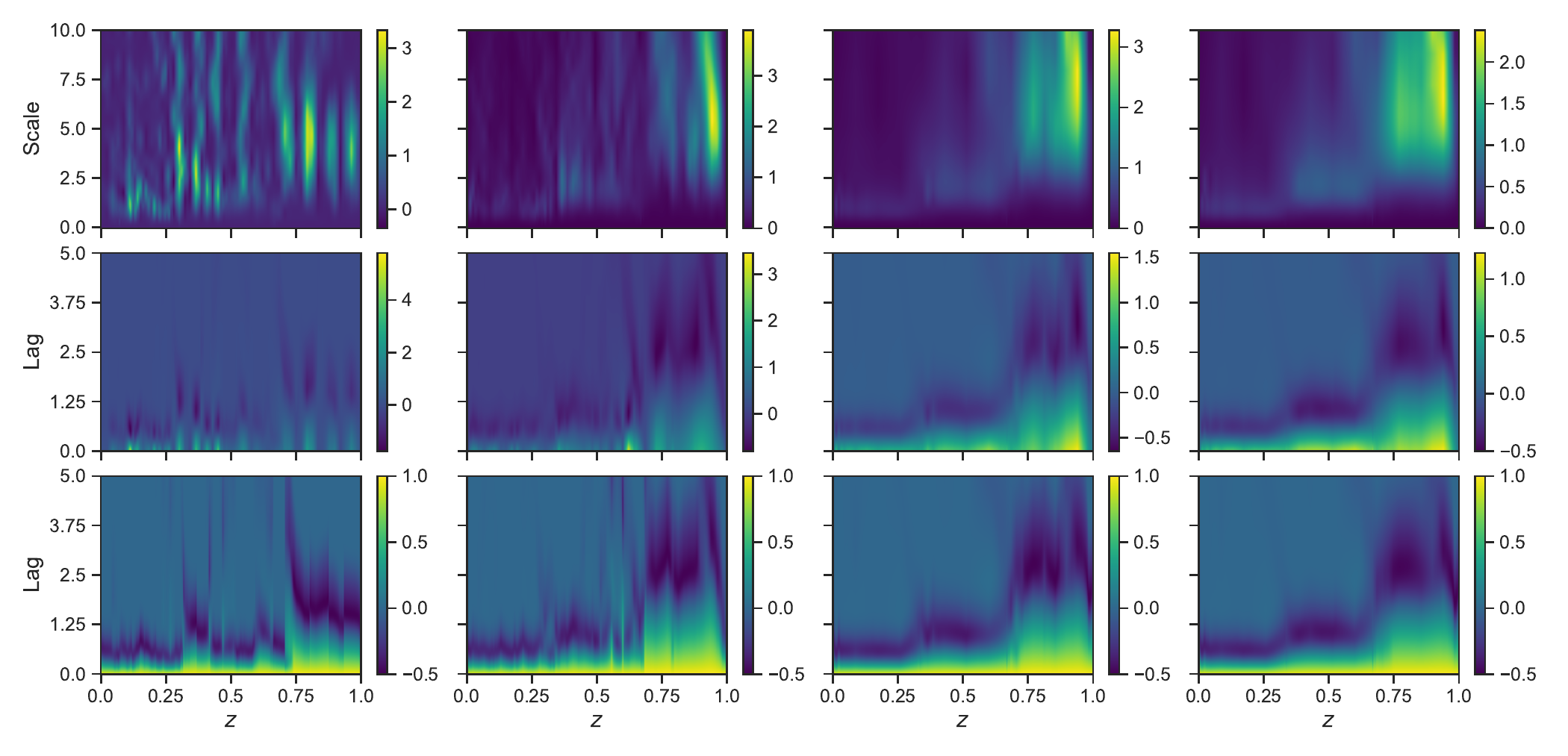}
    \caption{Haar moving average estimates computed using $R=1, 10, 100$ and $1000$ realizations, from left to right. Top row: raw wavelet periodogram. Central row: local autocovariance. Bottom row: local autocorrelation.}
    \label{fig:estims_sm}
\end{figure}

The top row of Figure~\ref{fig:estims_sm} displays estimates $\hat\beta$ of $\beta$ computed using equation~\eqref{eq:estimeq} for $R=1, 10, 100$ and $1000$ realizations. The other two rows of Figure~\ref{fig:estims_sm} show the corresponding estimates of the local autocovariance and autocorrelation computed from estimates of the continuous time evolutionary wavelet spectrum. 

As expected, the quality of the estimates improves as the number of samples used increases. In all cases, the three separate frequency regions can be distinguished from one another. The local autocorrelation resembles its theoretical counterpart far more than the local autocovariance, especially when a small number of samples is used.

\subsection{Missing Data}

\begin{figure}
    \centering
    \includegraphics[width=1\linewidth]{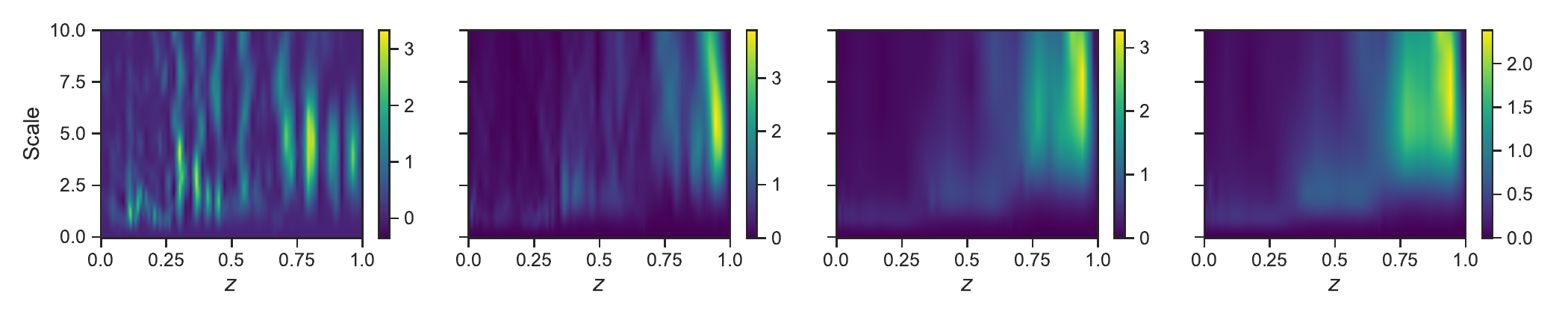}
    \caption{Haar moving average raw wavelet periodogram estimates with 25\% missing data. From left to right: estimates computed using $R=1, 10, 100$ and $1000$ realizations.}
    \label{fig:haar_missing}
\end{figure}

We conduct the same exercise as in Section~\ref{sec:haar} in the main
paper, but with missing data. We use the same set of non-stationary Haar moving average realizations, treating 25\% of observations as missing each time. We then compute an estimate $\hat\beta$ of the raw wavelet periodogram $\beta$ using equation~\eqref{eq:estimeq} for $R=1,10,100$ and $1000$.  Figure~\ref{fig:haar_missing} displays the results. The raw wavelet periodogram estimates look almost indistinguishable from those computed using the complete time series. Furthermore, these are defined on a regularly spaced grid, despite the data containing missing information. This facilitates the application of smoothing techniques or other tools, which are often developed for regularly spaced data.

\begin{remark}
    We kept the locations of the missing observations fixed throughout realizations. An advantage of using the continuous wavelet transform is that it can be computed at any scale and location, so the time points of missing values do not need to be identical between replicates. In fact, if the missing observations varied in position across realizations, our estimates would likely become more robust for increasing $R$, as the effects of the missingness would be averaged out.
\end{remark}

\subsection{Empirical Exploration of the Convergence of Estimates with Increasing Data}
The purpose of this final section is to examine the difference between the estimates of the raw wavelet periodogram and their theoretical counterpart as both $T$ and $R$ increase. To gain some more intuition as to what $T$ represents, let us assume that we have a finite realization of the stationary Haar \textsc{ma}$(\alpha)$ process
\begin{equation*}
    X(t)  = \alpha^{-1/2}\left[\left(B_t - B_{t-\alpha/2} \right) - \left(B_{t-\alpha/2} - B_{t-\alpha}\right)\right]
\end{equation*}
for $t\in(0, T)$, where $B_t$ is a standard Brownian motion and $\alpha$ is the order of the process. The rescaled-time stationary Haar \textsc{ma}$(\alpha)$ process is then
\begin{equation*}
    X(z) = \alpha^{-1/2}\left[\left(B_z - B_{z-\alpha/2} \right) - \left(B_{z-\alpha/2} - B_{z-\alpha}\right)\right]
\end{equation*}
for $z = t/T\in(0,1)$. Hence, as $T$ increases, we obtain increasing amounts of information about the local structure of the process. Since the standard Haar moving average process is stationary, this makes no difference to the asymptotic analysis as the local structure is the same across all time. However, for the non-stationary Haar moving average process, whose spectrum is displayed in Figure~\ref{fig:haar2}, this means that as $T$ increases, we obtain more and more information about the properties of each of the three distinct regions in the process.

\begin{figure}
    \centering
    \includegraphics[width=1\linewidth]{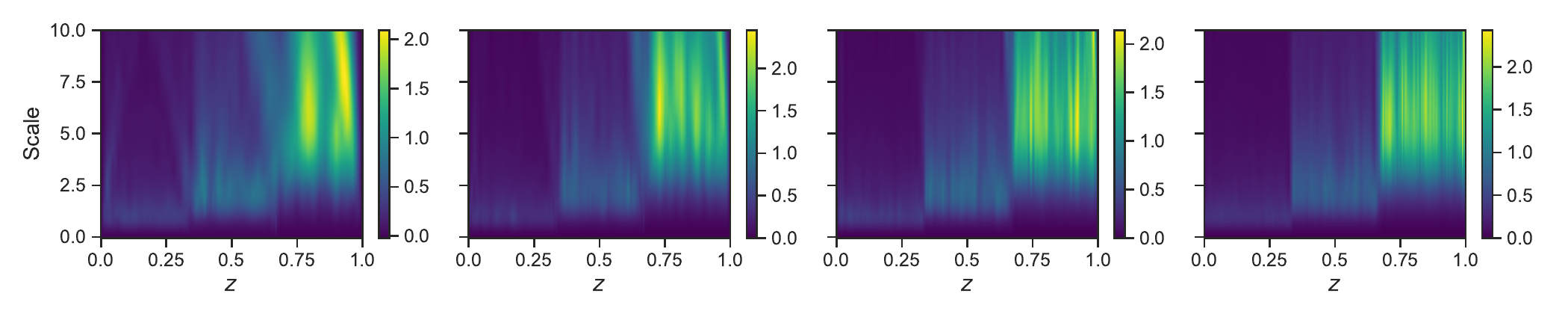}
    \caption{Haar moving average raw wavelet periodogram estimates using $R=100$ realizations. From left to right: amount of data $T=30$, $60$, $150$ and $300$.}
    \label{fig:increasingT}
\end{figure}

Figure~\ref{fig:increasingT} displays estimates $\hat\beta$ of the raw wavelet periodogram computed using $R=100$ realizations for increasing amounts of data $T$, namely $T=30$, $60$, $150$ and $300$. These values were selected 
so that the length of each region is $10$, $20$, $50$ and $100$. As expected, the estimates improve as $T$ increases, in accordance with the statement of Theorem~\ref{thm:expectation2}. Moreover, the boundary effects reduce with increasing data, further improving the quality of the estimate. This is due to the fact that, as $T$ grows, the width of the wavelet corresponding to the scale of that spectral information decreases as a proportion of the re-scaled time domain $(0, 1)$, thereby reducing the magnitude of the boundary effects on the estimate. 

\begin{figure}
    \centering
    \includegraphics[width=0.9\linewidth]{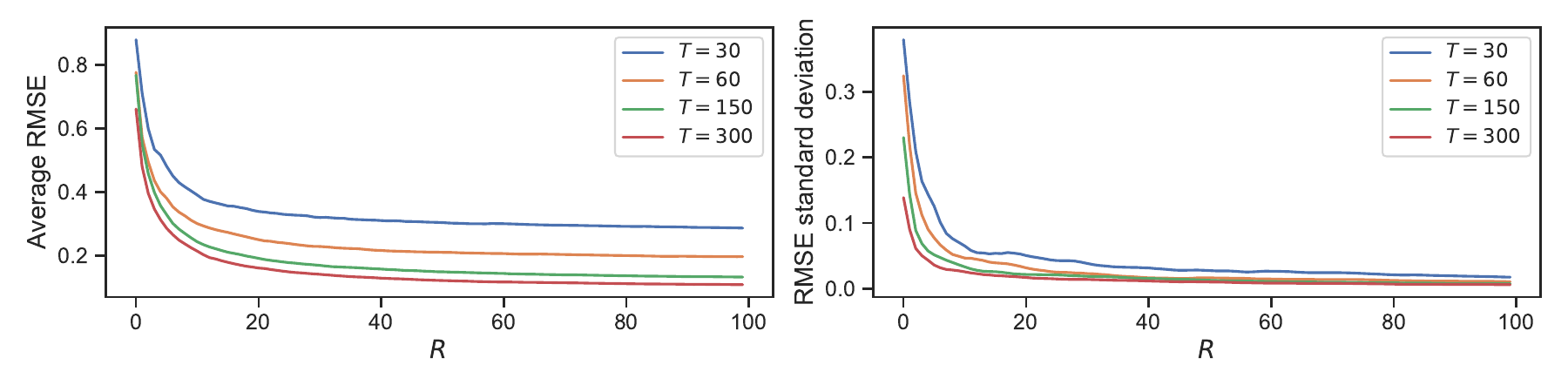}
    \caption{Left: average root mean squared error for raw wavelet periodograms estimated using increasing numbers of realizations $R$ and different amounts of data $T$. Corresponding standard deviations in the root mean squared errors across the $100$ runs.}
    \label{fig:rmse}
\end{figure}

Figure~\ref{fig:rmse} displays the root mean squared error between the estimates of the raw wavelet periodogram and their theoretical counterpart from $R = 1$ to $R=100$ and for different amounts of data $T$, averaged across $100$ runs. The right panel displays the corresponding standard deviations of the root mean squared errors across these runs. The root mean squared error here is 
\begin{equation*}
    \text{RMSE} = \left[\sum_{i=0}^{M_u}\sum_{j=0}^{M_v} \left\{\hat\beta(u_i, z_j) - \beta(u_i, z_j)\right\}^2\right]^{1/2},
\end{equation*}
using the same notation as Section~\ref{sec:discretization}. To obtain robust estimates of the root mean squared error, we averaged the results across $100$ runs. As expected, the root mean squared error decreases as the number of realizations $R$ used in equation~\eqref{eq:estimeq} increases. Moreover, the root mean squared error also decreases as the amount of data $T$ increases, again aligning with the results of Theorem~\ref{thm:expectation2}. Similarly, the standard deviation in the root mean squared errors across runs decreases with increasing $R$ and $T$, indicating that estimates become more robust with more data, as expected. 

\begin{figure}
    \centering
    \includegraphics[width=0.9\linewidth]{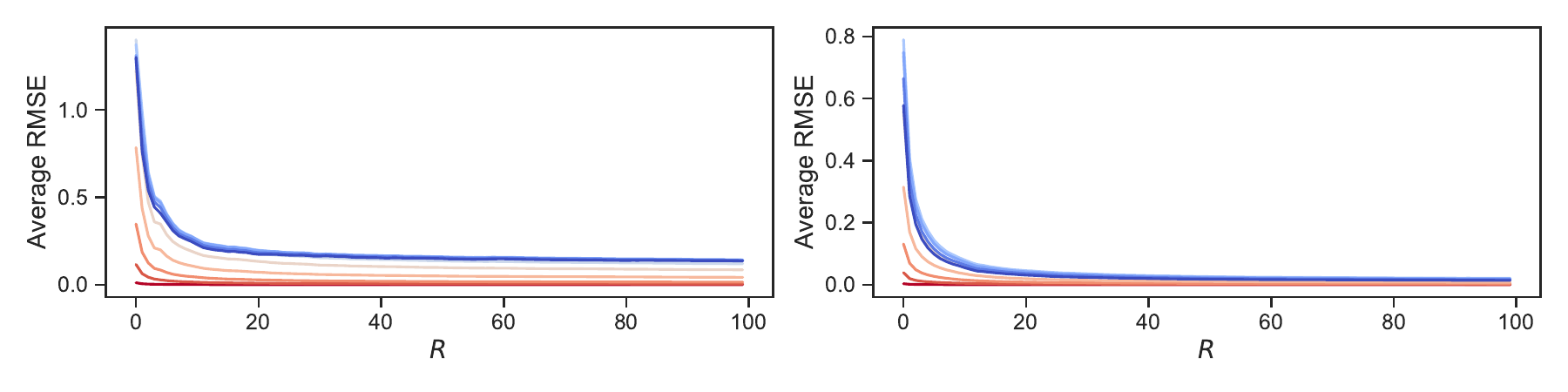}
    \caption{Average root mean squared errors by scale. The line colours gradually change from red to blue as the scales increase from finer to coarser, covering scales up to a width of $10$. Left: $T=30$. Right $T=300$.}
    \label{fig:per_scale_rmse}
\end{figure}

Figure~\ref{fig:per_scale_rmse} displays the average root mean squared errors across scales for $T=30$ and $T=300$. The errors are larger at the coarser scales, however, the highest errors occur between scales $5$ to $8$, which corresponds to the region with the strongest spectral information in the periodogram, as can be seen in the central panel of Figure~\ref{fig:haar2}. This is consistent with the standard notion that spectral information at coarser scales is inherently harder to  estimate.

\section{Proofs}
\subsection{Proof of Proposition \ref{prop:autocov}}
The autocovariance of $X(t)$ can be written as
\begin{align*}
    c_X(t,\tau) &= E\{X(t)X(t+\tau)\}\\
    &= \int_0^\infty\int_{-\infty}^\infty \int_0^\infty\int_{-\infty}^\infty W(u,v)W(x,y)\psi(u,v-t)\psi(x,y-t-\tau)du dv\\
    & = \int_0^\infty\int_{-\infty}^\infty |W(u,v)|^2\psi(u,v-t)\psi(u,v-t-\tau)du dv\\
    & = \int_0^\infty\int_{-\infty}^\infty S(u,v+t)\psi(u,v)\psi(u,v-\tau)du dv,
\end{align*}
where we used the fact that $E[L(du,dv)L(dx,dy)] = \delta(u-x)\delta(v-y)du dv dx dy$ and $\delta$ denotes the Dirac delta function. Since
    \begin{equation*}
        c(t,\tau) = \int_0^\infty S(u,t)\Psi(u,\tau)du = \int_0^\infty S(u,t)\int_{-\infty}^{\infty}\psi(u,v)\psi(u,v-\tau) du dv,
    \end{equation*}
we have
\begin{align*}
    |c_X(t, \tau) - c(t,\tau)| &=\left|\int_0^\infty\int_{-\infty}^\infty \left(S(u,v+t) - S(u, t)\right)\psi(u,v)\psi(u,v-\tau)du dv\right|\\
    &\leq \left| \int_0^\infty\int_{-\infty}^\infty \left|S(u,v+t) - S(u, t)\right|\psi(u,v)\psi(u,v-\tau)du dv\right|\\\numberthis \label{eq:p1a}
    &   \leq 2M \left|\int_0^\infty  K(u)\int_{-\infty}^\infty |v|\psi(u,v)\psi(u,v-\tau)du dv\right|\\\numberthis \label{eq:p1b}
    &\leq 2 M \omega \left|\int_0^\infty  u K(u) \int_{-\infty}^\infty  \psi(u,v)\psi(u,v-\tau)du dv\right|\\
    & = 2 M \omega  \left|\int_0^\infty u K(u) \Psi(u,\tau) du\right|\\ \numberthis\label{eq:p1d}
    & \leq   \gamma\int_0^\infty u K(u)du.
\end{align*}
In \eqref{eq:p1a}, we used that
\begin{align*}
    \left|S(u,v+t) - S(u, t)\right| &= \left||W(u,v+t)|^2 - |W(u, t)|^2\right|\\ &= \left|W(u,v+t) + W(u, t)\right| \left|W(u,v+t) - W(u, t)\right| \\
    &\leq 2 M K(u) |v|,
\end{align*}
for some $M\in\mathbb R^+$ such that $|W(u,v)| \leq M$ for all $u\in\mathbb R^+$ and $v\in\mathbb R$.

The existence of such an $M$ follows from the Lipshitz continuity of $W(u, \,\cdot\,)$ and assumption (ii) in definition~\ref{def:clswp}. The inequality in \eqref{eq:p1b} follows from the fact that the integral grows linearly in $u$, where $\omega\in\mathbb R^+$ is a constant that depends on the width and decay rate of $\psi$. Since $|\Psi(u,\tau)| \leq 1$ for all $u\in\mathbb R^+$ and $\tau\in\mathbb R$, and setting $\gamma = 2M\omega$, we obtain the final inequality in \eqref{eq:p1d}.

\subsection{Proof of Proposition \ref{prop:properties1}}
(i) Follows from the definition.

(ii) $\Psi(u,0)= ||\psi||^2 = 1$.

(iii) We have that
    \begin{align*}
         \Psi(u, \tau) &= \int_{-\infty}^{\infty}\psi(u,v)\psi(u,v-\tau)dv\\
         &= \frac{1}{u}\int_{-\infty}^{\infty}\psi\left(\frac{v}{u}\right)\psi\left(\frac{v-\tau}{u}\right)dv \\
         &=\int_{-\infty}^{\infty}\psi\left(y\right)\psi\left(y-\frac{\tau}{u}\right)dy \\
         & = \Psi\left(\frac{\tau}{u}\right), 
    \end{align*}
    where we made the substitution $v \mapsto uy$.
   
(iv) We have that 
    \begin{align*}
        \Psi(\tau) = \int_{-\infty}^{\infty}\psi(v)\psi(v-\tau)dv = (\psi *\psi)(\tau),
    \end{align*}
    where $*$ denotes the convolution operator. Hence
    \begin{align*}
        \hat\Psi(\omega) = \mathcal F\left(\Psi\right)(\omega)= |\hat\psi(\omega)|^2.
    \end{align*}
    Since, $\Psi(u, \tau) = \Psi\left(\tau/u\right)$, we have 
    \begin{align*}
        \hat\Psi(u,\omega) = \mathcal{F}\left\{\Psi\left(\frac{\cdot}{u}\right)\right\} (\omega)= u\hat\Psi(u\omega) =  u |\hat\psi(u\omega)|^2.
    \end{align*}

(v) We have 
    \begin{align*}
        \int_{-\infty}^\infty\Psi(u,\tau) d\tau &= \int_{-\infty}^\infty\int_{-\infty}^{\infty}\psi(u,v)\psi(u,v-\tau)dv d
        \tau\\
        &= \int_{-\infty}^\infty \psi(u,v)\int_{-\infty}^{\infty}\psi(u,v-\tau) d
        \tau dv\\
        & = 0,
    \end{align*}
    since the inner integral is zero.

\subsection{Proof of Proposition \ref{prop:properties2}}
Continuity follows from the continuity of the autocovariance wavelets. Symmetry is clear from the definition. We now show that $A$ is non-negative definite. Let $c_i,c_j\in\mathbb R$ and $x_i,x_j\in\mathbb R^+$ for $i,j=1,\ldots,n.$ Then
\begin{align*}
    \sum_{i=1}^n\sum_{j=1}^n c_ic_jA(x_i,x_j) &= \int_{-\infty}^\infty \sum_{i=1}^n\sum_{j=1}^n c_ic_j\Psi(x_i,y)\Psi(x_j,y)dy\\
    & =  \int_{-\infty}^\infty \left\{\sum_{i=1}^n  c_i\Psi(x_i,y)\right\}^2dy\\
    &\geq 0.
\end{align*}
For the final property, we have
\begin{align*}
    A(b x, b y) &= \int_{-\infty}^\infty \Psi(b x, \tau) \Psi(b y,\tau)  d\tau\\
    &= \int_{-\infty}^\infty \Psi\left( \frac{\tau}{b x}\right) \Psi\left(\frac{\tau}{b y}\right)  d\tau\\
    &= \int_{-\infty}^\infty \Psi\left(x, \frac{\tau}{b}\right) \Psi\left(y,\frac{\tau}{b}\right)  d\tau\\
    &=b \int_{-\infty}^\infty \Psi( x, t) \Psi( y,t)  dt\\
    &=b A(x,y),
\end{align*}
using property (iii) from Proposition~\ref{prop:properties1}.

\subsection{Proof of Theorem \ref{thm:expectation}}
We have that 
\begin{align*}
    \beta(u, v) &=  E\left\{I\left(u, v\right)\right\} = E\left[\left\{\int_{-\infty}^\infty X(t) \psi(u,t-v) dt\right\}^2\right]\\
    &=\int_{-\infty}^\infty \int_{-\infty}^\infty  E \{X(t)X(s)\}\psi(u,t-v)\psi(u,s-v)dt ds
\end{align*}
Now,
\begin{equation*}
    E \{X(t)X(s)\} = \int_0^\infty\int_{-\infty}^\infty |W(x,y)|^2\psi(x,y-t)\psi(x,y-s)dx dy,
\end{equation*}
where we used the fact that $ E\{L(dx,dy)L(d\eta,d\xi)\} = \delta(\eta-x)\delta(\xi-y)d\eta d\xi dx dy$ and $\delta$ denotes the Dirac delta function. Hence, 
\begin{align*}
    \beta(u, v) &=\int_{-\infty}^\infty \int_{-\infty}^\infty \int_0^\infty\int_{-\infty}^\infty S(x,y)\psi(x,y-t)\psi(x,y-s)\psi(u,t-v)\psi(u,s-v) dy dx dt ds\\  \numberthis\label{eq:thm1p2}
    & = \int_{-\infty}^\infty \int_{-\infty}^\infty \int_0^\infty\int_{-\infty}^\infty S(x,y+v)\psi(x,y-t)\psi(x,y-s)\psi(u,t)\psi(u,s) dy dx dt ds\\
    & =  \int_0^\infty\int_{-\infty}^\infty S(x,y+v)\left\{\int_{-\infty}^\infty\psi(x,y-t)\psi(u,t) dt\right\}^2 dy dx,
\end{align*}
where in \eqref{eq:thm1p2} we made the substitution $y\mapsto y+v$ followed by $t-v\mapsto t$ and $s-v\mapsto s$.
Now, 
\begin{align*}
    \int_0^\infty S\left(x, v\right)A(u,x)  dx &= \int_0^\infty S\left(x, v\right)\int_{-\infty}^\infty \Psi(x, \eta)\Psi(u, \eta) d\eta dx\\
    &= \int_0^\infty S\left(x,v\right)\int_{-\infty}^\infty  \int_{-\infty}^\infty \Psi(x, \eta)\psi(u,\eta+\gamma)\psi(u,\gamma) d\eta d\gamma dx\\ \numberthis\label{eq:thm1p3}
    & = \int_0^\infty S\left(x,v\right)\int_{-\infty}^\infty  \int_{-\infty}^\infty \Psi(x, t-s)\psi(u,t)\psi(u,s) dt ds dx,
\end{align*}
where in \eqref{eq:thm1p3} we made the substitutions $\eta\mapsto t-s$ and $\gamma\mapsto s$. Since 
\begin{align*}
    \Psi(x, t-s)= \int_{-\infty}^\infty\psi(x,r)\psi(x,r + t-s) dr =\int_{-\infty}^\infty\psi(x,y-t)\psi(x,y-s) dy,
\end{align*}
with $r\mapsto y-t$, we have 
\begin{align*}
    \int_0^\infty S\left(x, v\right) & A(u,x)  dx\\
    &= \int_{-\infty}^\infty \int_{-\infty}^\infty \int_0^\infty\int_{-\infty}^\infty S(x,v)\psi(x,y-t)\psi(x,y-s)\psi(u,t)\psi(u,s) dy dx dt ds\\
    & =  \int_0^\infty\int_{-\infty}^\infty S(x,v)\left\{\int_{-\infty}^\infty\psi(x,y-t)\psi(u,t) dt\right\}^2 dy dx.
\end{align*}
Using the fact that 
\begin{equation*}
    \left|S(x,y+v) - S(x,v)\right| \leq 2M K(x)  |y|,
\end{equation*}
as shown in the proof of Proposition~\ref{prop:autocov}, and following similar steps,
we obtain
\begin{align*}
    \left|\beta(u, v) - \int_0^\infty S\left(x, v\right)A(u,x)  dx \right|\hspace{-2.5cm}&\\
    & \leq \left| 2 M \int_0^\infty K(x) \int_{-\infty}^\infty |y|\left\{\int_{-\infty}^\infty\psi(x,y+t)\psi(u,t) dt\right\}^2 dy dx \right|\\
    & \leq \left| 2 M \omega  \int_0^\infty x K(x) \int_{-\infty}^\infty\left\{\int_{-\infty}^\infty\psi(x,y+t)\psi(u,t) dt\right\}^2 dy dx \right|\\
    & \leq \gamma \int_0^\infty x K(x) dx,
\end{align*}
setting again $\gamma = 2M\omega$. 

\subsection{Proof of Proposition \ref{prop:nontrivial}}

    First of all, we have that 
    \begin{equation*}
        A(x,x) = \int_{-\infty}^\infty \Psi(x,\tau)^2d\tau > 0.
    \end{equation*}
    Since A is continuous, we have that for each $x\in\mathbb R$, $A(x,y)>0$ for all $x-g(x) < y < x+g(x)$, for some function $g(x)>0$. Now, assume that there exists $f\in N(T_A)$ such that $f\neq 0$. This means that
    \begin{equation} \label{eq:null1}
        (T_A f)(x) =  \int_0^\infty A(x,y) f(y) dy = 0,\quad\forall x\in \mathbb R^+.
    \end{equation}
    Using the fact that $A(\beta x,\beta y) = \beta A(x,y)$ for all $\beta>0$, we have
    \begin{equation*}
        \int_0^\infty A(x,y) f(y) dy = \int_0^\infty x A\left(1,\frac{y}{x}\right) f(y) dy = \int_0^\infty x^2 A\left(1,u\right) f(xu) du,
    \end{equation*}
    where we made the substitution $y=xu$. Hence (\ref{eq:null1}) becomes
    \begin{equation*}
        \int_0^\infty A(1,u) f(xu) du = 0,\quad\forall x\in \mathbb R^+,
    \end{equation*}
    or equivalently
    \begin{equation*}
        \langle A(1,\,\cdot\,), f(x\,\cdot\,)\rangle = 0\quad\forall x\in \mathbb R^+,
    \end{equation*}
    which implies that $A(1,u)=0$ for all $u\in \mathbb R^+$, a contradiction. Hence $N(T_A)=\{0\}$. 

\subsection{Proof of Theorem \ref{thm:cov_inv}}
    Multiplying both sides of \begin{equation*}
    c(t,\tau) = \int_0^\infty S(u,t)\Psi(u,\tau)du,
    \end{equation*} 
    by $\Psi(u,\tau)$ and integrating, we have
    \begin{equation*}
        \int_{-\infty}^\infty c(t,\tau)\Psi(u,\tau) d\tau = \int_0^\infty S(x,t) A(u,x) dx = \{T_A S(\,\cdot\,,t)\}(u).
    \end{equation*}
    Using the inverse of \begin{equation}
    (T_A f)(x) = \sum_{n=1}^\infty \lambda_n \langle\varphi_n,f\rangle\varphi_n(x),
\end{equation}
we obtain
    \begin{align*}
        S(u, t) &= \sum_{n=1}^\infty \frac{1}{\lambda_n}\varphi_n(u)\int_{-\infty}^\infty c(t,\tau) \langle \Psi(\,\cdot\,,\tau), \varphi_n\rangle d\tau\\
        &= \int_{-\infty}^\infty c(t,\tau) \sum_{n=1}^\infty \frac{1}{\lambda_n}\langle \Psi(\,\cdot\,,\tau), \varphi_n\rangle\varphi_n(u)  d\tau,
    \end{align*}
    as required.  

\subsection{Proof of Example~\ref{ex:white_noise}}
The Haar wavelet coefficients of the Gaussian white noise process $G(t)$ are
\begin{align*}
    d(u,v) &= \int_{-\infty}^\infty G(t)\psi_H(u, t-v) dt\\
    & = u^{-1/2} \int_{-\infty}^\infty G(t)\left\{ \mathbbm{1}_{[0,u/2]}(t-v) - \mathbbm{1}_{[u/2, u]}(t-v)\right\} dt\\
    &= u^{-1/2}\left\{\int_{v}^{v+u/2} G(t) dt -\int_{v+u/2}^{v+u} G(t) dt \right\}\\
    & = u^{-1/2}\left\{G\left((v, v+u/2]\right) - G\left((v+u/2, v+u]\right)\right\}.
\end{align*}
Hence,
\begin{align*}
    \beta(u,v) &= E\{I(u,v)\} = E\left\{|d(u,v)|^2\right\}\\
    & = \frac{1}{u}E \left\{G((v, v+u/2])^2 - 2G((v, v+u/2])G((v+u/2, v+u]) + G((v+u/2, v+u])^2\right\}\\
    & = \frac{1}{u} \left(\frac{u\sigma^2}{2} + \frac{u\sigma^2}{2}\right) = \sigma^2,
\end{align*}
since
\begin{equation*}
    G((v,v+u/2]) \sim\mathcal{N}\left(0,\frac{u\sigma^2} {2}\right)\quad\text{and}\quad G((v+u/2,v+u]) \sim\mathcal{N}\left(0,\frac{u\sigma^2} {2}\right).
\end{equation*}
To obtain the spectrum, we need to solve
\begin{equation*}
    \sigma^2 = \int_0^\infty A_H(u,x) S(x, v)dx.
\end{equation*}
Guessing a solution of the form $S(u,v) = Cx^{-2}$ for some $C\in\mathbb{R}$, we have
\begin{align*}
    \int_0^\infty A_H(u,x) \frac{C}{x^2}dx = C\log(2)
\end{align*}
and therefore 
\begin{equation*}
    S(u,v) = \frac{\sigma^2}{\log(2)u^2},
\end{equation*}
for all $u\in\mathbb{R}^+$ and $v\in\mathbb R$.

\subsection{Proof of Proposition \ref{prop:autocov2}}
We have that 
\begin{align*}
    c_{X_T}(t,\tau) &=  E\{X_T(t)X_T(t+\tau)\}\\
    & = \int_0^{J_\psi(T)}\int_{-\infty}^\infty |w_T(u,v)|^2\psi(u,v-t)\psi(u,v-t-\tau) dv  du\\
    & = \int_0^{J_\psi(T)}\int_{-\infty}^\infty |w_T(u,v+t)|^2\psi(u,v)\psi(u,v-\tau) dv  du,
\end{align*}
where we used the fact that $E\{L(du,dv)L(dx,dy)\} = \delta(u-x)\delta(v-y) du  dv dx dy$ and $\delta$ denotes the Dirac delta function. Assumption (iii) in Definition~\ref{def:clswp2}, implies that
\begin{equation*}
    \left|w_T(u,v+t) - W\left(u,\frac{v+t}{T}\right)\right| \leq \frac{C(u)}{T}.
\end{equation*}
Furthermore, from the Lipschitz continuity of $W$, we also have that
\begin{equation*}
    \left|S\left(u,\frac{v+t}{T}\right) - S\left(u,\frac{t}{T}\right)\right| \leq 2MK(u)\left|\frac{v}{T}\right|
\end{equation*}
from which it follows that 
\begin{equation*}
    \left||w_T(u,v+t)|^2 - S\left(u,\frac{t}{T}\right)\right| \leq \frac{1}{T}\{C(u) + 2MK(u)|v|\}.
\end{equation*}
Using the expression for $c(t,\tau)$ derived in the proof of Proposition~\ref{prop:autocov} and replacing $t$ with $z=t/T$, we have
\begin{align*}
    \left|c_{X_T}(t,\tau) - c(z,\tau)\right| & = \left|\int_0^{J_\psi(T)}\int_{-\infty}^\infty \left\{|w_T(u,v+t)|^2 -S\left(u,\frac{t}{T}\right)\right\}\psi(u,v)\psi(u,v-\tau) dv  du\right|\\
    & \leq \left|\int_0^{J_\psi(T)}\int_{-\infty}^\infty \left||w_T(u,v+t)|^2 - S\left(u,\frac{t}{T}\right)\right|\psi(u,v)\psi(u,v-\tau) dv  du \right|\\
    & \leq \frac{1}{T}\left| \int_0^{J_\psi(T)}\int_{-\infty}^\infty \left\{C(u) + 2MK(u)|v|\right\}\psi(u,v)\psi(u,v-\tau) dv  du \right|\\
    & \leq \frac{1}{T}\left| \int_0^{J_\psi(T)} \left\{C(u) + \gamma u K(u)\right\} \int_{-\infty}^\infty \psi(u,v)\psi(u,v-\tau) dv  du \right|\\
    & = \frac{1}{T}\left| \int_0^{J_\psi(T)} \left\{C(u) + \gamma u K(u)\right\} \Psi(u,\tau)  du \right|\\
    & \leq \frac{1}{T}\int_0^{J_\psi(T)} \left\{C(u) + \gamma u K(u)\right\} du \\
    &= O\left(T^{-1}\right),
\end{align*}
since the integral is finite.

\subsection{Proof of Theorem \ref{thm:expectation2}}
We have that 
\begin{align*}
    \beta(u, z) &= \beta\left(u, \frac{v}{T}\right) = E\left[|d(u,v)|^2\right] = E\left[\left\{\int_{-\infty}^\infty X_T(t) \psi(u,t-v) dt\right\}^2\right]\\
    &=\int_{-\infty}^\infty \int_{-\infty}^\infty E \{X_T(t)X_T(s)\}\psi(u,t-v)\psi(u,s-v) dt ds
\end{align*}
Now,
\begin{equation*}
     E \{X_T(t)X_T(s)\} = \int_0^{J_\psi(T)}\int_{-\infty}^\infty |w_T(x,y)|^2\psi(x,t-y)\psi(x,s-y) dx dy,
\end{equation*}
where we used the fact that $ E\{L(dx,dy)L(d\eta,d\xi)\} = \delta(\eta-x)\delta(\xi-y) dx  dy  d\eta  d\xi$ and $\delta$ denotes the Dirac delta function. Hence,
\begin{align*}
    \beta & (u, z)\\
    &=\int_{-\infty}^\infty \int_{-\infty}^\infty \int_0^{J_\psi(T)}\int_{-\infty}^\infty |w_T(x,y)|^2\psi(x,y-t)\psi(x,y-s)\psi(u,t-v)\psi(u,s-v) dy dx dt ds\\  \label{eq:thm2p2}
    & = \int_{-\infty}^\infty \int_{-\infty}^\infty \int_0^{J_\psi(T)}\int_{-\infty}^\infty |w_T(x,y+v)|^2\psi(x,y-t)\psi(x,y-s)\psi(u,t)\psi(u,s) dy dx dt ds\\
    & =  \int_0^{J_\psi(T)}\int_{-\infty}^\infty |w_T(x,y+v)|^2\left\{\int_{-\infty}^\infty\psi(x,y-t)\psi(u,t) dt\right\}^2 dy dx,
\end{align*}
as in the proof of Theorem~\ref{thm:expectation}. From that proof we also have
\begin{align*}
    \int_0^{J_\psi(T)} S\left(x, z\right)A(u,x)  dx = \int_0^{J_\psi(T)}\int_{-\infty}^\infty S(x,z)\left\{\int_{-\infty}^\infty\psi(x,y-t)\psi(u,t) dt\right\}^2 dy dx,
\end{align*}
where we have replaced $v$ with it's rescaled-time counterpart $z$. Now, using the fact that 
\begin{equation*}
    \left||w_T(x,y+t)|^2 - S\left(x,\frac{t}{T}\right)\right| \leq \frac{1}{T}\left\{C(x) + 2MK(x)|y|\right\},
\end{equation*}
as shown in the proof of Proposition~\ref{prop:autocov2},
we have 
\begin{align*}
    \left|\beta(u, z) - \int_0^{J_\psi(T)} S\left(x, z\right) A(u,x) dx\right|\hspace{-3.5cm}& \\
    & \leq \frac{1}{T}\left|\int_0^{J_\psi(T)}\int_{-\infty}^\infty\left\{C(x) + 2MK(x)|y|\right\}\left\{\int_{-\infty}^\infty\psi(x,y+t)\psi(u,t) dt\right\}^2 dy dx \right|\\
    & \leq \frac{1}{T}\left| \int_0^{J_\psi(T)} \left\{C(x) + \gamma xK(x)\right\} \int_{-\infty}^\infty\left\{\int_{-\infty}^\infty\psi(x,y+t)\psi(u,t) dt\right\}^2 dy dx \right|\\
    &  \leq \frac{1}{T} \int_0^{J_\psi(T)} \left\{C(x) + \gamma x K(x)\right\}   dx\\
    & = O\left(T^{-1}\right).
\end{align*}
For the variance, we have that 
\begin{align*}
    \text{var}\left\{I(u,z)\right\} =  E\left\{I(u,z)^2\right\} -  E\left\{I(u,z)\right\}^2 = 2 E\left\{I(u,z)\right\}^2.
\end{align*}
This follows from 
\begin{align}\label{eq:exp2}
    E\left\{\prod_{i=1}^4 L(du_i, dv_i)\right\} = 3\delta(u_1-u_2)\delta(v_1-v_2)\delta(u_3-u_4)\delta(v_3-v_4)\prod_{i=1}^4  du_i  dv_i,
\end{align}
as there are 3 ways to arrange the Lévy bases so that there are exactly 2 overlapping pairs, and so 
\begin{equation*}
    E\left\{I(u,z)^2\right\} = 3 E\left\{I(u,z)\right\}^2.
\end{equation*}
Therefore
\begin{align*}
    \text{var}\left\{I(u,z)\right\} &= 2\left\{\int_0^{J_\psi(T)} A(u,x) S(x,z)  dx + O\left(T^{-1}\right) \right\}^2
\end{align*}
and so
\begin{align*}
    \left |\text{var}\left\{I(u,z)\right\} - 2\left\{\int_0^{J_\psi(T)} A(u,x) S(x,z)  dx \right\}^2\right| \hspace{-4cm} &\\
    &= 2 O\left(T^{-1}\right)\left[\int_0^{J_\psi(T)} A(u,x) S(x,z)  dx + O\left(T^{-1}\right) \right]\\
    &= O\left(uT^{-1}\right).
\end{align*}
This follows from the fact that
\begin{equation}
    \int_0^{J_\psi(T)} A(u,x) S(x,z)  dx \leq A(u,u) \int_0^{J_\psi(T)} S(x,z)  dx = \alpha  u,
\end{equation}
for some $\alpha\in\mathbb R^+$, as $A(u,u)$ is linear in $u$ and $\int_0^\infty S(x,z)  dx$ is finite by assumption.

\subsection{Proof of Theorem \ref{thm:smoothing}}
The proof is similar to that of~\cite{LSWP}. We have that
\begin{equation*}
    E\left(\tilde d^u_{\ell, m}\right) = E \left\{\int_0^1 I(u,z)\tilde{\psi}_{\ell, m}(z) dz\right\} = \int_0^1 \beta(u,z)\tilde{\psi}_{\ell, m}(z) dz.
\end{equation*}
Now, using Theorem~\ref{thm:expectation2}, 
\begin{align*}
    \left|E\left(\tilde d^u_{\ell, m}\right) - \int_0^1\int_0^\infty S\left(x, z\right) A(u,x) dx\, \tilde{\psi}_{\ell, m}(z) dz\right|\hspace{-4cm}&\\ 
    &\leq \left|\int_0^1 \frac{1}{T} \int_0^\infty \left\{C(u) + \gamma u K(u)\right\} du\,  \tilde{\psi}_{\ell, m}(z) dz \right|\\
    &\leq \frac{1}{T} \int_0^\infty \left\{C(u) + \gamma u K(u)\right\}   du   \int_0^1\left| \tilde{\psi}_{\ell, m}(z) \right| dz\\
    &= \frac{2^{\ell/2}}{T} \int_0^\infty \left\{C(u) + \gamma u K(u)\right\}   du\\
    &= O\left(2^{\ell/2}T^{-1}\right).
\end{align*}
The proof for the variance is very similar to that of~\cite{LSWP}. We substitute the spectrum $f$ into equation (4.7) of Theorem 4.3 of~\cite{vSS96} with the asymptotic limit $\int_0^\infty A(u,x) S(x,z) dx$. Since in our case $u$ is fixed, we replace the two-dimensional wavelets with one-dimensional ones. Noting that $C_h=1$ and $N$ is proportional to $u$, where $C_h$ and $N$ are defined in~\cite{vSS96},~\eqref{eq:smoothing_var} follows from (ii) in the proof of Theorem 4.3 of \cite{vSS96}, where the last term is exactly zero since we are assuming that the Lévy basis is Gaussian. This result is similar to that of~\cite{gao}.

\subsection{Proof of Theorem \ref{thm:smoothing_thresh}}
The proof closely follows the steps in \cite{vSS96} and \cite{gao}. It is easy to show that 
\begin{equation*}
    \tilde{d}^u_{\ell,m} = X^T_T B_T X_T,
\end{equation*}
where the superscript $T$ denotes the transpose operation, $X_T = \left\{X_T(t_0), \ldots, X_T(t_{n(T) - 1})\right\}^T$ and $B_T$ is an $n(T)\times n(T)$ matrix with entries
\begin{equation*}
    \left(B_T\right)_{i,j} = \frac{\tilde\Delta t_i \tilde\Delta t_j}{4 M_v}\sum_{k=0}^{M_v - 1} \psi(u, t_i - z_k T) \psi(u, t_j - z_k T) \tilde\psi_{\ell, m}(z_k).
\end{equation*}
Hence, we can get an equivalent expression for the wavelet coefficients $\tilde{d}^u_{\ell,m}$ to equation~(4.11) of \cite{vSS96} and so Proposition~4.8 of \cite{vSS96} holds for $\tilde{d}^u_{\ell,m}$ as well. Therefore, as in \cite{vSS96},
\begin{equation*}    
\text{pr}\left\{\left|\frac{\tilde{d}^u_{\ell,m} - d^u_{\ell,m}}{\tilde{\sigma}^u_{\ell,m}}\right| \geq \log(T)\right\} = O\left(T^{-\frac{1}{2}}\right),
\end{equation*}
where $d^u_{\ell,m}$ denote the true wavelet coefficients of $\beta(u, \,\cdot\,)$. We note that 
\begin{align*}
    \text{var} \left(\tilde{d}^u_{\ell,m}\right) &= \int_0^1\left\{\int_0^\infty S\left(x, z\right) A(u,x) dx\right\}^2 \tilde{\psi}_{\ell, m}^2(z) dz + O\left(2^\ell u T^{-2}\right)\\
    & = \int_0^1 \beta^2(u,z)\tilde{\psi}^2_{\ell, m}(z) dz + O\left(2^\ell u T^{-2}\right)\\
    &\leq 2T^{-1}\sup_{z\in(0,1)} \beta^2(u,z)
\end{align*}
and therefore 
\begin{equation*}
    \lambda(\ell, m, u, T) = \tilde\sigma^u_{\ell, m}\log\left(T\right)\leq 2^{1/2}\sup_{z\in(0,1)} |\beta(u,z)| \; T^{-\frac{1}{2}}\log(T).
\end{equation*}

Hence, the upper estimate of $\lambda(\ell, m, u, T)$ is equivalent to that of \cite{vSS96}. Since $\beta(u,\,\cdot\,)$ is Lipschitz continuous for each $u\in\mathbb R^+$, it is Holder continuous with $\sigma=1$, using the notation of \cite{vSS96}. The Lipschitz continuity of $\beta(u,\,\cdot\,)$ follows from the fact that it is a linear transformation of $S(u,\,\cdot\,)$, which is itself Lipschitz continuous. Further, we are smoothing each level $u$ separately and so the dimension in our setting is 1. Therefore~\eqref{eq:smoothing_rate} follows from a straightforward application of the results in  Section 4.5 of \cite{vSS96}.

\subsection{Proof of Proposition \ref{prop:haar1}}
We have that 
\begin{align*}
    \Psi_H(\tau) &= \int_{-\infty}^{\infty}\psi_H(x)\psi_H(x-\tau) dx\\
    &= \int_{-\infty}^{\infty}\left\{\mathbbm{1}_{\left[0,\frac{1}{2}\right]} (x) - \mathbbm{1}_{\left[\frac{1}{2}, 1\right]} (x)\right\}\left\{\mathbbm{1}_{\left[0,\frac{1}{2}\right]} (x-\tau) - \mathbbm{1}_{\left[\frac{1}{2}, 1\right]} (x-\tau)\right\}  dx\\
    &= \int_{-\infty}^{\infty}\left\{\mathbbm{1}_{\left[0,\frac{1}{2}\right]} (x) - \mathbbm{1}_{\left[\frac{1}{2}, 1\right]} (x)\right\}\left\{\mathbbm{1}_{\left[\tau,\frac{1}{2}+\tau\right]} (x) - \mathbbm{1}_{\left[\frac{1}{2}+\tau, 1+\tau\right]} (x)\right\}  dx\\
    &= \int_{-\infty}^{\infty}\mathbbm{1}_{\left[0,\frac{1}{2}\right]\cap \left[\tau,\frac{1}{2}+\tau\right]} (x) - \mathbbm{1}_{\left[\frac{1}{2}, 1\right]\cap \left[\tau,\frac{1}{2}+\tau\right]} (x)dx \\ &\hspace{4cm} - \int_{-\infty}^{\infty} \mathbbm{1}_{\left[0,\frac{1}{2}\right]\cap\left[\frac{1}{2}+\tau, 1+\tau\right]} (x) + \mathbbm{1}_{\left[\frac{1}{2}, 1\right]\cap\left[\frac{1}{2}+\tau, 1+\tau\right]} (x)  dx\\
    &= 2\left(\frac{1}{2}-|\tau|\right)\mathbbm{1}_{\left[-\frac{1}{2},\frac{1}{2}\right]}(\tau) - \left(\frac{1}{2}-\left|\tau-\frac{1}{2}\right|\right)\mathbbm{1}_{\left[0,1\right]}(\tau) - \left(\frac{1}{2}-\left|\tau+\frac{1}{2}\right|\right)\mathbbm{1}_{\left[-1,0\right]}(\tau) \\
    & = \left\{\begin{array}{lll}
     1-3|\tau|,    & & 0<|\tau| \leq 1/2,\\
     |\tau| - 1,    &&  1/2 < |\tau| \leq 1,\\
     0, && \text{otherwise.}
    \end{array}\right.
\end{align*}
Using property (iii) from Proposition~\ref{prop:properties1}, we have
\begin{align*}
    \Psi_H(u,\tau) = \left\{\begin{array}{lll}
     1-3|\tau|/u,   & &  0<|\tau| \leq u/2,\\
     |\tau|/u - 1, &  &  u/2 < |\tau| \leq u,\\
     0,& & \text{otherwise,}
    \end{array}\right.
\end{align*}
as required.

\subsection{Proof of Proposition \ref{prop:haar2}}
We have
\begin{align*}
    A_H(u,x) &= \int_{-\infty}^\infty \Psi_H(u, \tau) \Psi_H(x, \tau)  d\tau\\
    & = \int_{-\infty}^\infty \left\{\left(1-3\frac{|\tau|}{u} \right)\mathbbm{1}_{\left[0,\frac{u}{2}\right]}(|\tau|) + \left(\frac{|\tau|}{u} - 1\right)\mathbbm{1}_{\left[\frac{u}{2}, u\right]}(|\tau|) 
    \right\}\\
    &\hspace{3.5cm}
    \left\{\left(1-3\frac{|\tau|}{x} \right)\mathbbm{1}_{\left[0,\frac{x}{2}\right]}(|\tau|) + \left(\frac{|\tau|}{x} - 1\right)\mathbbm{1}_{\left[\frac{x}{2}, x\right]}(|\tau|) 
    \right\} d\tau 
\end{align*}
The result is easily obtained by evaluating each of the cross terms individually.

\subsection{Proof of Proposition \ref{prop:ricker1}}

We have that 
\begin{align*}
    \Psi_R(u,\tau) &= \int_{-\infty}^{\infty}\psi_R(u,x)\psi_R(u,x-\tau) dx\\
    &=\frac{4}{3u\pi^{1/2}} \int_{-\infty}^{\infty}\left(1 - \frac{x^2}{u^2}\right)\left\{1 - \frac{(x-\tau)^2}{u^2}\right\}\exp\left\{-\frac{x^2+(x-\tau)^2}{2u^2}\right\}dx\\
    &=\frac{4}{3u\pi^{1/2}} \int_{-\infty}^{\infty} \left\{1 -\frac{x^2+(x-\tau)^2}{u^2} + \frac{x^2(x-\tau)^2}{u^4}\right\} \exp\left\{-\frac{\left(x-\frac{\tau}{2}\right)^2+\frac{\tau^2}{4}}{u^2}\right\}dx\\
    &=\frac{4e^{-\frac{\tau^2}{4u^2}}}{3u\pi^{1/2}} \int_{-\infty}^{\infty} \left\{1 -2\frac{\left(x-\frac{\tau}{2}\right)^2+\frac{\tau^2}{4}}{u^2}+ \frac{x^2(x-\tau)^2}{u^4}\right\} \exp\left\{-\frac{\left(x-\frac{\tau}{2}\right)^2}{u^2}\right\}dx,
\end{align*}
completing the square in the exponential. Now, using the substitution
\begin{equation*}
    2^{1/2}\left(\frac{x-\tau/2}{u}\right) \mapsto y,
\end{equation*}
we have
\begin{align*}
    \Psi_R(u,\tau) &=\frac{4}{3}e^{-\frac{\tau^2}{4u^2}} \int_{-\infty}^{\infty} \left\{1 - \frac{\tau^2}{2u^2}+ \frac{\tau^4}{16u^4} - \left(1+\frac{\tau^2}{4u^2}\right)y^2 + \frac{1}{4}y^4  \right\}\frac{1}{(2\pi)^{1/2}} e^{-\frac{y^2}{2}}dy\\
    &=\frac{4}{3}e^{-\frac{\tau^2}{4u^2}} \left(1 - \frac{\tau^2}{2u^2}+ \frac{\tau^4}{16u^4} - 1-\frac{\tau^2}{4u^2} + \frac{3}{4} \right)\\
    &=\left(1 + \frac{\tau^4}{12u^4}-\frac{\tau^2}{u^2}\right)e^{-\frac{\tau^2}{4u^2}},
\end{align*}
where we used the fact that we are integrating against the probability density function of a standard normal distribution and so can use its moments. 

\subsection{Proof of Proposition \ref{prop:ricker2}}

We have
\begin{align*}
    A_R(u,x) &= \int_{-\infty}^\infty \Psi_R(u, \tau) \Psi_R(x, \tau)  d\tau\\
    &= \int_{-\infty}^\infty \left(1 + \frac{\tau^4}{12u^4}-\frac{\tau^2}{u^2}\right)\left(1 + \frac{\tau^4}{12x^4}-\frac{\tau^2}{x^2}\right)\exp\left\{-\frac{\tau^2}{4u^2}-\frac{\tau^2}{4x^2}\right\}d\tau\\
    &= \int_{-\infty}^\infty \left\{1 - \left(\frac{1}{x^2}+\frac{1}{u^2}\right)\tau^2 + \left(\frac{1}{12u^4} + \frac{1}{x^2u^2}+\frac{1}{12x^4}\right)\tau^4  \right.\nonumber\\ 
    &\hspace{1.5cm}\left. -\frac{1}{12}\left(\frac{1}{u^2x^4} + \frac{1}{u^4x^2}\right)\tau^6+ \frac{1}{144u^4x^4}\tau^8\right\}\exp\left\{-\frac{1}{4}\left(\frac{1}{u^2} + \frac{1}{x^2}\right)\tau^2\right\}d\tau.
    \end{align*}
    Making the substitution 
    \begin{equation*}
        \frac{1}{2}\left(\frac{1}{u^2}+\frac{1}{x^2}\right)\tau^2\mapsto y^2,
    \end{equation*}
    we have 
    \begin{align*}
        A_R(u,x) &= 2ux\left(\frac{\pi}{u^2+x^2}\right)^{1/2}\int_{-\infty}^\infty\frac{1}{2\pi}e^{-\frac{y^2}{2}} \left\{1-2y^2 + \frac{1}{3}y^4 + \frac{10u^2x^2}{3(u^2+x^2)^2}y^4\right.\\
        &\hspace{6cm}\left.- \frac{2u^2x^2}{3(u^2+x^2)^2}y^6 + \frac{u^4x^4}{9(u^2+x^2)^4} y^8\right\} dy\\
         &= 2ux\left(\frac{\pi}{u^2+x^2}\right)^{1/2}\left\{1-2+1+\frac{10u^2x^2}{(u^2+x^2)^2} - \frac{10u^2x^2}{(u^2+x^2)^2} +\frac{35u^4x^4}{3(u^2+x^2)^4}\right\}\\
         &=\frac{70\pi^{1/2} u^5x^5}{3(u^2+x^2)^{\frac{9}{2}}},
    \end{align*}
    where again we used the fact that we are integrating against the probability density function of a standard normal.

\subsection{Proof of Proposition~\ref{prop:morlet}}
From the definition, we have that
\begin{align*}
    \Psi_M(\tau) & = \int_{-\infty}^{\infty} \psi_M(x) \psi_M(x - \tau) dx\\
    & = \int_{-\infty}^{\infty} \frac{2}{\pi^{1/2}} \cos\left(\eta^{1/2}x\right)\cos\left\{\eta^{1/2}(x-\tau)\right\}e^{-\frac{x^2}{2}} e^{-\frac{(x-\tau)^2}{2}}dx\\  \numberthis\label{eq:m_ac1}
    & =\frac{1}{\pi^{1/2}} \int_{-\infty}^{\infty} \left[\cos\left\{\eta^{1/2}(2x-\tau)\right\} + \cos\left(\eta^{1/2}\tau\right)\right]e^{-\left(x-\frac{\tau}{2}\right)^2 -\frac{\tau^2}{4}}dx\\
    &=\frac{1}{\pi^{1/2}}e^{-\frac{\tau^2}{4}}\left[\int_{-\infty}^{\infty}\cos\left\{\eta^{1/2}(2x-\tau)\right\} e^{-\left(x-\frac{\tau}{2}\right)^2}dx + \cos\left(\eta^{1/2}\tau\right)\int_{-\infty}^{\infty}e^{-\left(x-\frac{\tau}{2}\right)^2}dx\right]\\\numberthis\label{eq:m_ac2}
    &=\frac{1}{\pi^{1/2}}e^{-\frac{\tau^2}{4}}\left[\int_{-\infty}^{\infty}\cos\left\{2\eta^{1/2}y\right\} e^{-y^2}dy + \pi^{1/2}\cos\left(\eta^{1/2}\tau\right)\right]\\\numberthis\label{eq:m_ac3}
    &=\frac{1}{\pi^{1/2}}e^{-\frac{\tau^2}{4}}\left\{\pi^{1/2}e^{-\eta} + \pi^{1/2}\cos\left(\eta^{1/2}\tau\right)\right\}\\
    &=\left\{e^{-\eta} + \cos\left(\eta^{1/2}\tau\right)\right\}e^{-\frac{\tau^2}{4}},
\end{align*}
where in \eqref{eq:m_ac1} we used the standard cosine identity
\begin{equation*}
    \cos(a)\cos(b) = \frac{\cos(a+b) + \cos(a-b)}{2} 
\end{equation*}
and completed the square in the exponential, in \eqref{eq:m_ac2} we made the substitution $x - \tau/2\mapsto y$, and in \eqref{eq:m_ac3} we used the fact that
\begin{equation}\label{eq:cos_exp}
    \int_{-\infty}^{\infty}e^{-ax^2}\cos(bx + c)dx = \left(\frac{\pi}{a}\right)^{1/2}e^{-\frac{b^2}{4a}}\cos(c),
\end{equation}
for some $a$, $b$ and $c$. Hence, by property (iii) of Proposition~\ref{prop:properties1},
\begin{equation*}
    \Psi_M(u,\tau) = \Psi_M\left(\frac{\tau}{u}\right) = \left\{e^{-\eta} + \cos\left(\frac{\eta^{1/2}\tau}{u}\right)\right\}e^{-\frac{\tau^2}{4u^2}}. 
\end{equation*}

\subsection{Proof of Proposition~\ref{prop:morlet2}}
We have that
\begin{align*}
    A_M(u,x) &= \int_{-\infty}^\infty \Psi_M(u, \tau) \Psi_M(x, \tau)  d\tau\\
    & = \int_{-\infty}^\infty \left\{e^{-\eta} + \cos\left(\frac{\eta^{1/2}\tau}{u}\right)\right\}
    \left\{e^{-\eta} + \cos\left(\frac{\eta^{1/2}\tau}{x}\right)\right\}
    e^{-\frac{\tau^2}{4u^2}}e^{-\frac{\tau^2}{4x^2}} d\tau\\
    & = \int_{-\infty}^\infty \left\{e^{-2\eta} + e^{-\eta}\cos\left(\frac{\eta^{1/2}\tau}{u}\right) + e^{-\eta}\cos\left(\frac{\eta^{1/2}\tau}{x}\right)\right. \\
    &\hspace{4.5cm}+ \left.\cos\left(\frac{\eta^{1/2}\tau}{u}\right)\cos\left(\frac{\eta^{1/2}\tau}{x}\right)\right\}
    e^{-\frac{\tau^2}{4}\left(\frac{1}{u^2} + \frac{1}{x^2}\right)} d\tau\\
    & = \int_{-\infty}^\infty \left[e^{-2\eta} + e^{-\eta}\cos\left(\frac{\eta^{1/2}\tau}{u}\right) + e^{-\eta}\cos\left(\frac{\eta^{1/2}\tau}{x}\right)\right. \\
    &\hspace{2.5cm}+ \left.
    \frac{1}{2}\cos\left\{\eta^{1/2}\left(\frac{1}{u} + \frac{1}{x}\right)\tau \right\} + \frac{1}{2}\cos\left\{\eta^{1/2}\left(\frac{1}{u} - \frac{1}{x}\right)\tau \right\}
    \right]
    e^{-\frac{\tau^2}{4}\left(\frac{1}{u^2} + \frac{1}{x^2}\right)} d\tau. 
\end{align*}
The first term in the integral above is simply
\begin{align*}
    e^{-2\eta}\int_{-\infty}^\infty e^{-\frac{\tau^2}{4}\left(\frac{1}{u^2} + \frac{1}{x^2}\right)}d\tau = 2 \pi^{1/2} e^{-2\eta} \left(\frac{1}{u^2} + \frac{1}{x^2}\right)^{-1/2},
\end{align*}
whereas the remaining terms can be evaluated using equation~\eqref{eq:cos_exp} to obtain
\begin{align*}
    A_M(u,x) &= 2 \pi^{1/2} e^{-\eta} \left(\frac{1}{u^2} + \frac{1}{x^2}\right)^{-1/2}\left(e^{-\eta} + e^{-\frac{\eta x^2}{u^2 + x^2}}+ e^{-\frac{\eta u^2}{u^2 + x^2}} + \frac{1}{2}e^{-\frac{2\eta ux}{u^2 + x^2}} + \frac{1}{2}e^{\frac{2\eta ux}{u^2 + x^2}}\right)\\
    &=2 \pi^{1/2} e^{-\eta} \left(\frac{1}{u^2} + \frac{1}{x^2}\right)^{-1/2}\left\{e^{-\eta} + e^{-\frac{\eta x^2}{u^2 + x^2}}+ e^{-\frac{\eta u^2}{u^2 + x^2}}+ \cosh\left(\frac{2\eta ux}{u^2 + x^2}\right)\right\},
\end{align*}
which follows from the identity
\begin{equation*}
    \cosh(a) = \frac{e^{a} + e^{-a}}{2}
\end{equation*}
for some $a$. Setting 
\begin{equation*}
    \lambda = \frac{u^2}{u^2 + x^2}
\end{equation*}
for simplicity, we also have
\begin{align*}
    e^{-\frac{\eta x^2}{u^2 + x^2}}+ e^{-\frac{\eta u^2}{u^2 + x^2}} 
    &= e^{-\eta(1-\lambda)}+e^{-\eta \lambda}\\
    &= e^{-\eta/2}\left\{e^{\eta \left(\lambda - \frac{1}{2}\right)}+ e^{-\eta \left(\lambda - \frac{1}{2}\right)}\right\}\\
    &= 2e^{-\eta/2}\cosh\left\{\eta\left(\lambda-\frac{1}2{}\right)\right\}\\
    & = 2e^{-\eta/2}\cosh\left\{\frac{\eta\left(u^2-x^2\right)}{2\left(u^2-x^2\right)}\right\}.
\end{align*}
Hence, putting everything together,
\begin{align*}
    A_M(u,x) =2 \pi^{1/2} e^{-2\eta} \left(\frac{1}{u^2} + \frac{1}{x^2}\right)^{-1/2}\left[1 + 2e^{\eta/2}\cosh\left\{\frac{\eta\left(u^2-x^2\right)}{2\left(u^2-x^2\right)}\right\} + e^{\eta}\cosh\left(\frac{2\eta ux}{u^2 + x^2}\right)\right].
\end{align*}

\subsection{Proof of Proposition \ref{prop:shannon1}}

To begin with, we have that
\begin{equation*}
    \Psi_S(u, \tau) = \int_{-\infty}^{\infty}\psi_S(u,x)\psi_S(u,x-\tau)dx = \int_{-\infty}^{\infty}\psi_S(u,x)\psi_S(u,\tau-x)dx = \psi_S(u,\,\cdot\,) *\psi_S(u,\,\cdot\,)
\end{equation*}
as Shannon wavelets are symmetric and where $*$ denotes the convolution operator. Hence
\begin{equation*}
    \hat\Psi_S(u, \omega) = \hat\psi_S(u, \omega) \hat\psi_S(u, \omega) = u^{1/2} \hat\psi_S(u, \omega),
\end{equation*}
since the square of an indicator function is simply equal to the indicator function and the cross term is zero since the intervals are non-overlapping. Hence, taking the inverse Fourier transform
\begin{equation*}
    \Psi_S(u, \tau) = u^{1/2} \psi_S(u, \tau),
\end{equation*}
as required.

\subsection{Proof of Proposition \ref{prop:shannon2}}

Using the convolution theorem, we have
\begin{align*}
    A_S(u,x) &= \int_{-\infty}^\infty \Psi_S(u, \tau) \Psi_S(x, \tau)  d\tau\\
    & =(ux)^{1/2}\int_{-\infty}^\infty \psi_S(u, \tau) \psi_S(x, \tau)  d\tau\\
    &=(ux)^{1/2}\int_{-\infty}^\infty \psi_S(u, \tau) \psi_S(x, -\tau)  d\tau\\
    &=(ux)^{1/2} \mathcal{F}^{-1}\left\{\hat\psi_S(u, \,\cdot\,) \hat\psi_S(x, \,\cdot\,)\right\}(0),
\end{align*}
where $\mathcal{F}^{-1}(f)(0)$ denotes the inverse Fourier transform of the function $f$ evaluated at zero, and
\begin{equation*}
        \hat\psi_S(u, \omega) = u^{1/2} \left\{\mathbbm{1}_{\left[-\frac{2\pi}{u}, -\frac{\pi}{u}\right]}(\omega) + \mathbbm{1}_{\left[\frac{\pi}{u}, \frac{2\pi}{u}\right]}(\omega) \right\}.
    \end{equation*}
is the Fourier transform of $\psi$. For $x > 2u$, the above is clearly zero as none of the intervals in the indicator functions overlap. For $u \leq x \leq 2u$, we have 
\begin{equation*}
    \hat\psi_S(u,\omega) \hat\psi_S(x, \omega) = (ux)^{1/2}\left\{\mathbbm{1}_{[-2\pi / x,-\pi/u]}(\omega)+\mathbbm{1}_{[\pi/u,2\pi /x]}(\omega)\right\}
\end{equation*}
The case $x < u$ is symmetrical. Taking the inverse Fourier transform of this, we obtain
\begin{align*}
    \mathcal{F}^{-1}\left\{\hat\psi_S(u, \,\cdot\,) \hat\psi_S(x, \,\cdot\,)\right\}(y) 
    &= \frac{(ux)^{1/2}}{2\pi}\left(\int_{-2\pi/x}^{-\pi/u}e^{i2\pi\omega y} d\omega + \int_{\pi/u}^{2\pi/x}e^{i2\pi\omega y} d\omega\right)\\
    & =  \frac{(ux)^{1/2}}{4\pi^2iy}\left( e^{i4\pi^2y/x} - e^{i2\pi^2y/u} + e^{-i2\pi^2y/u} + e^{-i4\pi^2y/x}  \right)\\
    & = \frac{(ux)^{1/2}}{2\pi^2y}\left\{\sin\left(\frac{4\pi^2 y}{x}\right) - \sin\left(\frac{2\pi^2 y}{x}\right) \right\}
\end{align*}
Now,
\begin{equation*}
     \lim_{y\rightarrow 0}\left[\frac{1}{2\pi^2y}\left\{\sin\left(\frac{4\pi^2 y}{x}\right) - \sin\left(\frac{2\pi^2 y}{x}\right) \right\}\right] = \frac{2}{x} - \frac{1}{u},
\end{equation*}
since
\begin{equation*}
    \lim_{y\rightarrow 0}\frac{\sin (\alpha y)}{y} = \alpha.
\end{equation*}
Hence, 
\begin{equation}
    A_S(u,x) = 2u-x,\quad\text{for }u\leq x\leq 2u,
\end{equation}
and the rest follows by symmetry.

\end{document}